\documentclass[a4paper,11pt,twoside]{article}
\usepackage{pst-fill,pst-grad}
\usepackage{textcomp}
\usepackage[english]{babel}
\usepackage[utf8x]{inputenc}
\usepackage{graphicx}
\usepackage{multicol}
\usepackage{float}
\usepackage{fancyhdr}
\usepackage[matrix,arrow,curve]{xy}
\usepackage{pstricks} 
\usepackage{amsmath,amsfonts,verbatim,afterpage,theorem,euscript,mathrsfs,amssymb}
\usepackage{array}
\usepackage{dsfont}
\usepackage{hyperref}

\numberwithin{equation}{section}
\newtheorem{Theorem}{Theorem}
\newtheorem{Lemma}{Lemma}
\newtheorem{Proposition}{Proposition}
\newtheorem{Corollary}{Corollary}
\newtheorem{Remark}{Remark}
\newtheorem{Definition}{Definition}
\newcommand{\vu}{\vec{u}}
\newcommand{\vw}{\vec{\omega}}
\newcommand{\vv}{\vec{v}}
\newcommand{\vf}{\vec{f}}
\newcommand{\vg}{\vec{g}}
\newcommand{\vV}{\vec{V}}
\newcommand{\vphi}{\vec{\phi}}
\newcommand{\vUU}{\vec{\mathcal{U}}}
\def \p{\mathfrak{p}}
\def \q{\mathfrak{q}}
\newcommand{\vn}{\vec{\nabla}}
\newcommand{\rot}{\vn \wedge}
\def \R{\mathbb{R}^3}
\def \M{\mathcal{M}}
\title{\bf Partial suitable solutions for the micropolar equations and regularity properties}
\usepackage{authblk}
\author{Diego Chamorro\footnote{Laboratoire de Math\'ematiques et Mod\'elisation d'Evry (LaMME), Universit\'e Paris-Saclay, France.}, David Llerena$^*$} 

\begin{document}
\maketitle
\begin{abstract}
The incompressible Micropolar system is given by two coupled equations: the first equation gives the evolution of the velocity field $\vu$ while the second equation gives the evolution of the microrotation field $\vw$. In this article we will consider regularity problems for weak solutions of this system. For this we will introduce the new notion of partial suitable solutions, which imposes a specific behavior for the velocity field $\vu$ only, and under some classical hypotheses over the pressure, we will obtain a h\"olderian gain for both variables $\vu$ and $\vw$. 
\end{abstract}
\section{Introduction}

We study here, under mild assumptions over only one variable, some general regularity properties for weak solutions of the 3D incompressible Micropolar equations. This system is composed of two coupled equations: the first one is based in the incompressible 3D Navier-Stokes problem, which gives the evolution of the velocity field $\vu$ with an internal pressure $p$, while the second one considers the evolution of a microrotation field $\vw$ representing the angular velocity of the rotation of the fluid particles. These equations are given by the following problem:
\begin{equation}\label{MicropolarEquationsIntro}
\begin{cases}
\partial_t \vu  = \Delta \vu  -(\vu \cdot \vn)\vu-\vn p +\frac{1}{2}\rot\vw ,\qquad div(\vu)=0, \\[3mm] 
\partial_t \vw = \Delta \vw +\vn div(\vw)-\vw -(\vu \cdot \vn)\vw+\frac{1}{2}\rot\vu ,\\[3mm] 
\vu (0,x)=\vu_0(x), \;\;\vw(0,x)= \vw_0(x)\quad \mbox{and}\quad  div(\vu _0)=0,\quad x\in \R,
\end{cases}
\end{equation}
here $\vu:[0, +\infty[\times \R \longrightarrow \R$ is the velocity field of the fluid, $p:[0, +\infty[\times \R \longrightarrow \mathbb{R}$ is the pressure and $\vw:[0, +\infty[\times \R \longrightarrow \R$ is the angular velocity. Micropolar equations were first introduced in 1966 by Eringen \cite{Eri66} and they are now used in some particular cases, such as in the study of polymers, blood, muddy fluids,  nematic liquid crystals and bubly liquids. We refer to the book \cite{ Luka99} for other applications of this model. From the mathematical point of view, this system was studied in \cite{Cru19, Nic22,  RemTic21,Yam05} where a variety of results were obtained.  \\

Let us start with some remarks about the equations (\ref{MicropolarEquationsIntro}). First note that when the microrotation field $\vw$ is null, we recover the usual 3D incompressible Navier-Stokes equations for irrotational fluids (\emph{i.e.} $\rot\vu=0$) which were studied for instance in \cite{Bat00} or in \cite{Gal11}. Next we observe that the angular velocity $\vw$ is \emph{not} a divergence free vector field and this makes the study of the properties of $\vw$ slightly more delicate to handle. Finally, it is important to observe that the evolution equation for $\vw$ is essentially linear and that there is a relatively mild coupling between the variables $\vu$ and $\vw$: in this article we will exploit this particular point to deduce our main regularity results for the system (\ref{MicropolarEquationsIntro}). \\

Existence of global weak solutions for this system were obtained in 1977 \cite{Gal97} and from now on we will always assume that $(\vu, \vw)\in L^\infty([0, +\infty[, L^2(\R))\cap L^2([0,+\infty[, \dot{H}^1(\R))$ is a weak solution of (\ref{MicropolarEquationsIntro}). Note that information over the pressure $p$ can be easily obtained from $\vu$: indeed, by applying the divergence operator in the first equation of (\ref{MicropolarEquationsIntro}) we obtain, since $div(\vu)=0$ and $div(\rot \vw)=0$, the usual equation for the pressure:
\begin{equation}\label{Eq_Pression}
\Delta p=-div((\vu\cdot \vn)\vu).
\end{equation}
As said before, in this article we are interested in studying regularity issues for the micropolar system. In the realm of fluid dynamics equations (and in particular for the Navier-Stokes equations) this topic is a challenging and often open problem which can be solved under some different sets of hypotheses such as the Serrin criterion (see \cite{ Serr62}, \cite{OLeary}), the Prodi-Serrin criterion (see \cite{Pro59}, \cite{Serr63}) or in the setting of the Caffarelli-Kohn-Nirenberg theory (see \cite{CKN}, \cite{Kukavica}).\\

Concerning the micropolar system (\ref{MicropolarEquationsIntro}), some recent results were obtained in  \cite{ChLl21, ChLl22,LoMelo} where almost all of the previous theories cited above were applied to obtain a regularity gain over the variables $\vu$ and $\vw$. Let us remark that in most of these references the additional information is asked for \emph{both} variables $\vu$ and $\vw$. However, as it was pointed out in \cite{ Cru22,Loa16,RagWU20} and \cite{Yua08} it is possible to make a separated study of each one of these variables. \\

In this article we are going one step further and in our main result (Theorem \ref{Teo_HolderRegularity} below) we will show that just some additional information over the velocity field $\vu$ is needed in order to deduce a gain of regularity for \emph{both} variables $\vu$ and $\vw$. In this sense, when studying regularity issues for the problem (\ref{MicropolarEquationsIntro}), we will say that the velocity field $\vu$ ``dominates'' the angular velocity field $\vw$.\\

To obtain a gain of regularity, we will work over small neighborhoods and for a point $(t,x)\in ]0,+\infty[\times \R$ let us consider the parabolic ball 
\begin{equation}\label{Def_BolaParabolica}
Q_R(t, x)= ]t-R^2, t+R^2[\times B_R(x),
\end{equation}
for some radius $0<R<1$ such that $t-R^2>0$. When the context is clear we will write $Q_R$ instead of $Q_R(t, x)$.\\ 

We introduce now the following concept:
\begin{Definition}[Partial suitable solutions]\label{Def_Partial_Suitable_Sol}
 Consider  $\vu,\vw \in L_t^{\infty}L_x^2(Q_R)\cap L_t^2\dot H_x^1(Q_R)$ two vector fields that satisfy the equation (\ref{MicropolarEquationsIntro}) in the weak sense over the set $Q_R$. Assume moreover that we have the following local information over the pressure: $p\in L^{\frac{3}{2}}_{t,x}(Q_R)$. We will say that $(\vu,p,\vw)$ is a \emph{partial suitable solution} for the micropolar equations (\ref{MicropolarEquationsIntro}) if the distribution $\mu$ given by the expression
\begin{eqnarray}\label{Def_Mu}
\mu=-\partial_t|\vu|^2+\Delta|\vu|^2-2|\vn\otimes\vu|^2-div\left((|\vu|^2+2p)\vu\right)+(\rot \vw)\cdot\vu,
\end{eqnarray}
is a non-negative locally finite measure on $Q_R$.
\end{Definition}
First note that in this local setting each term of the above expression is well defined. Remark also that the notion introduced above is only related to the evolution of the velocity field $\vu$ and that the action of the variable $\vw$ can be seen here as an external force. In a previous work \cite{ChLl21} we considered a non-negative measure involving the evolution of \emph{both} variables $\vu$ and $\vw$, but as we only consider here the equation related to the variable $\vu$ (and not the equation of $\vw$), this weaker notion of \emph{partial suitable solutions} is needed.\\

With all these notions above, we can now state our main result:
\begin{Theorem}\label{Teo_HolderRegularity}
Consider a parabolic ball $Q_R$ given by (\ref{Def_BolaParabolica}). Let $(\vu,p, \vw)$ be a partial suitable solution (in the sense of Definition \ref{Def_Partial_Suitable_Sol}) for the micropolar system (\ref{MicropolarEquationsIntro}) over $Q_R$. There exists a small constant  $0<\epsilon^*\ll1$  such that if for some point $(t_0,x_0)\in Q_R$ we have
\begin{equation}\label{Condition_PetitEpsilon}
\limsup_{r\to 0}\frac{1}{r}\int _{]t_0-r^2, t_0+r^2[\times B(x_0,r)}|\vn \otimes \vu|^2dx ds<\epsilon^*,
\end{equation}
then, the solution $(\vu,\vw)$ is H\"older continuous in time and space for some exponent $0<\alpha <\frac{1}{24}$ in a small neighborhood of $(t_0,x_0)$.
\end{Theorem} 
As it was mentioned before, observe that we only impose conditions on the variable $\vu$ (namely, the partial suitability condition given in (\ref{Def_Mu}) and a good behavior for the gradient of the velocity field given in (\ref{Condition_PetitEpsilon})) and no particular hypotheses are asked for the variable $\vw$. However, and despite of this fact, we will see here that we can deduce a h\"olderian gain of regularity for \emph{both} variables. Of course, the strategy of the proof of Theorem \ref{Teo_HolderRegularity} is  adapted to this setting: indeed, we will first perform a detailed study for the variable $\vu$ using the first equation of (\ref{MicropolarEquationsIntro}), next we will deduce some controls for the variable $\vw$ by studying the second equation of (\ref{MicropolarEquationsIntro}) and only then, once we have gathered enough information, we will obtain the wished gain of regularity for both variables by studying the evolution of the whole system (\ref{MicropolarEquationsIntro}). Finally, let us remark that the interval $0<\alpha <\frac{1}{24}$ for the index of h\"olderian regularity $\alpha$ given above is mainly technical and we do not claim any optimality on it.\\

We can give now the plan of the article: in Section \ref{Sec_Definitions} we present the main tools used in this article and in Section \ref{Sec_Propriedades_Vu1} we study the evolution of the variable $\vu$ to obtain some information on this variable. However, this information will not be enough and in Section \ref{Sec_Propriedades_Vu2} we will perform a more detailed analysis of the properties of the variables $\vu$. Then, in Section \ref{Sec_Propiedades_Vw} we will deduce from the previous sections some properties for the variable $\vw$. Finally, in Section \ref{Sec_Final}, we will gather all these results to give a proof of Theorem \ref{Teo_HolderRegularity}.

\section{Definitions and Useful results}\label{Sec_Definitions}
Before going any further, let us be more explicit about the H\"older regularity stated in Theorem \ref{Teo_HolderRegularity} above. Indeed, we will consider the homogeneous space $(\mathbb{R}\times \R, d, \mu)$ where $d$ is the parabolic distance given by 
$d\big((t,x), (s,y)\big)=|t-s|^{\frac{1}{2}}+|x-y|$ and where $\mu$ is the usual Lebesgue measure $d\mu=dxdt$. We then define the homogeneous (parabolic) H\"older spaces $\dot{\mathcal{C}}^\alpha(\mathbb{R}\times \R, \R)$ with $0<\alpha<1$ by the usual condition:
\begin{equation*}
\|\vphi\|_{\dot{\mathcal{C}}^\alpha}=\underset{(t,x)\neq (s,y)}{\sup}\frac{|\vphi(t,x)-\vphi(s,y)|}{\left(|t-s|^{\frac{1}{2}}+|x-y|\right)^\alpha}<+\infty,
\end{equation*}
and it is with respect to this functional space that we will obtain the regularity gain announced.\\

Let us now say few words about Morrey spaces: although completely absent in the statement of Theorem \ref{Teo_HolderRegularity}, they are a powerful tool when studying problems related to regularity in PDEs. This fact was particularly underlined in \cite{Robinson} and in \cite{PGLR1} for the Navier-Stokes problem since they provide a very natural framework as we shall see later on (see the key Lemma \ref{Lemma_parabolicHolder} below) and in this article we will use them in a systematic manner. Thus, for $1< p\leq q<+\infty$, the (parabolic) Morrey spaces $\mathcal{M}_{t,x}^{p,q}(\mathbb{R}\times \R)$ are defined as the set of measurable functions $\vphi:\mathbb{R}\times\R\longrightarrow \R$ that belong to the space $(L^p_{t,x})_{loc}$ such that $\|\vphi\|_{M_{t,x}^{p,q}}<+\infty$ where
\begin{equation}\label{DefMorreyparabolico}
\|\vphi\|_{\mathcal{M}_{t,x}^{p,q}}=\underset{x_{0}\in \R, t_{0}\in \mathbb{R}, r>0}{\sup}\left(\frac{1}{r^{5(1-\frac{p}{q})}}\int_{|t-t_{0}|<r^{2}}\int_{B(x_{0},r)}|\vphi(t,x)|^{p}dxdt\right)^{\frac{1}{p}}.
\end{equation}
We present now some well-known facts:
\begin{Lemma}[H\"older inequalities]\label{lemma_Product}
\begin{itemize}
\item[]
\item[1)]If $\vf, \vg:\mathbb{R} \times \R\longrightarrow \R$ are two functions such that $\vf\in \M_{t,x}^{p,q} (\mathbb{R} \times \R)$ and $\vg\in L^{\infty}_{t,x} (\mathbb{R} \times \R)$, then for all $1\leq p\leq q<+\infty$ we have
$\|\vf\cdot\vg\|_{\M_{t,x}^{p,q}}\leq  C\|\vf\|_{\M_{t,x}^{p, q}} \|\vg\|_{L^{\infty}_{t,x}}$.
\item[2)] Let $1\leq p_0 \leq q_0 <+\infty$, $1\leq p_1\leq q_1<+\infty$ and $1\leq p_2\leq q_2<+\infty$. If $\tfrac{1}{p_1}+\tfrac{1}{p_2}= \frac{1}{p_0}$ and $\tfrac{1}{q_1}+\tfrac{1}{q_2}=\tfrac{1}{q_0}$, then for two measurable functions $\vf, \vg:\mathbb{R} \times \R\longrightarrow \R$ such that $\vf\in \mathcal{M}^{p_1, q_1}_{t,x}(\mathbb{R} \times \R)$ and $\vg\in \mathcal{M}^{p_2, q_2}_{t,x}(\mathbb{R} \times \R)$, we have the following H\"older inequality in Morrey spaces 
$$\|\vf\cdot \vg\|_{\mathcal{M}^{p_0, q_0}_{t,x}}\leq \|\vf\|_{\mathcal{M}^{p_1, q_1}_{t,x}}\|\vg\|_{\mathcal{M}^{p_2, q_2}_{t,x}}.$$
\end{itemize}
\end{Lemma} 
\begin{Lemma}[Localization]\label{lemma_locindi}
Let $\Omega$ be a bounded set of $\mathbb{R} \times \R$. If we have $1\leq p_0 \leq q_0$, $1\leq p_1\leq q_1$ with the condition $q_0 \leq q_1<+\infty$ and if the function $\vf:\mathbb{R} \times \R\longrightarrow \R$  belongs to the space $\M_{t,x}^{p_1,q_1}(\mathbb{R} \times \R)$ then we have the following localization property 
$$\|\mathds{1}_{\Omega}\vf\|_{\M_{t,x}^{p_0, q_0}} \leq C\|\mathds{1}_{\Omega}\vf\|_{\M_{t,x}^{p_1,q_1}}\leq C\|\vf\|_{\M_{t,x}^{p_1,q_1}}.$$
\end{Lemma} 

In our work, the notion of parabolic Riesz potential (and its properties) will be crucial and for some index $0<\mathfrak{a}<5$ we define the parabolic Riesz potential $\mathcal{L}_{\mathfrak{a}}$ of a locally integrable function $\vec f: \mathbb{R}\times\mathbb{R}^3\to \mathbb{R}^3$ by
\begin{equation}\label{Def_ParabolicRieszPotential}
\mathcal{L}_{\mathfrak{a}}(\vf)(t,x)=\int_{\mathbb{R}}\int_{\mathbb{R}^3}
\frac{1}{(|t-s|^{\frac{1}{2}}+|x-y|)^{5-\mathfrak{a}}}\vec{f}(s,y)dyds.
\end{equation}
Then, we have the following property in Morrey spaces
\begin{Lemma}[Adams-Hedberg inequality]\label{Lemma_Hed}
If $0<\mathfrak{a}<\frac{5}{q}$, $1<p\le q<+\infty$ and $\vf\in \M_{t,x}^{p,q}(\mathbb{R} \times \R)$,
then for $\lambda=1-\frac{\mathfrak{a}q}{5}$ we have the following boundedness property in Morrey spaces:
\begin{equation*}
\|\mathcal{L}_{\mathfrak{a}}(\vf)\|_{\M_{t,x}^{\frac{p}{\lambda},\frac{q}{\lambda}}}\le C\|\vf\|_{\M_{t,x}^{p,q}}. 
\end{equation*}
\end{Lemma}
The three lemmas above constitute our main tools in Morrey spaces. For a more detailed study of these functional spaces we refer to the books \cite{PGLR1} and \cite{Adams}.
\section{A (first) partial gain of information for the variable $\vu$}\label{Sec_Propriedades_Vu1}
In this section we will only focus our study in the variable $\vu$ and its equation: 
\begin{equation*}
\partial_t \vu  = \Delta \vu  -(\vu \cdot \vn)\vu-\vn p +\frac{1}{2}\rot\vw ,\qquad div(\vu)=0. 
\end{equation*}
Here, the variable $\vw$ can be seen as an external force for which we have the information $\vw\in L^\infty_tL^2_x\cap L^2_t\dot{H}^1_x$. Note that at this stage, we will not obtain a gain of regularity for the variable $\vu$, instead, using the hypotheses given in Theorem \ref{Teo_HolderRegularity} above,  we will obtain a gain of \emph{integrability} for $\vu$ (stated, as we shall see, in terms of Morrey spaces). In this sense our first result is the following:
\begin{Proposition}\label{Propo_FirstMorreySpace}
Under the hypotheses of Theorem \ref{Teo_HolderRegularity} consider $(\vu,p,\vw)$ a partial suitable solution for the micropolar equations (\ref{MicropolarEquationsIntro}) over the set $Q_{R}$ given in (\ref{Def_BolaParabolica}). Then there exists a radius $0<R_1<\frac{R}{2}$ and an index $\tau_0>0$ with $\frac{5}{1-\alpha}<\tau_0<\frac{20}{3}$ such that we have the following local Morrey information:
\begin{equation}\label{Conclusion_FirstMorreySpace}
\mathds{1}_{Q_{R_1}(t_0,x_0)} \vu\in \mathcal{M}_{t,x}^{3,\tau_0}(\mathbb{R}\times \R),
\end{equation}
where the point $(t_0,x_0)\in Q_R$ is given by the hypothesis (\ref{Condition_PetitEpsilon}).
\end{Proposition}
{\bf Proof of the Proposition \ref{Propo_FirstMorreySpace}.} The proof of this result is rather technical and our starting point is given by the notion of partial suitable solution: indeed, from the Definition \ref{Def_Partial_Suitable_Sol} and exploiting the positivity of the quantity given in (\ref{Def_Mu}) we easily deduce the following \emph{partial} local energy inequality: for all $\phi\in \mathcal{D}_{t,x}(Q_R)$ (for which we have $\phi(0,x)=0$) we obtain
\begin{eqnarray}
\int_{\R} |\vu|^2\phi dx+2 \int_{\mathbb{R}} \int _{\R}|\vn \otimes \vu|^2\phi\, dx ds
&\leq& \int_{\R}(\partial_t\phi +\Delta \phi)|\vu|^2dxds+2\int_{\mathbb{R}}\int_{\R}  p (\vu\cdot \vn \phi )  dx ds\label{Ineq_Energie1}\\
&&+\int_{\mathbb{R}}\int _{\R} |\vu|^2(\vu \cdot \vn)\phi dx ds+\int_{\mathbb{R}}\int_{\R}(\rot \vw)\cdot (\phi\vu)dxds.\notag
\end{eqnarray}
Although this estimate is fundamental, it is necessary to fix a convenient test function $\phi$ which will allows us to perform some computations. A particular good choice has been given by Scheffer in \cite{Scheffer}:
\begin{Lemma}\label{Lemma_testfonc}
Let $0<r<\frac{\rho}{2}<R<1$. Let $\phi \in \mathcal{C}_0^{\infty}
(\mathbb{R}\times \R)$  be a function such that
\begin{equation*}
\phi(s,y)=r^2\omega\left(\frac{s-t}{\rho^2},\frac{y-x}{\rho}\right)
\theta\left(\frac{s-t}{r^2}\right)\mathfrak{g}_{(4r^2+t-s)}(x-y),
\end{equation*}
where $\omega \in \mathcal{C}_0^{\infty}(\mathbb{R}\times \R)$ is positive function whose support is in $Q_1(0,0)$ and equal to 1 in $Q_{\frac{1}{2}}(0,0)$. In addition $\theta$ is a smooth function non negative such that $\theta=1$ over $]-\infty,1[$ and $\theta=0$ over $]2,+\infty[$ and $\mathfrak{g}_t(\cdot)$ is the usual heat kernel. Then, we have the following points.
\begin{itemize}
\item[1)]the function $\phi$ is a bounded non-negative function, and its support is contained in the parabolic ball $Q_\rho$, and for all $(s,y)\in Q_r(t,x)$ we have the lower bound $\phi(s,y)\ge \frac{C}{r}$,
\item[2)] for all $(s,y)\in Q_\rho (t,x)$ with $0<s<t+r^2$ we have $\phi(s,y)\le \frac{C}{r}$,
\item[3)] for all $(s,y)\in Q_\rho(t,x) $ with $0<s<t+r^2$ we have $\vn \phi(s,y)\le \frac{C}{r^2}$,
\item [4)] moreover, for all $(s,y)\in Q_\rho(t,x) $ with $0<s<t+r^2$ we have $|(\partial_s+\Delta)\phi(s,y)|\le C\frac{r^2}{\rho^5}$.
\end{itemize}
\end{Lemma}
A detailed proof of this lemma can be found for instance in \cite{ChHe21} or in \cite{PGLR1}.\\ 

The strategy is thus the following: by a convenient use of the estimate (\ref{Ineq_Energie1}) and by the properties of the function $\phi$ given in the previous lemma, we will obtain -by controlling the information over small balls by the information over bigger balls- the wished Morrey information stated in Proposition \ref{Propo_FirstMorreySpace}.\\ 

To do so, it will be useful to introduce the following quantities:  for a point $(t,x)\in \mathbb{R}\times \R$ and for a real parameter $r>0$ we write
\begin{equation}\label{Def_Invariants}
\begin{split}
\mathcal{A}_r(t,x)&=\sup_{t-r^2<s<t+r^2}\frac{1}{r}\int_{B(x,r)}|\vu(s,y)|^2dy, 
\qquad\qquad \alpha_r(t,x)=\frac{1}{r} \int _{Q_r(t,x)}|\vn \otimes \vu(s,y)|^2dyds,\\
\lambda_r(t,x)&=\frac{1}{r^2} \int _{Q_r(t,x)}|\vu(s,y)|^3dyds, 
\qquad\qquad \qquad\mathcal{P}_r(t,x)=\frac{1}{r^2} \int _{Q_r(t,x)}|p (s,y)|^{\frac{3}{2}}dyds,
\end{split}
\end{equation}
and when the context is clear we will simply write $\mathcal{A}_r=\mathcal{A}_r(t,x)$. Note that the previous quantities correspond to the information $L^\infty_tL^2_x$, $L^2_t\dot{H}^1_x$, $L^3_{t,x}$ and $L^{\frac32}_{t,x}$. Note also that for $0<r<1$, we have the relationship between $\lambda_r$, $\mathcal{A}_r$ and $\alpha_r$
\begin{equation}\label{Estimation_normeLambda3}
\lambda_r^{\frac{1}{3}}\le C(\mathcal{A}_r+\alpha_r)^{\frac{1}{2}}.
\end{equation}
Indeed, using the definition of $ \lambda_r$ given in (\ref{Def_Invariants}) above and by H\"older inequality we have
$$\lambda_r^{\frac{1}{3}}=\frac{1}{r^{\frac{2}{3}}}\|\vu\|_{L_{t,x}^{3}(Q_r)}\le \frac{C}{r^\frac{1}{2}}\|\vu\|_{L_{t,x}^{\frac{10}{3}}(Q_r)}.$$
Since by interpolation we have $\|\vu(t,\cdot)\|_{L^{\frac{10}{3}}(B_r)}\leq \|\vu(t,\cdot)\|_{L^{2}(B_r)}^{\frac{2}{5}}\|\vu(t,\cdot)\|_{L^{6}(B_r)}^{\frac{3}{5}}$, we can easily deduce that $\|\vu\|_{L_{t,x}^{\frac{10}{3}} (Q_r)} \leq \|\vu\|_{L_t^{\infty}L_x^{2} (Q_r)}^{\frac25}\|\vu\|_{L_t^{2}L_x^{6} (Q_r)}^{\frac{3}{5}}$. Now, for the $L_t^2L_x^6$ norm of $\vu$, we use the classical Gagliardo-Nirenberg inequality (see \cite{Brezis}) to obtain $\|\vu\|_{L_t^{2}L_x^{6} (Q_r)} \leq C\big(\|\vn \otimes\vu\|_{L_t^{2}L_x^{2} (Q_r)} +\|\vu\|_{L_t^{\infty}L_x^{2} (Q_r)}\big)$ and using Young's inequalities we have
\begin{eqnarray*}
\|\vu\|_{L_{t,x}^{\frac{10}{3}} (Q_r)} &\leq& C \|\vu\|_{L_t^{\infty}L_x^{2} (Q_r)}^{\frac 25}\big(\|\vn \otimes \vu\|_{L_t^{2}L_x^{2} (Q_r)}^{\frac 35}+\|\vu\|_{L_t^{\infty}L_x^{2} (Q_r)}^{\frac35} \big) \leq C\big(\|\vu\|_{L_t^{\infty}L_x^{2} (Q_r)}+\|\vn\otimes \vu\|_{L_t^{2}L_x^{2} (Q_r)}\big).
\end{eqnarray*}
Noting that $\|\vu\|_{L_t^{\infty}L_x^{2} (Q_r)}=r^{\frac12}\mathcal{A}_r^{\frac12}$ and $\|\vn\otimes \vu\|_{L^2_tL_x^{2} (Q_r)}=r^{\frac12}\alpha_r^{\frac12}$, we finally obtain (\ref{Estimation_normeLambda3}).\\

We establish now a first relationship between the quantities given in (\ref{Def_Invariants}) that will be helpful to deduce by an iteration procedure the wished Morrey control.
\begin{Lemma}\label{lema_FirstEstimate}
Under the hypotheses of Theorem \ref{Teo_HolderRegularity} and with the notations given in (\ref{Def_Invariants}) we have for any radius $0<r<\frac{\rho}{2}<1$ the inequality
\begin{equation}\label{FirstEstimate}
\mathcal{A}_r+\alpha_r \le C\frac{r^2}{\rho ^2} \mathcal{A}_\rho
+\frac{\rho ^2}{r^2}\alpha_\rho^{\frac{1}{2}}\mathcal{A}_\rho+C\frac{\rho ^2}{r^2}
\mathcal{P}_\rho^{\frac{2}{3}}(\mathcal{A}_\rho+\alpha_\rho)^{\frac{1}{2}}+C\frac{\rho^{\frac{3}{2}}}{r} 
\|\rot \vw\|_{L_{t,x}^2(Q_\rho)}\alpha_{\rho}^{\frac{1}{2}}.
\end{equation}
\end{Lemma}
\textbf{Proof.}  With the support properties of the function $\phi$ stated in the Lemma \ref{Lemma_testfonc} and using the notations (\ref{Def_Invariants}) above we can rewrite the left hand side of the inequality (\ref{Ineq_Energie1}) in the following manner:
\begin{eqnarray}
\mathcal{A}_r+\alpha_r &\leq &\underbrace{\int_{\R}(\partial_t\phi +\Delta \phi)|\vu|^2dxds}_{(1)}+2\underbrace{\int_{\mathbb{R}}\int_{\R}  p (\vu\cdot \vn \phi )dxds}_{(2)}+\underbrace{\int_{\mathbb{R}}\int _{\R} |\vu|^2(\vu \cdot \vn)\phi dx ds}_{(3)}\notag\\
&&+\underbrace{\int_{\mathbb{R}}\int_{\R}(\rot \vw)\cdot (\phi\vu)dxds}_{(4)}. \label{Ineq_Energie2}
\end{eqnarray}
The terms of the right hand side above will be studied separately. Indeed, 
\begin{itemize}
\item[$\bullet$] For the quantity (1) in (\ref{Ineq_Energie2}), using the properties of the function $\phi$ given in Lemma \ref{Lemma_testfonc} and by the definition of the quantity $\mathcal{A}_\rho$ given in (\ref{Def_Invariants}) we have
$$\int_{\R}(\partial_t\phi +\Delta \phi)|\vu|^2dxds\leq C \frac{r^2}{\rho^5}\int_{Q_\rho}|\vu|^2dxds=C\frac{r^2}{\rho^5}\int_{t-\rho^2}^{t+\rho^2}\int_{B_\rho}|\vu|^2dxds\le C\frac{r^2}{\rho ^2} \mathcal{A}_\rho.$$
\item[$\bullet$] For the term (2) in (\ref{Ineq_Energie2}), by the properties of the function $\phi$ given in Lemma \ref{Lemma_testfonc} and by the H\"older inequality, we obtain
$$\int_{\mathbb{R}}\int_{\R}  p (\vu\cdot \vn \phi )dxds\leq \frac{C}{r^2}\int_{t-\rho^2}^{t+\rho^2}\int_{B_\rho} | p | |\vu|dxds\le \frac{C}{r^2}\|p\|_{L_{t,x}^{\frac{3}{2}}(Q_\rho)}\|\vu\|_{L_{t,x}^{3}(Q_\rho)},$$
noting that by (\ref{Def_Invariants}) we have $\|p\|_{L_{t,x}^{\frac{3}{2}}(Q_\rho)}=\rho^{\frac{4}{3}} \mathcal{P}_\rho^{\frac{2}{3}}$ and $\|\vu\|_{L_{t,x}^{3}(Q_\rho)}=\rho^{\frac{2}{3}}\lambda_\rho^{\frac{1}{3}}$, we can thus write
$$\int_{\mathbb{R}}\int_{\R}  p (\vu\cdot \vn \phi )dxds\leq  \frac{C}{r^2}\left(\rho^{\frac{4}{3}} \mathcal{P}_\rho^{\frac{2}{3}}\right)\left(\rho^{\frac{2}{3}}\lambda_\rho^{\frac{1}{3}}\right)\le C\frac{\rho ^2}{r^2} \mathcal{P}_\rho^{\frac{2}{3}}( \mathcal{A}_\rho+\alpha_\rho)^{\frac{1}{2}},$$
where in the last estimate we used the control (\ref{Estimation_normeLambda3}).
\item For the term (3) in (\ref{Ineq_Energie2}), let us first define the average $\displaystyle{(|\vu|^2)_\rho=\frac{1}{|B(x,\rho)|}\int_{B(x,\rho)}|\vu(t,y)|^2dy}$ and since $\vu$  is divergence free we have $\displaystyle{\int _{B_\rho} (|\vu|^2)_\rho(\vu \cdot \vn)\phi dx =0}$. Then, we can write by the properties of the function $\phi$ given in Lemma \ref{Lemma_testfonc} and by the H\"older inequality:
\begin{eqnarray*}
\int_{\mathbb{R}}\int_{\R} |\vu|^2(\vu \cdot \vn)\phi dxds&=&\int_{Q_\rho} [|\vu|^2-(|\vu|^2)_\rho](\vu\cdot 
\vn)\phi dx ds\leq \frac{C}{r^2}\int_{t-\rho^2}^{t+\rho^2}\int_{B_\rho} \left||\vu|^2
-(|\vu|^2)_\rho|\right|\vu|dx ds\\ 
&\le &\frac{C}{r^2}\int_{t-\rho^2}^{t+\rho^2}\||\vu|^2-(|\vu|^2)_\rho\|
_{L^{\frac{3}{2}}(B_\rho)}\|\vu(s,\cdot)\|_{L^{3}(B_\rho)}ds.
\end{eqnarray*}
Now, Poincare's inequality implies
\begin{align*}
&\le \frac{C}{r^2}\int_{t-\rho^2}^{t+\rho^2}\|\vn (|\vu(s,\cdot)|^2)\|_{L^{1}(B_\rho)}\|\vu(s,\cdot)\|_{L^{3}(B_\rho)}ds\\
&\le \frac{C}{r^2}\int_{t-\rho^2}^{t+\rho^2} \|\vu(s,\cdot)\|_{L^2(B_\rho)}\|\vn\otimes\vu(s,\cdot)\|_{L^2(B_\rho)}\|\vu(s,\cdot)\|_{L^{3}(B_\rho)}ds  \\
&\le \frac{C}{r^2} \|\vu\|_{L_t^6L_x^2(Q_\rho)}\|\vn \otimes \vu\|_{L_{t,x}
^2(Q_\rho)} \|\vu\|_{L_{t,x}^3(Q_\rho)},
\end{align*}
where in the last inequality we used the H\"older inequality in the time variable. We observe now that by the notations given in (\ref{Def_Invariants}) we can write 
$$\|\vu\|_{L_t^6L_x^2(Q_\rho)}\le C \rho ^{\frac{1}{3}}\|\vu\|_{L_t^{\infty}L_x^2(Q_\rho)}\le C\rho ^{\frac{5}{6}}\mathcal{A}_\rho^{\frac{1}{2}}, \quad \|\vn \otimes \vu\|_{L_{t,x}
^2(Q_\rho)}=\rho^{\frac{1}{2}}\alpha_\rho^{\frac{1}{2}}, \quad  \| \vu\|_{L_{t,x}
^3(Q_\rho)}=\rho^{\frac{2}{3}}\lambda_\rho^{\frac{1}{3}},$$
and we obtain, by (\ref{Estimation_normeLambda3}): 
$$\int_{\mathbb{R}}\int_{\R} |\vu|^2(\vu \cdot \vn)\phi dxds\leq C \frac{\rho^2}{r^2}\mathcal{A}_\rho^{\frac{1}{2}}\alpha_\rho^{\frac{1}{2}}\lambda_\rho^{\frac{1}{3}}\leq  C \frac{\rho^2}{r^2}\mathcal{A}_\rho^{\frac{1}{2}}\alpha_\rho^{\frac{1}{2}}(\mathcal{A}_\rho+\alpha_\rho)^{\frac{1}{2}}\leq C\frac{\rho^2}{r^2} \alpha_\rho^{\frac{1}{2}}
(\mathcal{A}_\rho+\alpha_\rho).$$
\item Finally, for the term (4) in (\ref{Ineq_Energie2}), by the H\"older inequality and by the properties of the function $\phi$ given in Lemma \ref{Lemma_testfonc} we write
\begin{eqnarray*}
\int_{\mathbb{R}}\int_{\R}(\rot \vw)\cdot(\phi \vu)dxds&\le &\int_{t-\rho^2}^{t+\rho^2}\|\phi(s,\cdot)\|_{L^3_x(B_\rho)}\|\rot \vw(s,\cdot)\|_{L_{x}^2(B_\rho)} \|\vu(s,\cdot)\|_{L_{x}^6(B_\rho)}ds\\
&\leq &C\frac{\rho}{r}\int_{t-\rho^2}^{t+\rho^2}\|\rot \vw(s,\cdot)\|_{L_{x}^2(Q_\rho)}\|\vu(s,\cdot)\|_{\dot{H}^1_x(Q_\rho)}ds\\
&\leq &C\frac{\rho}{r}\|\rot \vw\|_{L_{t,x}^2(Q_\rho)}\|\vu\|_{L_{t}^2 \dot{H}^1_x(Q_\rho)},
\end{eqnarray*}
where we applied the Sobolev inequalities (see Corollary 9.14 of \cite{LibroBrezis}) and the Cauchy-Schwartz inequality in the time variable.
Since by (\ref{Def_Invariants}) we have $\|\vu\|_{L_{t}^2 \dot{H}^1_x(Q_\rho)}=\rho^{\frac{1}{2}} \alpha_\rho^\frac{1}{2}$, we conclude 
$$\int_{\mathbb{R}}\int_{\R}(\rot \vw)\cdot(\phi \vu)dxds\le C\frac{\rho^{\frac{3}{2}}}{r} 
\|\rot \vw\|_{L_{t,x}^2(Q_\rho)}\alpha_{\rho}^{\frac{1}{2}}.$$
\end{itemize}
Gathering all these estimates we obtain the inequality (\ref{FirstEstimate}) and this ends the proof of the Lemma \ref{lema_FirstEstimate}. \hfill$\blacksquare$\\

The inequality (\ref{FirstEstimate}) is important, but it will not be enough for our purposes as we need to study more in detail the pressure $p$. This variable only appears in the first equation of the system (\ref{MicropolarEquationsIntro}) and since we have the condition $div(\vu)=0$ and the vectorial identity $div(\rot \vw)\equiv 0$, by applying the divergence operator to the equation of $\vu$ in (\ref{MicropolarEquationsIntro}), we can write
\begin{eqnarray*}
div(\partial_t \vu) &=& div(\Delta \vu)  - div((\vu \cdot \vn)\vu)-div(\vn p) +\frac{1}{2}div(\rot\vw)\\
0&=& - div((\vu \cdot \vn)\vu)-\Delta p,
\end{eqnarray*}
from which we obtain the following equation for the pressure (see also (\ref{Eq_Pression}) above):
\begin{equation}\label{Eq_Pression1}
-\Delta p=div((\vu \cdot \vn)\vu)=div\big(div(\vu\otimes \vu)\big)=\sum_{i,j= 1}^3\partial_i
\partial_j (u_i u_j).
\end{equation}
Note that this previous equation for the pressure $p$ is exactly the same for the system (\ref{MicropolarEquationsIntro}) than for the classical Navier-Stokes equation. Thus, by the same ideas given in Proposition 4.3 of our previous work \cite{ChLl21} (see also Proposition 4.2 of \cite{ChHe21} or Lemma 13.3 of \cite{PGLR1}) we obtain the following result for the pressure: 
\begin{Lemma}\label{Lema_EstimatePressure}
Under the hypotheses of Theorem \ref{Teo_HolderRegularity} and with the notations given in (\ref{Def_Invariants}) for any $0<r<\frac{\rho}{2}<R$ we have the inequality
\begin{equation}\label{SecondEstimate}
\mathcal{P}_r^{\frac{2}{3}}\le C\biggl(\left(\frac{\rho}{r}\right)\left(\mathcal{A}_\rho
\alpha_\rho\right)^{\frac{1}{2}}+\left(\frac{r}{\rho}\right)^{\frac{2}{3}} 
\mathcal{P}_\rho^{\frac{2}{3}} \biggr).
\end{equation}
\end{Lemma}
For the sake of completeness we give the proof of this result.\\

\noindent {\bf Proof.}  We will start by proving the following estimate
\begin{equation}\label{Estimation_Pression_Intermediaire}
\|p\|_{L_{t,x}^{\frac{3}{2}}(Q_\sigma)} \leq C \left( \sigma^{\frac{1}{3}}\|\vu\|_{L_t^{\infty}L_x^{2} (Q_1)} \|\vn \otimes \vu\|_{L_{t,x}^{2} (Q_1)} +  
\sigma^{2}\|p\|_{L_{t,x}^{\frac{3}{2}} (Q_1)}\right),
\end{equation}
where $Q_\sigma$ and $Q_1$ are parabolic balls of radius $\sigma$ and $1$ respectively (the definition of such balls given in (\ref{Def_BolaParabolica})). To obtain this inequality we introduce $\eta : \R \longrightarrow [0, 1]$ a smooth function supported in the ball $B_1$ such that $\eta \equiv 1$ on the ball $B_{\frac35}$ and  $\eta \equiv 0$ outside the ball $B_{\frac45}$. Note in particular that on $Q_\sigma$ we have the identity $p=\eta p$. Now a straightforward calculation shows that we have the identity
$$ - \Delta (\eta p) = -\eta \Delta p + (\Delta \eta)p - 2  \sum^3_{i= 1}\partial_i ((\partial_i \eta) p),$$
from which we deduce the inequality
\begin{equation}\label{FormuleEtaPression}
\|p\|_{L_{t,x}^{\frac{3}{2}}(Q_\sigma)} =\|\eta p\|_{L_{t,x}^{\frac{3}{2}}(Q_\sigma)} \leq \underbrace{\left\|\frac{\big(-\eta \Delta p\big)}{(-\Delta)}\right\|_{L_{t,x}^{\frac{3}{2}}(Q_\sigma)}}_{(p_1)} + \underbrace{\left\|\frac{(\Delta \eta) p}{(-\Delta)}\right\|_{L_{t,x}^{\frac{3}{2}}(Q_\sigma)} }_{(p_2)} +2\sum^3_{i= 1}\underbrace{\left\|\frac{\partial_i ( (\partial_i \eta) p)}{(-\Delta)}\right\|_{L_{t,x}^{\frac{3}{2}}(Q_\sigma)}}_{(p_3)}.
\end{equation}
For the first term of (\ref{FormuleEtaPression}), since we have the equation (\ref{Eq_Pression1}) $\Delta p=-\displaystyle{\sum^3_{i,j= 1}}\partial_i \partial_j (u_i u_j)$ on $Q_\sigma$, if  we denote by $N_{i,j} = u_i (u_j - (u_j)_1)$ where $ (u_j)_1$ is the average of $u_j$ over the ball of radius $1$, since $\vu$ is divergence free we have $\displaystyle{\sum^3_{i,j= 1}}\partial_i \partial_j (u_i u_j)=\displaystyle{\sum^3_{i,j= 1}}\partial_i \partial_j N_{i,j}$ and thus we can write
\begin{eqnarray}
(p_1)&=&\left\|\frac{\big(-\eta \Delta p\big)}{(-\Delta)}\right\|_{L_{t,x}^{\frac{3}{2}}(Q_\sigma)}\leq C\left\|\frac{1}{(-\Delta)}\Big(\eta \,  \sum^3_{i,j= 1}\partial_i \partial_j  N_{i,j}  \Big)\right\|_{L_{t,x}^{\frac 32}(Q_\sigma)}\notag\\
&\leq &C\sum^3_{i,j= 1}\left\|\frac{1}{(-\Delta)} \left(\partial_i \partial_j(\eta N_{i,j} ) - \partial_i \big((\partial_j \eta) N_{i,j} \big) - \partial_j \big((\partial_i \eta) N_{i,j} \big) + 2(\partial_i \partial_j \eta) N_{i,j}\right)\right\|_{L_{t,x}^{\frac 32}(Q_\sigma)}\label{FormulePression_I}
\end{eqnarray}
Denoting by $\mathcal{R}_i=\frac{\partial_i}{\sqrt{-\Delta}}$ the usual Riesz transforms on $\mathbb{R}^3$, by the boundedness of these operators in Lebesgue spaces and using the support properties of the auxiliary function $\eta$, we have for the first term above:
\begin{eqnarray*}
\left\|\frac{\partial_i \partial_j}{(-\Delta)} \eta N_{i,j}(t,\cdot) \right\|_{L^{\frac32} (B_\sigma)} &\leq &\|\mathcal{R}_i \mathcal{R}_j (\eta N_{i,j} )(t,\cdot) \|_{L^{\frac32} (\R)}  \leq  C\|\eta N_{i,j}(t,\cdot)\|_{L^{\frac32} (B_1)}\\ 
&\leq &C\|u_i(t,\cdot) \|_{L^{2}(B_1)} \|u_j(t,\cdot) - (u_j)_1\|_{L^{6} (B_1)}\\ &\leq &C\|\vu(t,\cdot)\|_{L^{2} (B_1)} \|\vn \otimes \vu(t,\cdot)\|_{L^{2} (B_1)},
\end{eqnarray*}
where we used H\"older and Poincaré inequalities in the last line. Now taking the $L^{\frac32}$-norm in the time variable of the previous inequality we obtain\\
\begin{equation}\label{p11}
\left\|\frac{\partial_i \partial_j}{(-\Delta)} \eta N_{i,j} \right\|_{L_{t,x}^{\frac32} (Q_\sigma)}\leq C\sigma^{\frac{1}{3}}\|\vu\|_{L^{\infty}_tL^{2}_x (Q_1)} \|\vn \otimes \vu\|_{L_{t,x}^{2} (Q_1)}. 
\end{equation}
The remaining terms of (\ref{FormulePression_I}) can all be studied in a similar manner. Indeed, noting that $\partial_i \eta$ vanishes on $B_{\frac35} \cup B^c_{\frac45}$ and since $B_\sigma \subset B_{\frac12}\subset B_{\frac{3}{5}}$, using the integral representation for the operator $\frac{\partial_i}{(- \Delta )}$ we have for the second term of (\ref{FormulePression_I}) the estimate
\begin{eqnarray}
\left\|\frac{\partial_i}{(- \Delta )}\big((\partial_j\eta)N_{i,j}\big)(t,\cdot)\right\|_{L^{\frac32}(B_\sigma)} &\leq &C\sigma^2\left\|\frac{\partial_i}{(- \Delta )}\big((\partial_j\eta)N_{i,j}\big)(t,\cdot)\right\|_{L^{\infty}(B_\sigma)}\notag\\
&\leq &C \, \sigma^2 \left\|\int_{\{\frac35<|y|<\frac45\}} \frac{x_i - y_i }{|x-y|^3} \big((\partial_j \eta) N_{i,j} \big)(t,y)  \, dy\right\|_{L^{\infty}(B_\sigma)}\notag\\
& \leq &C \, \sigma^2 \|N_{i,j}(t,\cdot)\|_{L^{1} (B_1)}\label{KernelEstimate1}\\
& \leq &C \, \sigma^2 \|u_i(t,\cdot) \|_{L^{2} (B_1)} \|u_j(t,\cdot) - (u_j)_1\|_{L^{2}(B_1)} \notag\\
&\leq& C \,  \|\vu(t,\cdot)\|_{L^{2} (B_1)}\|\vn \otimes \vu(t,\cdot)\|_{L^{2} (B_1)},\notag
\end{eqnarray}
where we used the same ideas as previously and the fact that $0<\sigma<1$, and with the same arguments as in (\ref{p11}) before, taking the $L^{\frac32}$-norm in the time variable, we obtain
\begin{equation}\label{p12}
\left\|\frac{\partial_i}{(- \Delta )}\big((\partial_j\eta)N_{i,j}\big)\right\|_{L^{\frac32}_{t,x}(Q_\sigma)} \leq C\sigma^{\frac{1}{3}}\|\vu\|_{L^{\infty}_tL^{2}_x (Q_1)} \|\vn \otimes \vu\|_{L_{t,x}^{2} (Q_1)}. 
\end{equation}
A symmetric argument gives 
\begin{equation}\label{p13}
\left\|\frac{\partial_{j}}{(- \Delta )}\big((\partial_i\eta)N_{i,j}\big)\right\|_{L^{\frac32}_{t,x}(Q_\sigma)} \leq C\sigma^{\frac{1}{3}}\|\vu\|_{L^{\infty}_tL^{2}_x (Q_1)} \|\vn \otimes \vu\|_{L_{t,x}^{2} (Q_1)}, \end{equation}
and observing that the convolution kernel associated to the operator $\frac{1}{(-\Delta)}$ is $\frac{C}{|x|}$, following the same ideas we have for the last term of (\ref{FormulePression_I}) the inequality
\begin{equation}\label{p14}
\left\|\frac{(\partial_i \partial_j \eta) N_{i,j}}{(-\Delta)}\right\|_{L_{t,x}^{\frac 32}(Q_\sigma)}\leq C\sigma^{\frac{1}{3}}\|\vu\|_{L^{\infty}_tL^{2}_x (Q_1)} \|\vn \otimes \vu\|_{L_{t,x}^{2} (Q_1)}.
\end{equation}
Therefore, combining the estimates \eqref{p11}, \eqref{p12},  \eqref{p13} and  \eqref{p14} and getting back to (\ref{FormulePression_I}) we finally have:
\begin{eqnarray}
(p_1)=\left\|\frac{\big(-\eta \Delta p\big)}{(-\Delta)}\right\|_{L_{t,x}^{\frac 32}(Q_\sigma)}\leq C\left(\sigma^{\frac{1}{3}}\|\vu\|_{L^{\infty}_tL^{2}_x (Q_1)} \|\vn \otimes \vu\|_{L_{t,x}^{2} (Q_1)}\right).\label{FormulePression_I1}
\end{eqnarray}
We continue our study of expression (\ref{FormuleEtaPression}) and for the  term $(p_2)$ we first treat the space variable. Recalling the support properties of the auxiliary function $\eta$ and properties of the convolution kernel associated to the operator $\frac{1}{(-\Delta)}$, we can write as before (see (\ref{KernelEstimate1})):
$$\left\|\frac{(\Delta \eta) p(t,\cdot)}{(-\Delta)}\right\|_{L^{\frac{3}{2}}(B_\sigma)}\leq C\sigma^{2} \| p(t,\cdot)\|_{L^1(B_1)}\leq C\sigma^{2} \| p(t,\cdot)\|_{L^{\frac{3}{2}}(B_1)},$$
and thus, taking the $L^{\frac{3}{2}}$-norm in the time variable we obtain:
\begin{equation}\label{FormulePression_II}
(p_2)=\left\|\frac{(\Delta \eta) p}{(-\Delta)}\right\|_{L^{\frac{3}{2}}_{t,x}(Q_\sigma)}\leq C\sigma^{2} \|p\|_{L^{q_0}_{t,x}(Q_1)}.
\end{equation}
For the last term of expression (\ref{FormuleEtaPression}), following the same ideas developed in (\ref{KernelEstimate1}) we can write
$$\left\|\frac{\partial_i}{(-\Delta)}  (\partial_i \eta) p(t,\cdot) \right\|_{L^{\frac{3}{2}}(B_\sigma)}\leq C\sigma^{2} \| p(t,\cdot)\|_{L^1(B_1)}\leq C\sigma^{2} \|p(t,\cdot)\|_{L^{\frac{3}{2}}(B_1)},$$
and we obtain 
\begin{equation}\label{FormulePression_III}
(p_3)=\left\|\frac{\partial_i ( (\partial_i \eta) p )}{(-\Delta)}\right\|_{L_{t,x}^{\frac{3}{2}}(Q_\sigma)}\leq C\sigma^{2} \|p\|_{L^{\frac{3}{2}}_{t,x}(Q_1)}.
\end{equation}
Now, gathering the estimates (\ref{FormulePression_I1}), (\ref{FormulePression_II}) and (\ref{FormulePression_III}) we obtain the inequality 
$$\|p\|_{L_{t,x}^{\frac{3}{2}}  (Q_\sigma)} \leq C \left( \sigma^{\frac{1}{3}}\|\vu\|_{L_t^{\infty}L_x^{2} (Q_1)} \|\vn \otimes \vu\|_{L_{t,x}^{2} (Q_1)} +  
\sigma^{2}\|p\|_{L_{t,x}^{\frac{3}{2}} (Q_1)}\right),$$
which is (\ref{Estimation_Pression_Intermediaire}). With estimate at hand, it is quite simple to deduce inequality \eqref{SecondEstimate}. Indeed, if we fix $\sigma = \frac{r}{\rho} \leq \frac12$ and if we introduce the functions $p_\rho(t,x)=p(\rho^2t, \rho x)$ and $\vu_\rho(t,x)=\vu(\rho^2t, \rho x)$ then from (\ref{Estimation_Pression_Intermediaire}) we have
$$\|p_\rho\|_{L_{t,x}^{\frac{3}{2}}  (Q_{\frac{r}{\rho}})} \leq C \left( \left(\frac{r}{\rho}\right)^{\frac{1}{3}}\|\vu_\rho\|_{L_t^{\infty}L_x^{2} (Q_1)} \|\vn \otimes b_\rho\|_{L_{t,x}^{2} (Q_1)} +\left(\frac{r}{\rho}\right)^{2}\|p_\rho\|_{L_{t,x}^{\frac{3}{2}} (Q_1)}\right),$$
and by a convenient change of variable we obtain
$$\|p\|_{L_{t,x}^{\frac{3}{2}}  (Q_{r})}\rho^{-\frac{10}{3}} \leq C \left( \left(\frac{r}{\rho}\right)^{\frac{1}{3}}\rho^{-\frac32}\|\vu\|_{L_t^{\infty}L_x^{2} (Q_\rho)} \rho^{-\frac32}\|\vn \otimes \vu\|_{L_{t,x}^{2} (Q_\rho)} +\left(\frac{r}{\rho}\right)^{2}\rho^{-\frac{10}{3}}\|p\|_{L_{t,x}^{\frac{3}{2}} (Q_\rho)}\right).$$
Now, recalling that by (\ref{Def_Invariants}) we have the identities 
$$r^{\frac{4}{3}}\mathcal{P}_r^{\frac{2}{3}}=\|p\|_{L_{t,x}^{\frac{3}{2}}  (Q_r)}, \quad \rho^{\frac{1}{2}}\mathcal{A}_\rho^{\frac{1}{2}}=\|\vu\|_{L_t^{\infty}L_x^{2} (Q_\rho)} \quad \mbox{and}\quad \rho^{\frac{1}{2}}\alpha_\rho^{\frac{1}{2}}=\|\vn \otimes \vu\|_{L_{t,x}^{2} (Q_\rho)},$$
we obtain
$$\mathcal{P}_r^{\frac{2}{3}}\leq C\left(\left(\frac{\rho}{r}\right)(\mathcal{A}_\rho\beta_\rho)^{\frac{1}{2}}+\left(\frac{r}{\rho}\right)^{\frac{2}{3}}\mathcal{P}_\rho^{\frac{2}{3}}\right),$$
and this finishes the proof of Lemma \ref{Lema_EstimatePressure}.\hfill $\blacksquare$\\

Now, with the estimates (\ref{FirstEstimate}) and (\ref{SecondEstimate}) obtained in the previous lines, we will set up a general inequality that will help us to deduce the gain of integrability stated in Proposition \ref{Propo_FirstMorreySpace}. For this, we introduce the notations
\begin{equation}\label{Def_QuantiteIteration}
\mathbb{A}_r=\frac{1}{r^{2(1-\frac{5}{\tau_0})}}\left(\mathcal{A}_r+\alpha_r\right),\quad \mathbb{P}_r=\frac{1}{r^{\frac{3}{2}(1-\frac{5}{\tau_0})}}\mathcal{P}_r\quad \mbox{and}\quad 
\mathbb{O}_r=\mathbb{A}_r+\left(\left(\frac{r}{\rho}\right)^{\frac{15}{\tau_0}-\frac{15}{2}}\mathbb{P}_r\right)^{\frac{4}{3}},
\end{equation}
and we have the following result:
\begin{Lemma}\label{Lema_iteration}
Under the hypotheses of Theorem \ref{Teo_HolderRegularity}, for $0<r<\frac{\rho}{2}<R$ there exists a constant $\epsilon>0$ such that
\begin{equation}\label{thetaite1}
\mathbb{O}_{r}(t_0,x_0)\le \frac{1}{2}\mathbb{O}_{\rho}(t_0,x_0)+\epsilon,
\end{equation}
where the point $(t_0,x_0)\in Q_R$ is given by the hypothesis (\ref{Condition_PetitEpsilon}).
\end{Lemma}
Note that this result allows us to control the information over smaller parabolic balls by the information over bigger parabolic balls and this will be the key to obtain the whished gain of integrability.\\

\noindent\textbf{Proof.} First, by the estimate (\ref{FirstEstimate}) we can write
\begin{eqnarray}
\mathbb{A}_r&=&\frac{1}{r^{2(1-\frac{5}{\tau_0})}}\left(\mathcal{A}_r+\alpha_r\right)\notag\\
&\leq& \frac{C}{r^{2(1-\frac{5}{\tau_0})}}\Big(\frac{r^2}{\rho ^2} \mathcal{A}_\rho+\frac{\rho^2}{r^2} \alpha_\rho^{\frac{1}{2}}\mathcal{A}_\rho+\frac{\rho ^2}{r^2}\mathcal{P}_\rho^{\frac{2}{3}}(\mathcal{A}_\rho+\alpha_\rho)^{\frac{1}{2}}+\frac{\rho^{\frac{3}{2}}}{r}\|\rot \vw\|_{L_{t,x}^2(Q_\rho)}\alpha_{\rho}^{\frac{1}{2}}\Big),\label{EstimationPourIteration}
\end{eqnarray}
and we will treat each one of the previous terms separately. Indeed, 
\begin{itemize}
\item For the first term of \eqref{EstimationPourIteration} we have
$$\frac{1}{r^{2(1-\frac{5}{\tau_0})}}\left(\frac{r^2}{\rho ^2}\mathcal{A}_\rho\right) \le \frac{1}{r^{2(1-\frac{5}{\tau_0})}}\frac{r^2}{\rho ^2}\rho^{2(1-\frac{5}{\tau_0})}\mathbb{A}_\rho =\left(\frac{r}{\rho}\right)^{\frac{10}{\tau_0}}\mathbb{A}_\rho.$$
\item For the second term of (\ref{EstimationPourIteration}), using the definition of $\mathbb{A}_\rho$ given in (\ref{Def_QuantiteIteration}), we obtain
$$\frac{1}{r^{2(1-\frac{5}{\tau_0})}}\left(\frac{\rho^2}{r^2} \alpha_\rho^{\frac{1}{2}}\mathcal{A}_\rho\right)\leq \frac{1}{r^{2(1-\frac{5}{\tau_0})}}\left(\frac{\rho^2}{r^2} \alpha_\rho^{\frac{1}{2}}\rho^{2(1-\frac{5}{\tau_0})}\mathbb{A}_\rho \right)=\left(\frac{\rho}{r}\right)^{4-\frac{10}{\tau_0}}\mathbb{A}_\rho \alpha_\rho^{\frac{1}{2}}.$$
\item The third term of (\ref{EstimationPourIteration}) follows essentially the same arguments as above and by the definition of the quantities  $\mathbb{A}_\rho$  and  $\mathbb{P}_\rho$ given in (\ref{Def_QuantiteIteration}) we can write
$$\frac{1}{r^{2(1-\frac{5}{\tau_0})}}\left(\frac{\rho ^2}{r^2}\mathcal{P}_\rho^{\frac{2}{3}}( \mathcal{A}_\rho+\alpha_\rho)^{\frac{1}{2}}\right)\leq \left(\frac{\rho}{r}\right)^{4-\frac{10}{\tau_0}}\mathbb{P}_\rho^{\frac{2}{3}}\mathbb{A}_\rho^{\frac{1}{2}}.$$
\item Finally, for the last term of (\ref{EstimationPourIteration}), we have
$$\frac{1}{r^{2(1-\frac{5}{\tau_0})}}\left(\frac{\rho^{\frac{3}{2}}}{r} 
\|\rot \vw\|_{L_{t,x}^2(Q_\rho)}\alpha_{\rho}^{\frac{1}{2}}\right)\leq \left(\frac{\rho}{r}\right)^{3-\frac{10}{\tau_0}} \rho^{\frac{10}{\tau_0}-\frac{3}{2}}\alpha_{\rho}^{\frac{1}{2}}\|\rot \vw\|_{L_{t,x}^2(Q_\rho)}.$$
\end{itemize}
Thus, gathering all these estimates, we have 
\begin{equation}\label{estimateA}
\mathbb{A}_r\leq C\Bigg(\left(\frac{r}{\rho}\right)^{\frac{10}{\tau_0}}\mathbb{A}_\rho+\left(\frac{\rho}{r}\right)^{4-\frac{10}{\tau_0}}\mathbb{A}_\rho \alpha_\rho^{\frac{1}{2}}+\left(\frac{\rho}{r}\right)^{4-\frac{10}{\tau_0}}\mathbb{P}_\rho^{\frac{2}{3}}\mathbb{A}_\rho^{\frac{1}{2}}+ \left(\frac{\rho}{r}\right)^{3-\frac{10}{\tau_0}} \rho^{\frac{10}{\tau_0}-\frac{3}{2}}\alpha_{\rho}^{\frac{1}{2}}\|\rot \vw\|_{L_{t,x}^2(Q_\rho)}\Bigg).
\end{equation}
Now, for the pressure, from the inequality \eqref{SecondEstimate} we can write
$$\mathbb{P}_r=\frac{1}{r^{\frac{3}{2}(1-\frac{5}{\tau_0})}}\mathcal{P}_r\le \frac{C}{r^{\frac{3}{2}(1-\frac{5}{\tau_0})}} \biggl(\left(\frac{\rho}{r}\right)^{\frac32}\mathcal{A}_\rho^{\frac{3}{4}}\alpha_\rho^{\frac{3}{4}}+ \left(\frac{r}{\rho}\right)\mathcal{P}_\rho \biggr),$$
and by the Young inequality and using the definition of $\mathbb{A}_\rho$ given in (\ref{Def_QuantiteIteration}) we obtain for the first term of the right-hand side above:
$$\frac{1}{r^{\frac{3}{2}(1-\frac{5}{\tau_0})}}\left(\frac{\rho}{r}\right)^{\frac{3}{2}}\mathcal{A}_\rho^{\frac{3}{4}} \alpha_\rho^{\frac{3}{4}}\leq\frac{1}{r^{\frac{3}{2}(1-\frac{5}{\tau_0})}}\left(\frac{\rho}{r}\right)^{\frac{3}{2}}\rho^{\frac{3}{2}(1-\frac{5}{\tau_0})}(\mathbb{A}_\rho\alpha_\rho)^{\frac{3}{4}}=\left(\frac{\rho}{r}\right)^{3-\frac{15}{2\tau_0}}(\mathbb{A}_\rho\alpha_\rho)^{\frac{3}{4}},$$
and using the fact that $\frac{1}{r^{\frac{3}{2}(1-\frac{5}{\tau_0})}} \left(\frac{r}{\rho}\right)\mathcal{P}_\rho= \left(\frac{\rho}{r}\right)^{\frac{1}{2}-\frac{15}{2\tau_0}}\mathbb{P}_\rho$ (by the definition of $\mathbb{P}_\rho$ given in (\ref{Def_QuantiteIteration})), we conclude that
\begin{align}\label{estimateP}
\mathbb{P}_r\le C\biggl(\left(\frac{\rho}{r}\right)^{3-\frac{15}{2\tau_0}}
(\mathbb{A}_\rho \alpha_\rho)^{\frac{3}{4}}+\left(\frac{\rho}{r}\right)^{\frac{1}{2}-\frac{15}{2\tau_0}}\mathbb{P}_\rho\biggr).
\end{align}
With the estimates (\ref{estimateA}) and (\ref{estimateP}) at hand, we will now introduce a relationship between the parameters $r$ and $\rho$: indeed, let us fix $0<\kappa\ll \frac{1}{2}$ a real number and consider $r=\kappa \rho$, then, by the definition of the quantity $\mathbb{O}_r$ given in (\ref{Def_QuantiteIteration}) we obtain:
\begin{eqnarray}
\mathbb{O}_r&=&\mathbb{A}_r+\left(\kappa^{\frac{15}{\tau_0}-\frac{15}{2}}\mathbb{P}_r\right)^{\frac{4}{3}}\leq C\Bigg(\underbrace{\kappa^{\frac{10}{\tau_0}}\mathbb{A}_\rho+\kappa^{\frac{10}{\tau_0}-4}\mathbb{A}_\rho \alpha_\rho^{\frac{1}{2}}}_{(1)}+\underbrace{\kappa^{\frac{10}{\tau_0}-4}\mathbb{P}_\rho^{\frac{2}{3}}\mathbb{A}_\rho^{\frac{1}{2}}}_{(2)}+\underbrace{\kappa^{\frac{10}{\tau_0}-3} \rho^{\frac{10}{\tau_0}-\frac{3}{2}}\alpha_{\rho}^{\frac{1}{2}}\|\rot \vw\|_{L_{t,x}^2(Q_\rho)}}_{(3)}\Bigg)\notag\\
&&+C\underbrace{\biggl(\kappa^{\frac{45}{2\tau_0}-\frac{21}{2}}(\mathbb{A}_\rho \alpha_\rho)^{\frac{3}{4}}+\kappa^{\frac{45}{2\tau_0}-8}\mathbb{P}_\rho\biggr)^{\frac43}}_{(4)}.\label{EstimationTheta}
\end{eqnarray}
We will rewrite now each one of the previous terms:
\begin{itemize}
\item Since by (\ref{Def_QuantiteIteration}) we have $\mathbb{A}_\rho\leq \mathbb{O}_\rho$, it is then easy to see that the term (1) above can be controlled in the following manner:
$$\kappa^{\frac{10}{\tau_0}}\mathbb{A}_\rho+\kappa^{\frac{10}{\tau_0}-4}\mathbb{A}_\rho \alpha_\rho^{\frac{1}{2}}\leq \big(\kappa^{\frac{10}{\tau_0}}+\kappa^{\frac{10}{\tau_0}-4} \alpha_\rho^{\frac{1}{2}}\big)\mathbb{O}_\rho.$$
\item For the quantity (2) in (\ref{EstimationTheta}), using Young's inequality and the relationships given in (\ref{Def_QuantiteIteration}), we observe that
\begin{eqnarray*}
\kappa^{\frac{10}{\tau_{0}}-4}  \mathbb{P}_\rho^{\frac{2}{3}} \mathbb{A}_\rho^{\frac{1}{2}}&=&\kappa^{\frac{10}{\tau_{0}}-4} \left(\kappa^{5(\frac{1}{\tau_{0}} -\frac{1}{2})} \mathbb{P}_\rho^{\frac{2}{3}} \times\kappa^{5(\frac{1}{2} - \frac{1}{\tau_{0}} )} \mathbb{A}_\rho^{\frac{1}{2}}\right)\leq\kappa^{\frac{10}{\tau_{0}}-4} \left( \kappa^{10(\frac{1}{2} - \frac{1}{\tau_{0}} )} \mathbb{A}_\rho + \kappa^{10(\frac{1}{\tau_{0}} - \frac{1}{2} )} \mathbb{P}_\rho^{\frac{4}{3}} \right)\notag \\
&\leq & \kappa \left(\mathbb{A}_\rho+ \left(\kappa^{\frac{15}{\tau_{0}} -\frac{15}{2}}\mathbb{P}_\rho\right)^{\frac{4}{3}} \right) \leq \kappa \, \mathbb{O}_\rho.
\end{eqnarray*}
\item For the term $(3)$ of (\ref{EstimationTheta}), using the fact that $\frac{10}{\tau_0}>\frac{3}{2}$ (recall the hypothesis of Proposition \ref{Propo_FirstMorreySpace}: we have $\frac{5}{1-\alpha}<\tau_0<\frac{20}{3}$) and that $0<\rho<R<1$, we obtain $ \rho^{\frac{10}{\tau_0}-\frac32}<1$, and thus
\begin{equation*}
\kappa^{\frac{10}{\tau_0}-3} \rho^{\frac{10}{\tau_0}-\frac32}\alpha_{\rho}^{\frac{1}{2}}\|\rot \vw\|_{L_{t,x}^2(Q_\rho)}\le  \kappa^{\frac{10}{\tau_0}-3} \alpha_{\rho}^{\frac{1}{2}}\|\rot \vw\|_{L_{t,x}^2(Q_\rho)}.
\end{equation*} 
\item For the last term of (\ref{EstimationTheta}), since $\left(\kappa^{\frac{15}{\tau_0}-\frac{15}{2}}\mathbb{P}_\rho\right)^{\frac{4}{3}}\leq \mathbb{O}_\rho$ and $\mathbb{A}_\rho\leq \mathbb{O}_\rho$, we have
\begin{eqnarray*}
\biggl(\kappa^{\frac{45}{2\tau_0}-\frac{21}{2}}(\mathbb{A}_\rho \alpha_\rho)^{\frac{3}{4}}+\kappa^{\frac{45}{2\tau_0}-8}\mathbb{P}_\rho\biggr)^{\frac43}\leq C\biggl(\kappa^{\frac{30}{\tau_0}-14}\mathbb{A}_\rho \alpha_\rho+(\kappa^{\frac{45}{2\tau_0}-8}\mathbb{P}_\rho)^{\frac43}\biggr)\leq C\bigg(\kappa^{\frac{30}{\tau_0}-14} \alpha_\rho+\kappa^{\frac{10}{\tau_0}-\frac{2}{3}}\bigg)\mathbb{O}_\rho.
\end{eqnarray*}
\end{itemize}
Gathering these estimates we finally obtain
\begin{equation}\label{thetaite2}
\mathbb{O}_r \le \Biggl(\kappa^{\frac{10}{\tau_0}}+\kappa^{\frac{10}{\tau_0}-4}\alpha_\rho^{\frac{1}{2}}+\kappa+\kappa^{\frac{30}{\tau_0}-14}\alpha_\rho+\kappa^{\frac{10}{\tau_0}-\frac{2}{3}}\Biggr)\mathbb{O}_\rho+\kappa^{\frac{10}{\tau_0}-3} \alpha_{\rho}^{\frac{1}{2}}\|\rot \vw\|_{L_{t,x}^2(Q_\rho)}.
\end{equation}
Futhermore, we claim that we have
\begin{equation}\label{LimsupHypo1}
\Biggl(\kappa^{\frac{10}{\tau_0}}+\kappa^{\frac{10}{\tau_0}-4}\alpha_\rho^{\frac{1}{2}}+\kappa+\kappa^{\frac{30}{\tau_0}-14}\alpha_\rho+\kappa^{\frac{10}{\tau_0}-\frac{2}{3}} \Biggr)\le\frac{1}{2}.
\end{equation}
Indeed, since $\kappa=\frac{r}{\rho}\ll\frac{1}{2}$ is a fixed small parameter and since $\frac{10}{\tau_0}-\frac{2}{3}>0$ (recall again that $\frac{5}{1-\alpha}<\tau_0<\frac{20}{3}$), then the quantities $\kappa^{\frac{10}{\tau_0}}$, $\kappa$ and $\kappa^{\frac{10}{\tau_0}-\frac{2}{3}}$ in the previous formula are small. Now, using the fact that we have the control $\alpha_\rho \le \epsilon^*$ which is given in the hypothesis (\ref{Condition_PetitEpsilon}) where $\epsilon^*>0$ is small enough, then the terms $\kappa^{\frac{10}{\tau_0}-4}\alpha_\rho^{\frac{1}{2}}$ and $\kappa^{\frac{30}{\tau_0}-14}\alpha_\rho$   can be made small enough and thus we obtain \eqref{LimsupHypo1}.
To continue, noting that the quantity $\|\rot \vw\|_{L_{t,x}^2(Q_\rho)}$ is bounded since $\vw\in L_t^{\infty}L_x^2(Q_R)\cap L_t^2\dot H_x^1(Q_R)$, we can apply the same ideas used previously (\emph{i.e.} $\alpha_\rho \le \epsilon^*\ll1$) to obtain
$$\kappa^{\frac{10}{\tau_0}-3} \alpha_{\rho}^{\frac{1}{2}}\|\rot \vw\|_{L_{t,x}^2(Q_\rho)}<\epsilon.$$
Then, with these estimates at hand and coming back to (\ref{thetaite2}) we conclude that $\mathbb{O}_r\le \frac{1}{2}\mathbb{O}_\rho+\epsilon$ and Lemma \ref{Lema_iteration} is proven. \hfill $\blacksquare$\\

Lemma \ref{Lema_iteration} paved the way to obtain some Morrey information for the velocity $\vu$ that will be crucial. Indeed, from the definition of Morrey spaces given in (\ref{DefMorreyparabolico}) we only need to prove that for all radius $r>0$ such that $r<R_1\le \frac{R}{2}$ and  $(t,x)\in Q_{R_1}(t_0,x_0)$, we have
\begin{equation}\label{morrey1}
\int_{Q_r(t,x)} |\vu|^3 dy ds \le C r^{5(1-\frac{3}{\tau_0})},
\end{equation}
and this will imply that $\mathds{1}_{Q_{R_1}}\vu\in \M^{3,\tau_0}(\mathbb{R}\times \R)$. In order to obtain the control (\ref{morrey1}), by the definitions given in (\ref{Def_Invariants}) and by the estimate (\ref{Estimation_normeLambda3}), we observe that
\begin{equation*}
\int_{Q_r(t,x)} |\vu|^3dy ds=r^2\lambda_r(t,x) \le  r^2(\mathcal{A}_r(t,x)+\alpha_r(t,x))^{\frac{3}{2}}.
\end{equation*}
Hence, it is then enough to prove for all $0<r< {R_1}<\frac{R}{2}<R<1$ and  $(t,x)\in Q_{R_1}$ that one has the control
\begin{equation*}
\mathcal{A}_r(t,x)+\alpha_r(t,x) \le C r^{2(1-\frac{5}{\tau_0})}.
\end{equation*}
Recalling the definition of the quantity $\mathbb{A}_r$ given in (\ref{Def_QuantiteIteration}), we easily  see that the condition (\ref{morrey1}) above is equivalent to prove that there exists some $R_1$ and  $0<\kappa\ll \frac12$ such that for all $n\in \mathbb{N}$ and $(t,x)\in Q_{R_1}(t_0, x_0)$, we have estimates:
\begin{equation}\label{iterative}
\mathbb{A}_{\kappa^n R_1}(t,x)\le C.
\end{equation}
Note that, for any radius $r$ such that $0<r<R_1<\min\{\frac{R}{2},dist(\partial Q_R, (t_0,x_0))\}$ (and since we have $Q_{R_1}(t_0,x_0)\subset Q_{R}$) by the hypotheses of the Theorem \ref{Teo_HolderRegularity}, we have the bounds 
$$\|\vu\|_{L_t^\infty L_x^2(Q_r(t_0,x_0))}\le \|\vu\|_{L_t^\infty L_x^2(Q_R)}<+\infty,\quad \|\vn \otimes\vu\|_{L_{t,x}^2(Q_r(t_0,x_0))}\le \|\vn\otimes\vu\|_{L_{t,x}^2(Q_R)}<+\infty,$$
and $ \|p\|_{L_{t,x}^{\frac{3}{2}}(Q_r(t_0,x_0))}\le \|p\|_{L_{t,x}^{\frac{3}{2}}(Q_R)}<+\infty$. Then, by the notations introduced in (\ref{Def_Invariants}), we have the uniform bounds $\underset{0<r<R}{\sup}\biggl\{ r \mathcal{A}_r,r \alpha_r, r^2\mathcal{P}_r\biggr\}<+\infty$ from which we can deduce by the definition of the quantities $\mathbb{A}_\rho(t_0,x_0)$ and $\mathbb{P}_\rho (t_0,x_0)$ given in (\ref{Def_QuantiteIteration}), the uniform bounds
\begin{equation*}
\sup_{0<r<R}r^{3-\frac{10}{\tau_0}}\mathbb{A}_r(t_0,x_0)<+\infty,\quad\text{and} \quad \sup_{0<r<R}r^{5-\frac{3}{2}(1+\frac{5}{\tau_0})}\mathbb{P}_r (t_0,x_0)<+\infty.
\end{equation*}
Thus, there exists a radius $0<r_0< R$ small such that, by the estimates above, the quantities $\mathbb{A}_{r_0}$ and $\mathbb{P}_{r_0}$ are bounded: indeed, recall that we have $\tau_0>\frac{5}{1-\alpha}>5$ (where $0<\alpha\ll 1$) and this implies that all the powers of $r$ in the expression above are positive. As a consequence of this fact, by (\ref{Def_QuantiteIteration}) the quantity $\mathbb{O}_{r_0}$ is itself bounded. Remark also that, if $r_0$ is small enough, then the inequality (\ref{thetaite1}) holds true and we can write $\mathbb{O}_{\kappa r_0}(t_0,x_0)\le \frac{1}{2}\mathbb{O}_{r_0}(t_0,x_0)+\epsilon$. We can iterate this process and we obtain for all $n>1$,
$$\mathbb{O}_{\kappa^nr_0}(t_0,x_0)\le \frac{1}{2^n}\mathbb{O}_{r_0}(t_0,x_0)+\epsilon\sum_{j=0}^{n-1}2^{-j},$$
and therefore there exists $N\ge 1$ such that for all $n\ge N$ we have
$\mathbb{O}_{\kappa^nr_0}(t_0,x_0)\le 4\epsilon$ from which we obtain (using the definition of $\mathbb{O}_r$ given in (\ref{Def_QuantiteIteration})) that
$$\mathbb{A}_{\kappa^Nr_0}(t_0,x_0)\le \frac{1}{8}C \quad \text{and}\quad
\mathbb{P}_{\kappa^Nr_0}(t_0,x_0)\le \frac{1}{32}C.$$
This information is centered at the point $(t_0,x_0)$, in order to treat
the uncentered bound, we can let $\frac{1}{2}\kappa^Nr_0$ to be the radius $R_1$
we want to find, thus for all points $(t,x)\in Q_{R_1}(t_0,x_0)	$ we have that
$Q_{R_1}\subset Q_{2R_1}(t_0,x_0)$, which implies 
\begin{equation*}
\mathbb{A}_{R_1}(t,x)\le 2^{3-\frac{10}{\tau_0}} \mathbb{A}_{2R_1}(t_0,x_0)\le 8 \mathbb{A}_{2R_1}(t_0,x_0)\le 8 \mathbb{A}_{\kappa^N \rho }(t_0,x_0)<C,
\end{equation*}
and $\mathbb{P}_{R_1}(t,x)\le 2^{5-\frac{3}{2}(1+\frac{5}{\tau_0})}\mathbb{P}_{2R_1}(t_0,x_0)\le 32\mathbb{P}_{2R_1}(t_0,x_0)\le 8 \mathbb{P}_{\kappa^N r }(t_0,x_0)<C$. Having obtained these bounds, by the definition of $\mathbb{O}_{R_1}$, we thus get $\mathbb{O}_{R_1}(t,x)\le C$. Applying the Lemma \ref{Lema_iteration} and iterating once more, we find that the same will be true for $\kappa R_1$ and then, for all $\kappa^n R_1$, $n\in \mathbb{N}$. Since by definition we have  $\mathbb{A}_{\kappa^n R_1}(t,x)\le  \mathbb{O}_{\kappa^n R_1}(t,x)$ we have finally obtained the estimate $\mathbb{A}_{\kappa^n R_1}(t,x)\le C$ and the inequality (\ref{iterative}) is proven which implies the Proposition \ref{Propo_FirstMorreySpace}.\hfill$\blacksquare$\\
\begin{Corollary}\label{corolarioMorrey}
Under the hypotheses of Proposition \ref{Propo_FirstMorreySpace}, we also have the following local control:
\begin{equation}\label{ConlusioncorolarioMorrey}
\mathds{1}_{Q_{R_1}(t_0,x_0)}\vn\otimes \vu\in \mathcal{M}^{2,\tau_1}_{t,x}(\mathbb{R}\times \R),\quad \mbox{with}\quad  \frac{1}{\tau_1}=\frac{1}{\tau_0}+\frac{1}{5}.
\end{equation}
\end{Corollary}
\textbf{Proof.} In the previous results we have proved the estimate (\ref{iterative}). Let us recall now that, by the definition of the quantity $\mathbb{A}_r$ given in (\ref{Def_QuantiteIteration}), we can easily deduce for all $0<r\le R_1$ and $(t,x)\in Q_{R_1}$ the control $\alpha_r\leq Cr^{2(1-\frac{5}{\tau_0})}$ which can we rewritten as
$$\frac{1}{r}\bigg(\int _{Q_r(t,x)}|\vn \otimes \vu|^2dyds\bigg)\leq Cr^{2(1-\frac{5}{\tau_0})}.$$
Thus, since $\frac{1}{\tau_1}=\frac{1}{\tau_0}+\frac{1}{5}$, for all $0<r\le R_1$ and $(t,x)\in Q_{R_1}(t_0,x_0)$, we have the estimate
$$\int _{Q_r}|\vn \otimes \vu|^2dyds\leq C r^{3-\frac{10}{\tau_0}}= C r^{5(1-\frac{2}{\tau_1})},$$
and by the definition of Morrey spaces given in (\ref{DefMorreyparabolico}), we obtain that $\mathds{1}_{Q_{R_1}(t_0,x_0)}\vn \otimes \vu\in \mathcal{M}^{2,\tau_1}_{t,x}(\mathbb{R}\times \R)$.\hfill$\blacksquare$
\section{A (second) partial gain of information for the variable $\vu$}\label{Sec_Propriedades_Vu2}
This first gain of integrability information stated in Proposition \ref{Propo_FirstMorreySpace} is fundamental for our theory to work, however it is not enough since we only obtain a ``small'' control\footnote{In terms of the indexes of the Morrey spaces involved in Proposition \ref{Propo_FirstMorreySpace}.} for the variable $\vu$ and without any information on the variable $\vw$ we can not go very far: now we will see how to obtain some further control on $\vw$ and how it is possible to reinject this information in the study of the variable $\vu$. Indeed, in our recent article \cite{ChLl22} we proved the following result which gives some mild control over the variable $\vw$: 
\begin{Theorem}\label{Teo_InterdepenciaMicropolar}
Let  $(\vu,p, \vw)$ be a weak solution of the micropolar equations (\ref{MicropolarEquationsIntro}) over a parabolic ball $Q_{R}$ of the form (\ref{Def_BolaParabolica}) for some fixed radius $R>0$. Assume that $\vu,\vw\in L_t^\infty L_x^2 \cap L_t^2 \dot H_x ^1(Q_R)$ and $p\in \mathcal{D}'_{t,x}(Q_R)$.
Suppose in addition that for some $0<R_1<R$ we have
\begin{equation}\label{Hypothese1}
\mathds{1}_{Q_{R_1}}\vu \in \M_{t,x}^{p_0,q_0}(\mathbb{R}\times\R)\quad with\;\;
2<p_0\le q_0,\; 5<q_0\le 6,
\end{equation}
then
\begin{itemize}
\item[1)] for a parabolic ball $Q_{\mathfrak{r}_1}$, with $0<\mathfrak{r}_1<R_1$ we have
\begin{equation*}\label{Conclusion1}
\mathds{1}_{Q_{\mathfrak{r}_1}}\vu \in L_{t,x}^{q_0}(\mathbb{R}\times\R),\qquad 5<q_0\le 6,
\end{equation*}
\item[2)] for a parabolic ball  $Q_{\mathfrak{r}_2}$, with  $0<\mathfrak{r}_2<\mathfrak{r}_1<R_1$ we have  
\begin{equation*}\label{Conclusion2}
\mathds{1}_{Q_{\mathfrak{r}_2}}\vw \in L_{t,x}^{q_0}(\mathbb{R}\times\R),
\end{equation*}
for  $5<q_0\leq 6$.
\end{itemize}
\end{Theorem}
As we can see, this result gives an interesting improvement of integrability for \emph{both} variables $\vu$ and $\vw$ as long as we have the hypothesis (\ref{Hypothese1}), but this is precisely the conclusion of Proposition \ref{Propo_FirstMorreySpace}: indeed, over a small parabolic ball $Q_{R_1}(t_0,x_0)$ we do have $\mathds{1}_{Q_{R_1}(t_0,x_0)} \vu\in \mathcal{M}_{t,x}^{3,\tau_0}(\mathbb{R}\times \R)$ and it is enough to remark that we have here $p_0=3$ and $q_0=\tau_0$ with $\frac{5}{1-\alpha}<\tau_0<\frac{20}{3}$ and this last parameter can be chosen such that $\tau_0= 6<\frac{20}{3}$. Thus, we deduce that 
\begin{equation}\label{Gain_IntegrabiliteUW}
\mathds{1}_{Q_{\mathfrak{r}_1}(t_0,x_0)}\vu \in L_{t,x}^{6}(\mathbb{R}\times\R)\quad \mbox{and}\quad \mathds{1}_{Q_{\mathfrak{r}_2}(t_0,x_0)}\vw \in L_{t,x}^{6}(\mathbb{R}\times\R),
\end{equation}
where $\mathfrak{r}_2<\mathfrak{r}_1<R_1<R<1$.\\

Note that from the initial setting $\vu, \vw\in L^\infty_tL^2_x\cap L^2_t\dot{H}^1_x$, the controls stated in (\ref{Gain_IntegrabiliteUW}) provide a better integrability information and we will see now how to improve the Morrey information given in Proposition \ref{Propo_FirstMorreySpace} for the variable $\vu$:
\begin{Proposition}\label{Propo_GainenMorreySigma}
Under the hypotheses of Theorem \ref{Teo_HolderRegularity} and within the framework of Proposition \ref{Propo_FirstMorreySpace}, there exists a radius $R_2$ with $0<R_2<\mathfrak{r}_2<\mathfrak{r}_1<R_1<R<1$ such that
\begin{equation*}
\mathds{1}_{Q_{R_2}(t_0,x_0)}\vu \in \M_{t,x}^{3,\sigma}(\mathbb{R}\times \R),
\end{equation*}
for some $\sigma$ close to $\tau_0=6$ such that $\tau_0<\sigma$.
\end{Proposition}
\textbf{Proof of the Proposition \ref{Propo_GainenMorreySigma}.} In order to obtain this small additional gain of integrability we will first localize the variable $\vu$ in a suitable manner and then we will study its evolution: the wished result will then be deduced from the Duhamel formula and from all the available information over $\vu$. Let us start fixing the parameters $\mathfrak{R}_c, \mathfrak{R}_b, \mathfrak{R}_a$ such that
$$0<R_2<\mathfrak{R}_c<\mathfrak{R}_b<\mathfrak{R}_a<\mathfrak{r}_2<\mathfrak{r}_1<R_1,$$ 
with the associated parabolic balls $Q_{R_2}\subset Q_{\mathfrak{R}_c}\subset Q_{\mathfrak{R}_b}\subset Q_{\mathfrak{R}_a}\subset Q_{R_1}$ (all centered in the point $(t_0,x_0)$). Consider now $ \phi, \psi: \mathbb{R}\times\mathbb{R}^3\longrightarrow \mathbb{R}$ two  non-negative functions such that $\phi, \psi\in \mathcal{C}_0^{\infty}
(\mathbb{R}\times \mathbb{R}^3)$ and such that
\begin{equation}\label{ProprieteLocalisation_IterationVU}
 \phi \equiv 1\;\; \text{over}\; \; Q_{\mathfrak{R}_c},\; \; supp( \phi)\subset Q_{\mathfrak{R}_b} \quad \mbox{and}\quad \psi \equiv 1\;\; \text{over}\; \; Q_{\mathfrak{R}_a},\; \; supp(\psi)\subset Q_{ \mathfrak{r}_2}.
\end{equation}
Using these auxiliar functions we will study the evolution of the variable $\vv= \phi\,\vu$ given by the system
\begin{equation}\label{equationV}
\begin{cases}
\partial_t\vv=\Delta\vv+\vV,\\[3mm]
\vv(0,x)=0,
\end{cases}
\end{equation}
where we have  
\begin{equation}\label{Formule_equationV}
\vV=(\partial_t\phi- \Delta \phi)\vu-2\sum_{i=1}^3 (\partial_i \phi)(\partial_i \vu)-\phi(\vu \cdot \vn)\vu-2\phi\vn p+\phi(\rot\vw).
\end{equation}
We will now rewrite the term $\phi\vn p$ above in order to avoid a direct derivative over the pressure. Indeed, as we have the identity $p= \psi p$ over $Q_{\mathfrak{R}_a}$, then over the smaller ball $Q_{R_2}$  (recalling that $ \psi=1$ over $Q_{R_2}$ by (\ref{ProprieteLocalisation_IterationVU}) since $Q_{R_2}\subset Q_{\mathfrak{R}_a}$), we can write $-\Delta( \psi p)=- \psi \Delta p+(\Delta  \psi)p-2\displaystyle{\sum_{i=1}^3}\partial_i((\partial_i \psi)p)$
from which we deduce the identity 
\begin{equation}\label{ExpressionPression6}
 \phi \vn p= \phi\frac{\vn (- \psi \Delta p)}{(-\Delta)}+ \phi\frac{\vn ((\Delta  \psi)p)}{(-\Delta)}-2\sum_{i=1}^3 \phi\frac{\vn (\partial_i((\partial_i \psi)p))}{(-\Delta)}.
\end{equation}
At this point we recall that we have by (\ref{Eq_Pression}) the following equation for the pressure $\Delta p = -\displaystyle{\sum_{i,j=1}^3}\partial_i\partial_j\left(u_i u_j\right)$ and thus, the first term of the right-hand side of the previous formula can be written in the following manner:
\begin{eqnarray*}
 \phi\frac{\vn (- \psi \Delta p)}{(-\Delta)}&= & \phi\frac{\vn}{(-\Delta)}\left( \sum_{i,j=1}^3 \psi\big(\partial_i\partial_j u_i u_j\big)\right)\notag\\
&=& \sum_{i,j=1}^3  \phi\frac{\vn}{(-\Delta)}\bigg(\partial_i\partial_j( \psi u_i u_j)\bigg)-\sum_{i,j=1}^3  \phi\frac{\vn}{(-\Delta)}\bigg(\partial_i((\partial_j  \psi)u_i u_j)+\partial_j((\partial_i  \psi)u_i u_j)-(\partial_i \partial_j  \psi)(u_i u_j)\bigg),
\end{eqnarray*}
Recalling that by construction of the auxiliar functions $ \phi$, $ \psi$ given in (\ref{ProprieteLocalisation_IterationVU}) we have the identity $ \phi  \psi= \phi$, we can write for the first term above:
$$ \phi\frac{\vn}{(-\Delta)}\partial_i\partial_j( \psi u_i u_j)=\left[ \phi, \frac{\vn\partial_i\partial_j}{(-\Delta)}\right]( \psi u_i u_j)+ \frac{\vn\partial_i\partial_j}{(-\Delta)}( \phi u_i u_j),$$
and we finally obtain the following expression for (\ref{ExpressionPression6}):
\begin{eqnarray*}
 \phi \vn p&=& \sum_{i,j=1}^3\left[ \phi, \frac{\vn\partial_i\partial_j}{(-\Delta)}\right]( \psi u_i u_j)+ \sum_{i,j=1}^3\frac{\vn\partial_i\partial_j}{(-\Delta)}( \phi u_i u_j)\\
&&-\sum_{i,j=1}^3  \phi\frac{\vn}{(-\Delta)}\bigg(\partial_i((\partial_j  \psi)u_i u_j)+\partial_j((\partial_i  \psi)u_i u_j)-(\partial_i \partial_j  \psi)(u_i u_j)\bigg)\\
&&+ \phi\frac{\vn ((\Delta  \psi)p)}{(-\Delta)}-2\sum_{i=1}^3 \phi\frac{\vn (\partial_i((\partial_i \psi)p))}{(-\Delta)}.
\end{eqnarray*}
With this expression for the term that contains the pressure $p$, we obtain the (lengthy) formula for (\ref{Formule_equationV}):
\begin{align}
&\vV= \underbrace{(\partial_t  \phi - \Delta  \phi)\vu}_{(1)}-2\sum_{i=1}^3 \underbrace{(\partial_i  \phi)(\partial_i \vu)}_{(2)}-\underbrace{ \phi(\vu\cdot \vn)\vu}_{3}-\sum_{i,j=1}^3\underbrace{\left[ \phi,\frac{\vn\partial_i\partial_j}{(-\Delta)}\right]( \psi u_i u_j)}_{(4)}+\sum_{i,j=1}^3 \underbrace{\frac{\vn\partial_i\partial_j}{(-\Delta)}( \phi u_i u_j)}_{(5)}\notag\\
&- \sum_{i,j=1}^3 \frac{ \phi\vn}{(-\Delta)}\big[\underbrace{\partial_i((\partial_j  \psi)u_i u_j)}_{(6)}+\underbrace{\partial_j((\partial_i  \psi)u_i u_j)}_{(7)}-\underbrace{(\partial_i \partial_j  \psi)(u_i u_j)}_{(8)}\big]\label{FormulePourN}\\
&+2\underbrace{ \phi\frac{\vn ((\Delta  \psi)p)}{(-\Delta)}}_{(9)}-2\sum_{i=1}^3\underbrace{ \phi\frac{\vn (\partial_i((\partial_i \psi)p))}{(-\Delta)}}_{(10)}+\underbrace{\phi(\rot\vw)}_{(11)}.\notag
\end{align}
Thus, by the Duhamel formula, the solution $\vv$ of the equation  (\ref{equationV}) is given by
$$\vv=\int_0^{t}e^{(t-s)\Delta}\vV(s,\cdot)ds=\sum_{k=1}^{11}\int_0^{t}e^{(t-s)\Delta}\vV_k(s,\cdot)ds=\sum_{k=1}^{11}\vec{\mathbb{V}}_k.$$
Since $\vv=\phi\vu$, and due to the support properties of $\phi$ (see (\ref{ProprieteLocalisation_IterationVU})), we have $\mathds{1}_{Q_{R_2}}\vv=\mathds{1}_{Q_{R_2}}\vu$ and to conclude that $\mathds{1}_{Q_{R_2}}\vu \in \M_{t,x}^{3,\sigma}(\mathbb{R}\times \R) $ we will study $\mathds{1}_{Q_{R_2}}\vec{\mathbb{V}}_k$ for all $1\leq k\leq 11$. 
\begin{itemize}
\item For $\vec{\mathbb{V}}_1$, by the term (1) in (\ref{FormulePourN}) we have
\begin{equation}\label{Estimation_Ponctuelle_Riesz0}
|\mathds{1}_{Q_{R_2}}\vec{\mathbb{V}}_1(t,x)|=\left|\mathds{1}_{Q_{R_2}}\int_{0}^{t} e^{(t-s)\Delta}[(\partial_t \phi - \Delta \phi)\vu](s,x)ds\right|,
\end{equation}
since the convolution kernel of the semi-group $e^{(t-s)\Delta}$ is the  usual 3D heat kernel $\mathfrak{g}_{t}$, we can write by the decay properties of the heat kernel as well as the properties of the test function $\phi$ (see (\ref{ProprieteLocalisation_IterationVU})), the estimate
$$|\mathds{1}_{Q_{R_2}}\vec{\mathbb{V}}_1(t,x)|\leq C\mathds{1}_{Q_{R_2}}\int_{\mathbb{R}} \int_{\mathbb{R}^3} \frac{1}{(|t-s|^{\frac{1}{2}}+|x-y|)^3}\left| \mathds{1}_{Q_{\mathfrak{R}_b}}\vu(s,y)
\right| \,dy \,ds,$$
Now, recalling the definition of the parabolic Riesz potential given in (\ref{Def_ParabolicRieszPotential}) and since $Q_{R_2}\subset Q_{\mathfrak{R}_b}$ we obtain the pointwise estimate 
\begin{equation}\label{Estimation_Ponctuelle_Riesz1}
|\mathds{1}_{Q_{R_2}}\vec{\mathbb{V}}_1(t,x)|\leq C\mathds{1}_{Q_{\mathfrak{R}_b}}\mathcal{L}_{2}(|\mathds{1}_{Q_{\mathfrak{R}_b}}\vu|)(t,x),
\end{equation}
and taking Morrey $\mathcal{M}_{t,x}^{3, \sigma}$ norm we obtain
$$\|\mathds{1}_{Q_{R_2}}\vV_1(t,x)\|_{\mathcal{M}_{t,x}^{3, \sigma}}\leq C\|\mathds{1}_{Q_{\mathfrak{R}_b}}\mathcal{L}_{2}(  |\mathds{1}_{Q_{\mathfrak{R}_b}} \vu|)\|_{\mathcal{M}_{t,x}^{3, \sigma}}.$$
Now, for some $2<q<\frac{5}{2}$ we set $\lambda=1-\frac{2q}{5}$. Then, we have $ 3\le\frac{3}{\lambda}$ and $\sigma\leq \frac{q}{\lambda}$. Thus, by Lemma \ref{lemma_locindi} and by Lemma \ref{Lemma_Hed} we can write:
\begin{eqnarray*}
\|\mathds{1}_{Q_{\mathfrak{R}_b}}\mathcal{L}_{2}(  |\mathds{1}_{Q_{\mathfrak{R}_b}}  \vu |)\|_{\mathcal{M}_{t,x}^{3, \sigma}}&\leq &C\|\mathcal{L}_{2}(  |\mathds{1}_{Q_{\mathfrak{R}_b}} \vu|)\|_{\mathcal{M}_{t,x}^{\frac{3}{\lambda},\frac{q}{\lambda}}}\notag\\
&\leq &C\|\mathds{1}_{Q_{\mathfrak{R}_b}}  \vu\|_{\mathcal{M}_{t,x}^{3, q}}\leq C\|\mathds{1}_{Q_{R_1}} \vu\|_{\mathcal{M}_{t,x}^{3, \tau_0}}<+\infty,
\end{eqnarray*}
where in the last estimate we applied again Lemma \ref{lemma_locindi}  (noting that $q<\tau_0=6$) and we used the estimates over $\vu$ available in (\ref{Conclusion_FirstMorreySpace}).

\item For $\vec{\mathbb{V}}_2$, using the expression (2) in (\ref{FormulePourN}) we write $(\partial_{i}\phi) (\partial_i\vu)=\partial_i((\partial_i\phi)\vu)-(\partial_i^2\phi)\vu$ and we have
\begin{equation}\label{EstimationPonctuelleVV2}
|\mathds{1}_{Q_{R_2}}\vec{\mathbb{V}}_2(t,x)|\leq \sum_{i=1}^{3}\left|\mathds{1}_{Q_{R_2}}\int_{0}^{t} e^{ (t-s)\Delta} \partial_i\big((\partial_{i}\phi) \vu\big)ds\right|+\left|\mathds{1}_{Q_{R_2}}\int_{0}^{t} e^{ (t-s)\Delta} [(\partial_{i}^2\phi)\vu] ds\right|.
\end{equation}
Remark that the second term of the right-hand side of  (\ref{EstimationPonctuelleVV2}) can be treated in the same manner as the term $\vec{\mathbb{V}}_1$ so we will only study the first term: by the properties of the heat kernel and by the definition of the Riesz potential $\mathcal{L}_{1}$ (see (\ref{Def_ParabolicRieszPotential})), we obtain
\begin{eqnarray}
A_2:=\left|\mathds{1}_{Q_{R_2}}\int_{0}^{t} e^{ (t-s)\Delta} \partial_i\big((\partial_{i}\phi) \vu \big)ds\right|=\left|\mathds{1}_{Q_{R_2}}\int_{0}^{t} \int_{\mathbb{R}^3}\partial_i\mathfrak{g}_{t-s}(x-y)(\partial_{i}\phi) \vu(s,y)dyds\right|\notag\\
\leq C\mathds{1}_{Q_{R_2}}\int_{\mathbb{R}}\int_{\mathbb{R}^3} \frac{|\mathds{1}_{Q_{\mathfrak{R}_b}} \vu(s,y)|}{(|t-s|^{\frac{1}{2}}+|x-y|)^4}dyds
\leq C\mathds{1}_{Q_{R_2}}(\mathcal{L}_1(|\mathds{1}_{Q_{\mathfrak{R}_b}} \vu |))(t,x).\label{Estimation_Ponctuelle_Riesz2}
\end{eqnarray}
Taking the Morrey $\mathcal{M}_{t,x}^{3, \sigma}$ norm we obtain $\|A_2\|_{\mathcal{M}_{t,x}^{3, \sigma}}\leq C\|\mathds{1}_{Q_{R_2}}(\mathcal{L}_1(|\mathds{1}_{Q_{\mathfrak{R}_b}} \vu|))\|_{\mathcal{M}_{t,x}^{3, \sigma}}$.
Now, for some $4\leq q <5<\tau_0=6$ we define $\lambda=1-\frac{q}{5}$, noting that $3\leq \frac{3}{\lambda}$ and $\sigma \leq\frac{q}{\lambda}$, by Lemma \ref{Lemma_Hed}, we can write
\begin{eqnarray*}
\|\mathds{1}_{Q_{R_2}}(\mathcal{L}_1(|\mathds{1}_{Q_{\mathfrak{R}_b}} \vu|))\|_{\mathcal{M}_{t,x}^{3, \sigma}}&\leq &C\|\mathcal{L}_1(|\mathds{1}_{Q_{\mathfrak{R}_b}} \vu|)\|_{\mathcal{M}_{t,x}^{\frac{3}{\lambda}, \frac{q}{\lambda}}}\leq C\|\mathds{1}_{Q_{\mathfrak{R}_b}} \vu\|_{\mathcal{M}_{t,x}^{3, q}}\\
&\leq & C\|\mathds{1}_{Q_{R_1}} \vu\|_{\mathcal{M}_{t,x}^{3, \tau_0}}<+\infty,
\end{eqnarray*}
from which we deduce that $\|\mathds{1}_{Q_{R_2}}\vec{\mathbb{V}}_2\|_{\mathcal{M}_{t,x}^{3, \sigma}}<+\infty$.

\item For the term $\vec{\mathbb{V}}_3$, by the same arguments given to obtain the pointwise estimate (\ref{Estimation_Ponctuelle_Riesz1}), we have
\begin{eqnarray*}
|\mathds{1}_{Q_{R_2}}\vec{\mathbb{V}}_3(t,x)|&=&\left|\mathds{1}_{Q_{R_2}}\int_{0}^{t} \int_{\mathbb{R}^3}  \mathfrak{g}_{t-s}(x-y)\left[\phi \left( (\vu\cdot\vn)\vu\right)\right](s,y)dyds\right|\\
&\leq & C\mathds{1}_{Q_{R_2}}\mathcal{L}_2\left(\left|\mathds{1}_{Q_{\mathfrak{R}_b}}\left( (\vu\cdot\vn)\vu\right)\right|\right)(t,x),
\end{eqnarray*}
(recall (\ref{ProprieteLocalisation_IterationVU})) from which we deduce 
\begin{equation}\label{Decomposition2termesVV3}
\|\mathds{1}_{Q_{R_2}}\vec{\mathbb{V}}_3\|_{\mathcal{M}_{t,x}^{3, \sigma}}\leq C\left\|\mathds{1}_{Q_{R_2}}\mathcal{L}_2\left(|\mathds{1}_{Q_{\mathfrak{R}_b}} (\vu\cdot\vn)\vu|\right)\right\|_{\mathcal{M}_{t,x}^{3, \sigma}}.
\end{equation}
We set now $\frac{5}{3-\alpha}<q<\frac{5}{2}$ and $\lambda=1-\frac{2q}{5}$. Since $3\leq \frac{6}{5\lambda}$ and $\tau_0=6<\sigma\leq \frac{q}{\lambda}$, applying Lemma \ref{lemma_locindi} and Lemma \ref{Lemma_Hed} we have
$$\left\|\mathds{1}_{Q_{R_2}}\mathcal{L}_2\left(|\mathds{1}_{Q_{\mathfrak{R}_b}} (\vu\cdot\vn)\vu|\right)\right\|_{\mathcal{M}_{t,x}^{3, \sigma}}\leq C\left\|\mathds{1}_{Q_{R_2}}\mathcal{L}_2\left(|\mathds{1}_{Q_{\mathfrak{R}_b}} (\vu\cdot\vn)\vu|\right)\right\|_{\mathcal{M}_{t,x}^{\frac{6}{5\lambda},\frac{q}{\lambda}}}\leq C\left\|\mathds{1}_{Q_{\mathfrak{R}_b}} (\vu\cdot\vn)\vu\right\|_{\mathcal{M}_{t,x}^{\frac{6}{5}, q}}.$$
Recall that we have $\tau_0=6<\sigma$ and by the H\"older inequality in Morrey spaces (see Lemma \ref{lemma_Product}) we obtain
$$\left\|\mathds{1}_{Q_{\mathfrak{R}_b}} (\vu\cdot\vn)\vu\right\|_{\mathcal{M}_{t,x}^{\frac{6}{5}, q}}\leq \left\|\mathds{1}_{Q_{R_1}}\vu\right\|_{\mathcal{M}_{t,x}^{3, \tau_0}}\left\|\mathds{1}_{Q_{R_1}}\vn \otimes \vu\right\|_{\mathcal{M}_{t,x}^{2, \tau_1}}<+\infty,$$
where $\frac{1}{q}=\frac{1}{\tau_0}+\frac{1}{\tau_1}=\frac{2}{\tau_0}+\frac{1}{5}$. These two last quantities are bounded by (\ref{Conclusion_FirstMorreySpace}) and (\ref{ConlusioncorolarioMorrey}). Note that the condition $\tau_0=6<\sigma$ and the relationship $\frac{1}{q}=\frac{2}{\tau_0}+\frac{1}{5}$ are compatible with the fact that $\frac{5}{3-\alpha}<q < \frac{5}{2}$ (recall that $0<\alpha<\frac{1}{24}$). 
\item The term $\vec{\mathbb{V}}_4$ is the most technical one. Indeed, by the expression of $\vV_4$ given in (\ref{FormulePourN}), we write
\begin{eqnarray*}
|\mathds{1}_{Q_{R_2}}\vec{\mathbb{V}}_4|\leq \sum^3_{i,j= 1} \mathds{1}_{Q_{R_2}}\int_{\mathbb{R}} \int_{\mathbb{R}^3}\frac{\left|\left[\phi, \, \frac{ \vn \partial_i\partial_j}{(-\Delta)}\right] (\psi u_i u_j)(s,y)\right|}{(|t-s|^{\frac{1}{2}}+|x-y|)^3}dyds\leq \sum^3_{i,j= 1}\mathds{1}_{Q_{R_2}}\mathcal{L}_{2}\left(\left|\left[\phi, \, \frac{ \vn \partial_i\partial_j}{(-\Delta)}\right] (\psi u_i u_j)\right|\right),
\end{eqnarray*}
and taking the $\mathcal{M}_{t,x}^{3, \sigma}$-norm we have
$\|\mathds{1}_{Q_{R_2}}\vec{\mathbb{V}}_4\|_{\mathcal{M}_{t,x}^{3, \sigma}}\leq \sum^3_{i,j= 1}\left\|\mathds{1}_{Q_{R_2}}\mathcal{L}_{2}\left(\left|\left[\phi, \, \frac{ \vn \partial_i\partial_j}{(-\Delta)}\right] (\psi u_i u_j)\right|\right)\right\|_{\mathcal{M}_{t,x}^{3, \sigma}}$. If we set $\frac{1}{q}=  \frac{2}{\tau_0} + \frac{1}{5}$ and $\lambda=1-\frac{2q}{5}$ then we have $3\leq \frac{3}{2\lambda}$ and for
\begin{equation}\label{UpperBound_sigma}
\sigma \leq \frac{q}{\lambda} = \frac{5 \tau_0}{10- \tau_0},
\end{equation}
by Lemmas \ref{lemma_locindi} and \ref{Lemma_Hed} we obtain:
\begin{eqnarray*}
\left\|\mathds{1}_{Q_{R_2}}\mathcal{L}_{2}\left(\left|\left[\phi, \, \frac{ \vn \partial_i\partial_j}{(-\Delta)}\right] (\psi u_i u_j)\right|\right)\right\|_{\mathcal{M}_{t,x}^{3, \sigma}}&\leq &C\left\|\mathds{1}_{Q_{R_2}}\mathcal{L}_{2}\left(\left|\left[\phi, \, \frac{ \vn \partial_i\partial_j}{(-\Delta)}\right] (\psi u_i u_j)\right|\right)\right\|_{\mathcal{M}_{t,x}^{\frac{3}{2\lambda}, \frac{q}{\lambda}}}\\
&\leq &C\left\|\left[\phi, \, \frac{ \vn \partial_i\partial_j}{(-\Delta)}\right] (\psi u_i u_j)\right\|_{\mathcal{M}_{t,x}^{\frac{3}{2}, q}},
\end{eqnarray*}
We will study this norm and by the definition of Morrey spaces (\ref{DefMorreyparabolico}), if we introduce a threshold $\mathfrak{r}=\frac{\mathfrak{R}_b-R_2}{2}$, we have
\begin{equation}\label{CommutatorEstimatevVV4}
\begin{split}
\left\|\left[\phi, \, \frac{ \vn \partial_i\partial_j}{(-\Delta)}\right] (\psi u_i u_j)\right\|_{\mathcal{M}_{t,x}^{\frac{3}{2}, q}}^{\frac{3}{2}}&\leq\underset{\underset{0<r<\mathfrak{r}}{(\mathfrak{t},\bar{x})}}{\sup}\;\frac{1}{r^{5(1-\frac{3}{2q})}}\int_{Q_r(\mathfrak{t},\bar x)}\left|\left[\phi, \, \frac{ \vn \partial_i\partial_j}{(-\Delta)}\right] (\psi u_i u_j)\right|^{\frac{3}{2}}dxdt\qquad\\
&+\underset{\underset{\mathfrak{r}< r}{(\mathfrak{t},\bar{x})}}{\sup}\;\frac{1}{r^{5(1-\frac{3}{2q})}}\int_{Q_r(\mathfrak{t},\bar{x})}\left|\left[\phi, \, \frac{ \vn \partial_i\partial_j}{(-\Delta)}\right] (\psi u_i u_j)\right|^{\frac{3}{2}}dxdt.\qquad
\end{split}
\end{equation}
Now, we study the second term of the right-hand side above, which is easy to handle as we have $\mathfrak{r}<r$ and we can write
$$\underset{\underset{\mathfrak{r}< r}{(\mathfrak{t},\bar{x}) \in \mathbb{R}\times \R}}{\sup}\;\frac{1}{r^{5(1-\frac{3}{2q})}}\int_{Q_r(\mathfrak{t},\bar{x})}\left|\left[\phi, \, \frac{ \vn \partial_i\partial_j}{(-\Delta)}\right] (\psi u_i u_j)\right|^{\frac{3}{2}}dxdt\leq C_{\mathfrak{r}}\left\|\left[\phi, \, \frac{ \vn \partial_i\partial_j}{(-\Delta)}\right] (\psi u_i u_j)\right\|_{L^{\frac{3}{2}}_{t,x}}^{\frac{3}{2}},$$
and since $\bar\phi$ is a regular function and $\frac{ \vn \partial_i\partial_j}{(-\Delta)}$ is a Calder\'on-Zydmund operator, by the Calder\'on commutator theorem (see the book \cite{PGLR0}), we have that the operator $\left[\phi, \, \frac{ \vn \partial_i\partial_j}{(-\Delta)}\right]$ is bounded in the space $L^{\frac{3}{2}}_{t,x}$ and we can write (using the support properties of $\psi$ given in (\ref{ProprieteLocalisation_IterationVU}) and the information given in (\ref{Conclusion_FirstMorreySpace})):
\begin{eqnarray*}
\left\|\left[\bar\phi, \, \frac{ \vn \partial_i\partial_j}{(-\Delta)}\right] (\psi u_i u_j)\right\|_{L^{\frac{3}{2}}_{t,x}}&\leq &C\left\|\psi u_i u_j\right\|_{L^{\frac{3}{2}}_{t,x}}\leq C\|\mathds{1}_{Q_{R_1}} u_i u_j\|_{\mathcal{M}^{\frac{3}{2}, \frac{3}{2}}_{t,x}}\\
&\leq & C\|\mathds{1}_{Q_{R_1}} \vu\|_{\mathcal{M}^{3,3}_{t,x}}\|\mathds{1}_{Q_{R_1}} \vu\|_{\mathcal{M}^{3, 3}_{t,x}}\leq C\|\mathds{1}_{Q_{R_1}} \vu\|_{\mathcal{M}^{3,\tau_0}_{t,x}}\|\mathds{1}_{Q_{R_1}} \vu\|_{\mathcal{M}^{3, \tau_0}_{t,x}}<+\infty,
\end{eqnarray*}
where in the last line we used H\"older inequalities in Morrey spaces and we applied Lemma \ref{lemma_locindi}.\\

The first term of the right-hand side of  (\ref{CommutatorEstimatevVV4}) requires some extra computations: indeed, as we are interested to obtain information over the parabolic ball $Q_{r}(\mathfrak{t}, \bar{x})$ we can write
for some  $0<r<\mathfrak{r}$:
\begin{equation}\label{CommutatorEstimatevVV401}
\mathds{1}_{Q_{r}}\left[\phi, \, \frac{ \vn \partial_i\partial_j}{(-\Delta)}\right] (\psi u_i u_j))=\mathds{1}_{Q_{r}}\left[\phi, \, \frac{ \vn \partial_i\partial_j}{(-\Delta)}\right] (\mathds{1}_{Q_{2r}}\psi u_i u_j)+\mathds{1}_{Q_{r}}\left[\phi, \, \frac{ \vn \partial_i\partial_j}{(-\Delta)}\right] ((\mathbb{I}-\mathds{1}_{Q_{2r}})\psi u_i u_j),
\end{equation}
and as before we will study the $L^{\frac{3}{2}}_{t,x}$ norm of these two terms. For the first quantity in the right-hand side of (\ref{CommutatorEstimatevVV401}), by the Calder\'on commutator theorem, by the definition of Morrey spaces and by the H\"older inequalities we have 
\begin{eqnarray*}
\left\|\mathds{1}_{Q_{r}}\left[\phi, \, \frac{ \vn \partial_i\partial_j}{(-\Delta)}\right] (\mathds{1}_{Q_{2r}}\psi u_i u_j)\right\|_{L^{\frac{3}{2}}_{t,x}}^{\frac{3}{2}}&\leq& C\|\mathds{1}_{Q_{2r}}\psi u_i u_j\|_{L^{\frac{3}{2}}_{t,x}}^{\frac{3}{2}}\leq Cr^{5 (1-\frac{3}{\tau_0})} \|\mathds{1}_{Q_{R_1}}u_i u_j\|_{\mathcal{M}^{\frac{3}{2}, \frac{\tau_0}{2}}_{t,x}}^{\frac{3}{2}}\\
&\leq & Cr^{5 (1-\frac{3}{\tau_0})} \|\mathds{1}_{Q_{R_1}}\vu\|_{\mathcal{M}^{3, \tau_0}_{t,x}}^{\frac{3}{2}}\|\mathds{1}_{Q_{R_1}}\vu\|_{\mathcal{M}^{3, \tau_0}_{t,x}}^{\frac{3}{2}},
\end{eqnarray*}
for all $0<r<\mathfrak{r}$, from which we deduce that 
$$\underset{\underset{0<r<\mathfrak{r}}{(\mathfrak{t},\bar{x}) }}{\sup}\;\frac{1}{r^{5(1-\frac{3}{2q})}}\int_{Q_r(\mathfrak{t},\bar{x})}\left|\mathds{1}_{Q_{r}}\left[\phi, \, \frac{ \vn \partial_i\partial_j}{(-\Delta)}\right] (\mathds{1}_{Q_{2r}}\psi u_i u_j)\right|^{\frac{3}{2}}dxdt\leq C \|\mathds{1}_{Q_{R_1}}\vu\|_{\mathcal{M}^{3, \tau_0}_{t,x}}^{\frac{3}{2}}\|\mathds{1}_{Q_{R_1}}\vu\|_{\mathcal{M}^{3, \tau_0}_{t,x}}^{\frac{3}{2}}<+\infty.$$
We study now the second term of the right-hand side of (\ref{CommutatorEstimatevVV401}) and for this we consider the following operator:
$$T: f \mapsto  \left(\mathds{1}_{Q_{r}}  \left[\phi, \, \frac{ \vn \partial_i \partial_j}{- \Delta }\right] (\mathbb{I} - \mathds{1}_{Q_{2r}})  \psi
\right) f,$$
and by the properties of the convolution kernel of the operator $\frac{1}{(-\Delta)}$ we obtain
$$|T(f)(x)|\leq C\mathds{1}_{Q_{r}}(x)\int_{\mathbb{R}^3}\frac{(\mathbb{I} - \mathds{1}_{Q_{2r}})(y) \mathds{1}_{Q_{R_1}}(y) |f(y)| | \phi(x)- \phi(y)|}{|x-y|^4} dy.$$
Recalling that $0<r<\mathfrak{r}=\frac{\mathfrak{R}_b-R_2}{2}$, by the support properties of the test function $\phi$ (see (\ref{ProprieteLocalisation_IterationVU})), the integral above is meaningful if $|x-y|>r$ and thus we can write
\begin{eqnarray*}
\left\|\mathds{1}_{Q_{r}}\left[\phi, \, \frac{ \vn \partial_i\partial_j}{(-\Delta)}\right] ((\mathbb{I}-\mathds{1}_{Q_{2r}})\psi u_i u_j)\right\|_{L^\frac{3}{2}_{t,x}}^\frac{3}{2}\leq C\left\|\mathds{1}_{Q_{r}} \int_{\mathbb{R}^3} \frac{\mathds{1}_{|x-y| > r}}{|x-y|^4}(\mathbb{I} - \mathds{1}_{Q_{2r}})(y) \mathds{1}_{Q_{R_1}}(y)|u_i u_j |dy\right\|_{L^\frac{3}{2}_{t,x}}^\frac{3}{2}\\
\leq C\left(\int_{|y|>r}\frac{1}{|y|^4}\| \mathds{1}_{Q_{R_1}}
|u_i u_j |(\cdot-y)\|_{L^\frac{3}{2}_{t,x}(Q_r)}dy\right)^\frac{3}{2}\leq Cr^{-\frac{3}{2}}\| \mathds{1}_{Q_{R_1}}u_i u_j \|_{L^\frac{3}{2}_{t,x}(Q_{r})}^\frac{3}{2},
\end{eqnarray*}
with this estimate at hand and using the definition of Morrey spaces, we can write
\begin{eqnarray*}
\int_{Q_r(\mathfrak{t},\bar{x})}\left|\mathds{1}_{Q_{r}}\left[\phi, \, \frac{ \vn \partial_i\partial_j}{(-\Delta)}\right] ((\mathbb{I}-\mathds{1}_{Q_{2r}})\psi u_i u_j)\right|^{\frac{3}{2}}dxdt&\leq &Cr^{-\frac{3}{2}}r^{5 (1-\frac{3}{\tau_0})}\| \mathds{1}_{Q_{R_1}}u_i u_j \|_{\mathcal{M}^{\frac{3}{2}, \frac{\tau_0}{2}}_{t,x}}^\frac{3}{2}\\
&\leq &Cr^{5(1-\frac{3}{2q})}\| \mathds{1}_{Q_{R_1}}u_i u_j \|_{\mathcal{M}^{\frac{3}{2}, \frac{\tau_0}{2}}_{t,x}}^\frac{3}{2},
\end{eqnarray*}
where in the last inequality we used the fact that $\frac{1}{q}=  \frac{2}{\tau_0} + \frac{1}{5}$, which implies  $r^{-\frac{3}{2}}r^{5 (1-\frac{3}{\tau_0})}= r^{5(1-\frac{3}{2q})}$. Thus we finally obtain
$$\underset{\underset{0<r<\mathfrak{r}}{(\mathfrak{t},\bar{x}) }}{\sup}\;\frac{1}{r^{5(1-\frac{3}{2q})}}\int_{Q_r(\mathfrak{t},\bar{x})}\left|\mathds{1}_{Q_{r}}\left[\phi, \, \frac{ \vn \partial_i\partial_j}{(-\Delta)}\right] ((\mathbb{I}-\mathds{1}_{Q_{2r}})\psi u_i u_j)\right|^{\frac{3}{2}}dxdt\leq C \|\mathds{1}_{Q_{R_1}}\vu\|_{\mathcal{M}^{3, \tau_0}_{t,x}}^{\frac{3}{2}}\|\mathds{1}_{Q_{R_1}}\vu\|_{\mathcal{M}^{3, \tau_0}_{t,x}}^{\frac{3}{2}}<+\infty.$$
We have proven that all the term in (\ref{CommutatorEstimatevVV4}) are bounded and we can conclude that $\|\mathds{1}_{Q_{R_2}}\vec{\mathbb{V}}_4\|_{\mathcal{M}_{t,x}^{3, \sigma}}<+\infty$.
\begin{Remark}
The condition \eqref{UpperBound_sigma} implies an upper bound for $\sigma$ depending on the current Morrey information of $\vu$, which a priori is close to $\tau_0=6$. Nevertheless it is clear that whether we obtain a better Morrey information on integrability for $\vu$, the value of $\sigma$ can increase.
\end{Remark}
\item For the quantity $\vec{\mathbb{V}}_5$, based in the expression (\ref{FormulePourN}) we write
\begin{eqnarray*}
|\mathds{1}_{Q_{R_2}}\vec{\mathbb{V}}_5(t,x)|&\leq &C\sum^3_{i,j= 1} \mathds{1}_{Q_{R_2}}\int_{\mathbb{R}}\int_{\mathbb{R}^3} \frac{|\mathcal{R}_i\mathcal{R}_j( \phi u_i u_j )(s,y)|}{(|t-s|^{\frac{1}{2}}+|x-y|)^4}dyds\leq C\sum^3_{i,j= 1} \mathds{1}_{Q_{R_2}} \mathcal{L}_{1}\left(|\mathcal{R}_i\mathcal{R}_j( \phi u_i u_j )|\right)(t,x),
\end{eqnarray*}
where we used the decaying properties of the heat kernel (recall that $\mathcal{R}_i=\frac{\partial_i}{\sqrt{- \Delta}}$ are the Riesz transforms). 
Now taking the Morrey $\mathcal{M}^{3, \sigma}_{t,x}$ norm and by Lemma \ref{lemma_locindi} (with $\nu=\frac{4\tau_0+5}{5\tau_0}$, $p=3$, $q=\tau_0$ such that $\frac{p}{\nu}>3$ and $\frac{q}{\nu}>\sigma$ which is compatible with the condition $\tau_0<\sigma$) we have
\begin{eqnarray*}
\|\mathds{1}_{Q_{R_2}}\vV_5\|_{\mathcal{M}^{3, \sigma}_{t,x}}&\leq &C \sum^3_{i,j= 1}\|  \mathds{1}_{Q_{R_2}}\mathcal{L}_{1}\left(|\mathcal{R}_i\mathcal{R}_j( \phi u_i u_j )|\right)\|_{\mathcal{M}^{\frac{p}{\nu}, \frac{q}{\nu}}_{t,x}}
\end{eqnarray*}
Then by Lemma \ref{Lemma_Hed} with $\lambda= 1 - \tfrac{\tau_0 /2}{5}$ (recall $\tau_0=6 <10$ so that $\nu > 2 \lambda$) and by the boundedness of Riesz transforms in Morrey spaces we obtain:
\begin{eqnarray*}
\|\mathds{1}_{Q_{R_2}} \mathcal{L}_{1}\left(|\mathcal{R}_i\mathcal{R}_j( \phi u_i u_j )|\right)\|_{\mathcal{M}^{\frac{p}{\nu}, \frac{q}{\nu}}_{t,x}}
&\leq&C\|\mathcal{L}_{1}\left(|\mathcal{R}_i\mathcal{R}_j( \phi u_i u_j )|\right)\|_{\mathcal{M}^{\frac{p}{2\lambda}, \frac{q}{2\lambda}}_{t,x}}
\leq C \|\mathcal{R}_i\mathcal{R}_j( \phi u_i u_j )\|_{\mathcal{M}^{\frac{3}{2},\frac{\tau_0}{2}}_{t,x}}\\
&\leq &\| \mathds{1}_{Q_{R_1}} u_i u_j \|_{\mathcal{M}^{\frac{3}{2},\frac{\tau_0}{2}}_{t,x}}
\leq  C\|\mathds{1}_{Q_{R_1}} \vu\|_{\mathcal{M}^{3,\tau_0}_{t,x}}\|  \mathds{1}_{Q_{R_1}} \vu\|_{\mathcal{M}^{3,\tau_0}_{t,x}}<+\infty,
\end{eqnarray*}
and we obtain $\|\mathds{1}_{Q_{R_2}}\vec{\mathbb{V}}_5\|_{\mathcal{M}_{t,x}^{3, \sigma}}<+\infty$.
\item For the term $\vec{\mathbb{V}}_6$ and following the same ideas we have
$$|\mathds{1}_{Q_{R_2}}\vec{\mathbb{V}}_6|\leq C\sum^3_{i,j= 1}  \mathds{1}_{Q_{R_2}}\int_{\mathbb{R}}  \int_{\mathbb{R}^3}\frac{\left|\frac{\phi \vn\partial_i}{(- \Delta )} (\partial_j \psi) u_i u_j(s,y)\right|}{(|t-s|^{\frac{1}{2}}+|x-y|)^3} dyds=C\sum^3_{i,j= 1}  \mathds{1}_{Q_{R_2}}\mathcal{L}_{2}\left(\left|\frac{\phi \vn\partial_i}{(- \Delta )} (\partial_j \psi) u_i u_j\right|\right).$$
For $2<q<\frac{5}{2}$, define $\lambda=1-\frac{2q}{5}$, we thus have $3\leq \frac{3}{2\lambda}$ and $\sigma \leq \frac{q}{\lambda}$. Then, by Lemma \ref{lemma_locindi} and Lemma \ref{Lemma_Hed} we can write
\begin{eqnarray*}
\left\|\mathds{1}_{Q_{R_2}}\mathcal{L}_{2}\left|\frac{\phi \vn\partial_i}{(- \Delta )} (\partial_j \psi) u_i u_j\right|\right\|_{\mathcal{M}^{3, \sigma}_{t,x}}\leq C\left\|\mathds{1}_{Q_{R_2}}\mathcal{L}_{2}\left|\frac{\phi \vn\partial_i}{(- \Delta )} (\partial_j \psi) u_i u_j\right|\right\|_{\mathcal{M}^{\frac{3}{2\lambda}, \frac{q}{\lambda}}_{t,x}}\leq  C\left\|\frac{\phi \vn\partial_i}{(- \Delta )} (\partial_j \psi) u_i u_j\right\|_{\mathcal{M}^{\frac{3}{2}, q}_{t,x}},
\end{eqnarray*}
but since the operator $\frac{\phi \vn\partial_i}{(- \Delta )}$ is bounded in Morrey spaces and since $2<q<\frac{5}{2}< \tfrac{\tau_0}{2}=3$ (since $\tau_0=6$), one has by Lemma \ref{lemma_locindi} and by the H\"older inequalities
\begin{eqnarray*}
\left\|\frac{\phi \vn\partial_i}{(- \Delta )} (\partial_j \psi) u_i u_j\right\|_{\mathcal{M}^{\frac{3}{2}, q}_{t,x}}\leq C\left\| \mathds{1}_{Q_{R_1}}u_i u_j\right\|_{\mathcal{M}^{\frac{3}{2}, q}_{t,x}}\leq C\| \mathds{1}_{Q_{R_1}}u_i u_j\|_{\mathcal{M}^{\frac{3}{2}, \frac{\tau_0}{2}}_{t,x}}\leq C\|\mathds{1}_{Q_{R_1}} \vu\|_{\mathcal{M}^{3,\tau_0}_{t,x}}\|  \mathds{1}_{Q_{R_1}} \vu\|_{\mathcal{M}^{3,\tau_0}_{t,x}},
\end{eqnarray*}
from which we deduce $\|\mathds{1}_{Q_{R_2}}\vec{\mathbb{V}}_6\|_{\mathcal{M}^{3, \sigma}_{t,x}}<+\infty$. Note that the same computations can be performed to obtain that $\|\mathds{1}_{Q_{R_2}}\vec{\mathbb{V}}_7\|_{\mathcal{M}^{3, \sigma}_{t,x}}<+\infty$.

\item The quantity $\vec{\mathbb{V}}_8$ based in the term (8) of (\ref{FormulePourN}) is treated in the following manner: we first write
$$\|\mathds{1}_{Q_{R_2}}\vec{\mathbb{V}}_8\|_{\mathcal{M}^{3, \sigma}_{t,x}}\leq C \sum^3_{i,j= 1}\left\| \mathds{1}_{Q_{R_2}}\left(\mathcal{L}_2\left|\phi \frac{\vn}{(-\Delta )}(\partial_i \partial_j \psi) (u_i u_j)\right|\right)\right\|_{\mathcal{M}^{3,\sigma}_{t,x}}.$$
We set $1<\nu<\frac{3}{2}$, $2\nu <q<\frac{5\nu}{2}$ and $\lambda=1-\frac{2q}{5\nu}$, thus we have $3\leq \frac{\nu}{\lambda}$ and $\sigma \leq \frac{q}{\lambda}$, then, by Lemma \ref{lemma_locindi} and by Lemma \ref{Lemma_Hed} we can write
\begin{eqnarray}
\left\|\mathds{1}_{Q_{R_2}}\left(\mathcal{L}_2\left|\phi \frac{\vn}{(-\Delta )}(\partial_i \partial_j \psi) (u_i u_j)\right|\right)\right\|_{\mathcal{M}^{3,\sigma}_{t,x}}\leq C\left\| \mathds{1}_{Q_{R_2}}\left(\mathcal{L}_2\left|\phi \frac{\vn}{(-\Delta )}(\partial_i \partial_j \psi) (u_i u_j)\right|\right)\right\|_{\mathcal{M}^{\frac{\nu}{\lambda},\frac{q}{\lambda}}_{t,x}}\notag\\
\leq C\left\| \phi \frac{\vn}{(-\Delta )}(\partial_i \partial_j \psi) (u_i u_j)\right\|_{\mathcal{M}^{\nu,q}_{t,x}}\leq C\left\| \phi \frac{\vn}{(-\Delta )}(\partial_i \partial_j \psi) (u_i u_j)\right\|_{\mathcal{M}^{\nu,\frac{5\nu}{2}}_{t,x}}\notag\\
\leq C\left\| \phi \frac{\vn}{(-\Delta )}(\partial_i \partial_j \psi) (u_i u_j)\right\|_{L^{\nu}_tL^{\infty}_x}\label{Formula_intermediairevVV80}
\end{eqnarray}
where in the last estimate we used the space inclusion $L^{\nu}_tL^{\infty}_x\subset \mathcal{M}^{\nu,\frac{5\nu}{2}}_{t,x}$. 
\begin{Remark}\label{Remarque_Iteration}
Note that if the parameter $q$ above is close to the value $\frac{5\nu}{2}$, then $\lambda=1-\frac{2q}{5\nu}$ is close to $0$ and thus the value $\frac{q}{\lambda}$ can be made very big: in the estimates (\ref{Formula_intermediairevVV80}) we can consider a Morrey space $\mathcal{M}^{3,\sigma}_{t,x}$ with $\sigma\gg 1$.
\end{Remark}
Let us focus now in the $L^\infty$ norm above (\emph{i.e.} without considering the time variable). Remark that due to the support properties of the auxiliary function $\psi$ given in (\ref{ProprieteLocalisation_IterationVU}) we have $supp(\partial_i \partial_j \psi) \subset Q_{R_1}\setminus Q_{\mathfrak{R}_a}$ and recall by (\ref{ProprieteLocalisation_IterationVU}) we have $supp\;  \phi = Q_{\mathfrak{R}_b}$ where $\mathfrak{R}_b<\mathfrak{R}_a<R_1$, thus by the properties of the kernel of the operator $\frac{\vn}{(-\Delta)}$ we can write
\begin{eqnarray}
\left| \phi \frac{\vn}{(-\Delta )}(\partial_i \partial_j \psi) (u_i u_j)\right|\leq C\left|\int_{\mathbb{R}^3} \frac{1}{|x-y|^2}\mathds{1}_{Q_{\mathfrak{R}_b}}(x)\mathds{1}_{Q_{R_1}\setminus Q_{\mathfrak{R}_a}}(y)(\partial_i \partial_j \psi) (u_i u_j)(\cdot,y)dy\right|\notag\\
\leq  C\left|\int_{\mathbb{R}^3} \frac{\mathds{1}_{|x-y|>(\mathfrak{R}_a-\mathfrak{R}_b)}}{|x-y|^2}\mathds{1}_{Q_{\mathfrak{R}_b}}(x)\mathds{1}_{Q_{R_1}\setminus Q_{\mathfrak{R}_a}}(y)(\partial_i \partial_j \psi) (u_i u_j)(\cdot,y)dy\right|,\label{Formula_intermediairevVV801}
\end{eqnarray}
and the previous expression is nothing but the convolution between the function $(\partial_i \partial_j \psi) (u_i u_j)$ and a $L^\infty$-function, thus we have 
\begin{equation}\label{Formula_intermediairevVV81}
\left\| \phi \frac{\vn}{(-\Delta )}(\partial_i \partial_j \psi) (u_i u_j) (t,\cdot)\right\|_{L^\infty}\leq C\|(\partial_i \partial_j \psi) (u_i u_j)(t,\cdot)\|_{L^1}\leq C\|\mathds{1}_{Q_{R_1}}(u_i u_j)(t,\cdot)\|_{L^{\nu}},
\end{equation}
and taking the $L^\nu$-norm in the time variable we obtain
\begin{eqnarray*}
\left\| \phi \frac{\vn}{(-\Delta )}(\partial_i \partial_j \psi) (u_i u_j)\right\|_{L^{\nu}_t L^{\infty}_{x}}
&\leq &C\|\mathds{1}_{Q_{R_1}}u_i u_j\|_{L^{\nu}_{t,x}}\leq C\|\mathds{1}_{Q_{R_1}}\vu\|_{\mathcal{M}^{3,\tau_0}_{t,x}}\|\mathds{1}_{Q_{R_1}}\vu\|_{\mathcal{M}^{3,\tau_0}_{t,x}}<+\infty,
\end{eqnarray*}
where we used the fact that $1<\nu<\frac{3}{2}<\frac{\tau_0}{2}$ and we applied Hölder's inequality. Gathering together all these estimates we obtain 
$\|\mathds{1}_{Q_{R_2}}\vec{\mathbb{V}}_8\|_{\mathcal{M}^{3, \sigma}_{t,x}}<+\infty$.\\

\item The quantity $\vec{\mathbb{V}}_{9}$ based in the term (9) of (\ref{FormulePourN})  can be treated in a similar manner. Indeed, by the same arguments displayed to deduce (\ref{Formula_intermediairevVV80}), we can write (recall that $1<\nu<\frac{3}{2}$): 
$$\|\mathds{1}_{Q_{R_2}}\vec{\mathbb{V}}_{9}\|_{\mathcal{M}^{3, \sigma}_{t,x}}\leq C\left\| \phi \frac{\vn}{(-\Delta )}( (\Delta \psi) p)\right\|_{L^{\nu}_t L^{\infty}_{x}},$$ 
and if we study the $L^\infty$-norm in the space variable of this term, by the same ideas used in (\ref{Formula_intermediairevVV801})-(\ref{Formula_intermediairevVV81}) we obtain $\left\| \phi \frac{\vn}{(-\Delta )}( (\Delta \psi) p) (t,\cdot)\right\|_{L^{\infty}}\leq C\|(\Delta \psi) p (t,\cdot)\|_{L^1}\leq C\|\mathds{1}_{Q_{R_1}}p (t,\cdot)\|_{L^{\nu}}$. Thus, taking the $L^{\nu}$-norm in the time variable we have
$$\|\mathds{1}_{Q_{R_2}}\vec{\mathbb{V}}_{9}\|_{\mathcal{M}^{3, \sigma}_{t,x}}
\leq C\left\| \phi \frac{\vn}{(-\Delta )}( (\Delta \psi) p)\right\|_{L^{\nu}_t L^{\infty}_{x}}\leq C \|\mathds{1}_{Q_{R_1}}p\|_{L^{\nu}_{t,x}}\leq C \|\mathds{1}_{Q_{R_1}}p\|_{L^{\frac32}_{t,x}}<+\infty.$$

\item The study of the quantity $\vec{\mathbb{V}}_{10}$ follows almost the same lines as the terms $\vec{\mathbb{V}}_{8}$ and $\vec{\mathbb{V}}_{9}$. However instead of (\ref{Formula_intermediairevVV801}) we have
$$\left|\phi \frac{\vn\partial_i }{(-\Delta)} ( (\partial_i \psi)p)\right|
\leq  C\left|\int_{\mathbb{R}^3} \frac{\mathds{1}_{|x-y|>(\mathfrak{R}_a-\mathfrak{R}_b)}}{|x-y|^3}\mathds{1}_{Q_{\mathfrak{R}_b}}(x)\mathds{1}_{Q_{R_1}\setminus Q_{\mathfrak{R}_a}}(y)(\partial_i  \psi) p(t,y)dy\right|,$$
and thus we can write:
$$\|\mathds{1}_{Q_{R_2}}\vec{\mathbb{V}}_{10}\|_{\mathcal{M}^{3, \sigma}_{t,x}}\leq  \left\|\phi \frac{\vn\partial_i }{(-\Delta)} ( (\partial_i \psi) p)\right\|_{L^{\nu}_t L^{\infty}_{x}}\leq C \|\mathds{1}_{Q_{R_1}}p\|_{L^{\nu}_{t,x}}\leq C \|\mathds{1}_{Q_{R_1}}p\|_{L^{\frac32}_{t,x}}<+\infty.$$

Note that, by the same reason given in the Remark \ref{Remarque_Iteration}, in the study of the terms that contain the pressure (\emph{i.e.} $\vec{\mathbb{V}}_{9}$ and $\vec{\mathbb{V}}_{10}$) we can consider a Morrey space $\mathcal{M}^{3,\sigma}_{t,x}$ with $\sigma\gg 1$. But this is not the case anymore for the last term below.

\item Finally, for the term $\vec{\mathbb{V}}_{11}$ based in the term (11) of (\ref{FormulePourN}) we write:
\begin{eqnarray*}
|\mathds{1}_{Q_{R_2}}\vec{\mathbb{V}}_{11}|&=&\left|\mathds{1}_{Q_{R_2}}\int_{0}^{t} e^{(t-s)\Delta}[\phi(\rot\vw)](s,x)ds\right|\\
&\leq &\left|\mathds{1}_{Q_{R_2}}\int_{0}^{t} e^{(t-s)\Delta}\rot(\phi \vw)(s,x)ds\right|+\left|\mathds{1}_{Q_{R_2}}\int_{0}^{t} e^{(t-s)\Delta}(\vn\phi)\wedge\vw(s,x)ds\right|\\
&\leq & \mathbb{V}_{a}+\mathbb{V}_{b},
\end{eqnarray*}
for the first term above, and following the ideas given in (\ref{Estimation_Ponctuelle_Riesz2}), we have the following estimate with the Riesz potential $\mathcal{L}_1$, and by Lemma \ref{lemma_locindi} we can write
$$\| \mathbb{V}_{a}\|_{\mathcal{M}^{3, \sigma}_{t,x}}\leq C\|\mathds{1}_{Q_{R_2}}(\mathcal{L}_1(|\mathds{1}_{Q_{\mathfrak{R}_b}} \vw |))\|_{\mathcal{M}^{3, \sigma}_{t,x}}\leq C\|\mathcal{L}_1(|\mathds{1}_{Q_{\mathfrak{R}_b}} \vw |)\|_{\mathcal{M}^{120, 120}_{t,x}}=\|\mathcal{L}_1(|\mathds{1}_{Q_{\mathfrak{R}_b}} \vw |)\|_{\mathcal{M}^{\frac{q}{\lambda}, \frac{q}{\lambda}}_{t,x}},$$
where $q=\frac{24}{5}$ and $\lambda=\frac{1}{25}$. Thus, since $1<\frac{5}{q}$ and since $\lambda=1-\frac{q}{5}$ we can apply the Lemma \ref{Lemma_Hed} to obtain that 
$$\|\mathcal{L}_1(|\mathds{1}_{Q_{\mathfrak{R}_b}} \vw |)\|_{\mathcal{M}^{\frac{q}{\lambda}, \frac{q}{\lambda}}_{t,x}}\leq C\|\mathds{1}_{Q_{\mathfrak{R}_b}} \vw \|_{\mathcal{M}^{q,q}_{t,x}}=C\|\mathds{1}_{Q_{\mathfrak{R}_b}} \vw \|_{\mathcal{M}^{\frac{24}{5},\frac{24}{5}}_{t,x}}\leq C\|\mathds{1}_{Q_{\mathfrak{R}_b}} \vw \|_{\mathcal{M}^{6,6}_{t,x}}=C \|\mathds{1}_{Q_{\mathfrak{R}_b}} \vw \|_{L^{6}_{t,x}}<+\infty.$$
Now, for the term $\mathbb{V}_{b}$ above, using the same ideas as in  (\ref{Estimation_Ponctuelle_Riesz0})-(\ref{Estimation_Ponctuelle_Riesz1}) and applying again the Lemma \ref{lemma_locindi},  we obtain
\begin{eqnarray}
\| \mathbb{V}_{b}\|_{\mathcal{M}^{3, \sigma}_{t,x}}&\leq &C\|\mathds{1}_{Q_{R_2}}(\mathcal{L}_2(|\mathds{1}_{Q_{\mathfrak{R}_b}} \vw |))\|_{\mathcal{M}^{3, \sigma}_{t,x}}\notag\\
&\leq &C\|\mathds{1}_{Q_{R_2}}(\mathcal{L}_2(|\mathds{1}_{Q_{\mathfrak{R}_b}} \vw |))\|_{\mathcal{M}^{60, 60}_{t,x}}\label{FormuleMaxSigma}\\
&\leq &C\|\mathcal{L}_2(|\mathds{1}_{Q_{\mathfrak{R}_b}} \vw |)\|_{\mathcal{M}^{\frac{q}{\lambda}, \frac{q}{\lambda}}_{t,x}},\notag
\end{eqnarray}
where this time $q=\frac{12}{5}$ and $\lambda=\frac{1}{25}$. Since we have $2<\frac{5}{q}$ and $\lambda=1-\frac{2q}{5}$, we apply Lemma \ref{Lemma_Hed} and we have
$$C\|\mathcal{L}_2(|\mathds{1}_{Q_{\mathfrak{R}_b}} \vw |)\|_{\mathcal{M}^{\frac{q}{\lambda}, \frac{q}{\lambda}}_{t,x}}\leq C\|\mathds{1}_{Q_{\mathfrak{R}_b}} \vw \|_{\mathcal{M}^{q, q}_{t,x}}= C\|\mathds{1}_{Q_{\mathfrak{R}_b}} \vw \|_{\mathcal{M}^{\frac{12}{5}, \frac{12}{5}}_{t,x}}\leq \|\mathds{1}_{Q_{\mathfrak{R}_b}} \vw \|_{\mathcal{M}^{6,6}_{t,x}}=\|\mathds{1}_{Q_{\mathfrak{R}_b}} \vw \|_{L^{6}_{t,x}}<+\infty.$$
We can thus conclude that 
$$\|\mathds{1}_{Q_{R_2}}\vec{\mathbb{V}}_{11}\|_{\mathcal{M}^{3, \sigma}_{t,x}}<+\infty.$$
\end{itemize}
With all these estimates Proposition \ref{Propo_GainenMorreySigma} is now proven. \hfill$\blacksquare$
\begin{Remark}\label{Rem_LimitIteration_vu}
Note that the value of the index $\sigma$ of the Morrey space $\mathcal{M}_{t,x}^{3, \sigma}(\mathbb{R}\times \R)$ is potentially bounded by the information available over $\vw$ and the maximal possible value for this parameter is close to $\sigma=60$ (see the expression (\ref{FormuleMaxSigma}) above).
\end{Remark}
This result gives a small gain of integrability as we pass from an information on the Morrey space $\mathcal{M}_{t,x}^{3, \tau_0}$ to a control over the space $\mathcal{M}_{t,x}^{3, \sigma}$ with $\tau_0<\sigma$ with $\sigma$ close to $\tau_0$. This is of course not enough and we need to repeat the arguments above in order to obtain a better control. In this sense we have the following proposition:
\begin{Proposition}\label{MorreyPropo2}
Under the hypotheses of Theorem \ref{Teo_HolderRegularity} and within the framework of Proposition \ref{Propo_FirstMorreySpace}, there exists a radius $\bar R_2$ with $0<\bar R_2<R_2$ such that
\begin{equation}\label{Conclusion_SecondMorreySpace}
\mathds{1}_{Q_{\bar R_2}(t_0,x_0)}\vu\in \M_{t,x}^{3,60}(\mathbb{R}\times \R),
\end{equation}
\end{Proposition}
\noindent \textbf{Proof.}
By the Proposition \ref{Propo_GainenMorreySigma} above it follows that $\mathds{1}_{Q_{R_2}} \vu \in \mathcal{M}_{t,x}^{3,\sigma}(\mathbb{R}\times \R)$ with $\sigma$ very close to $\tau_0$ (say $\sigma=\tau_0+\epsilon$). Hence, with the information $\mathds{1}_{Q_{R_2}} \vu \in \mathcal{M}_{t,x}^{3, \tau_0+\epsilon}(\mathbb{R}\times \R)$ at hand, we can reapply the Proposition \ref{Propo_GainenMorreySigma} to obtain for some smaller radius $\bar {R}_2< R_2$ that $\mathds{1}_{Q_{\bar{R}_2}} \vu \in \mathcal{M}_{t,x}^{3, \sigma_1}(\mathbb{R}\times \R)$ where $\sigma_1=\sigma+\epsilon=\tau_0+2\epsilon$. Iterating these arguments as long as necessary, we obtain the information $\mathds{1}_{Q_{R_2}} \vu \in \mathcal{M}_{t,x}^{3, 60}(\mathbb{R}\times \R)$ where the value $\sigma=60$ is fixed by the information available for the quantity $\vw$ which is the only term that is fixed: see the computation leading to the estimate (\ref{FormuleMaxSigma}) and Remark \ref{Rem_LimitIteration_vu}. Let us note that a slight abuse of language has been used for the radius $\bar R_2$: at each iteration this radius is smaller and smaller, but in order to maintain the notations we still denote the final radius by $\bar R_2$. \hfill $\blacksquare$
\section{A first gain of information for the variable $\vw$}\label{Sec_Propiedades_Vw}

Note that the Proposition \ref{MorreyPropo2} and the Corollary \ref{corolarioMorrey} give interesting control (on a small neighborhood of a point $(t_0, x_0)$) for the variable $\vu$. Remark also that Theorem \ref{Teo_InterdepenciaMicropolar} gives some information for the variable $\vw$:
\begin{equation}\label{Information_PremiereEtape}
\begin{split}
\mathds{1}_{Q_{\bar R_2}(t_0,x_0)} \vu\in \mathcal{M}_{t,x}^{3,60}(\mathbb{R}\times \R),&\quad \mathds{1}_{Q_{\mathfrak{r}_1}(t_0,x_0)} \vu\in L_{t,x}^{6}(\mathbb{R}\times \R), \quad
\mathds{1}_{Q_{R_1}(t_0,x_0)}\vn\otimes \vu\in \mathcal{M}^{2,\tau_1}_{t,x}(\mathbb{R}\times \R),\\[3mm]
&\qquad\mathds{1}_{Q_{\mathfrak{r}_2}(t_0,x_0)} \vw\in L_{t,x}^{6}(\mathbb{R}\times \R),
\end{split}
\end{equation}
where 
\begin{equation}\label{Info_Radios}
0<\bar R_2< R_2<\mathfrak{r}_2<\mathfrak{r}_1<R_1<R<1,
\end{equation}
with $\tau_0=6$ and $\tau_1=\frac{30}{11}$ (which is given by the condition $\frac{1}{\tau_1}=\frac{1}{\tau_0}+\frac{1}{5}$, see the Corollary \ref{corolarioMorrey}). Note that we have $\frac{120}{45}<\tau_1=\frac{30}{11}$.\\

We will exploit all this information in order to derive some Morrey control for the variable $div(\vw)$, indeed, we have:
\begin{Proposition}\label{Propo_divw}
Under the general hypotheses of Theorem \ref{Teo_HolderRegularity}, if we have the controls (\ref{Information_PremiereEtape}) over $\vu$ and $\vw$ then we have, for some radius $0<R_3<\bar R_2$, we have
$$\mathds{1}_{Q_{R_3}(t_0,x_0)}div(\vw) \in\M_{t,x}^{\frac{6}{5},\frac{60}{11}}(\mathbb{R}\times\R).$$ 
\end{Proposition}
\textbf{Proof.} We first apply the divergence operator to the equation satisfied by $\vw$ (see the system (\ref{MicropolarEquationsIntro})) and since we have the identities $div(\vn div(\vw))=\Delta \vw$ and $div(\rot \vu)\equiv 0$, we obtain
$$\partial_t div(\vw)=2\Delta div(\vw)-div(\vw)-div((\vu\cdot\vn)\vw).$$
Consider now $\bar \phi: \mathbb{R}\times\mathbb{R}^3\longrightarrow \mathbb{R}$ a non-negative function such that $\bar \phi\in \mathcal{C}_0^{\infty}(\mathbb{R}\times \mathbb{R}^3)$ and such that
\begin{equation}\label{Def_FuncAuxiliarPourDivW}
\bar \phi \equiv 1\;\; \text{over}\; \; Q_{\rho_b}(t_0,x_0),\; \; supp(\bar \phi)\subset Q_{\rho_a}(t_0,x_0),
\end{equation}
where we have
\begin{equation}\label{Info_Radios1}
0<R_3<\rho_b<\rho_a<\bar R_2, 
\end{equation} 
where the radius $\bar R_2$ is fixed in (\ref{Info_Radios}).  With the help of this auxiliar function we define the variable $\mathcal{W}$ by 
$$\mathcal{W}=\bar \phi div(\vw),$$ 
note that, due to the support properties of the function $\bar \phi$ we have $\mathds{1}_{Q_{R_3}}\mathcal{W}=\mathds{1}_{Q_{R_3}}div(\vw)$. 
If we study the evolution of $\mathcal{W}$ we obtain:
\begin{eqnarray*}
\partial_t\mathcal{W}&=&(\partial_t \bar \phi )div(\vw)+ \bar \phi \bigg(2\Delta div(\vw)-div(\vw)-div((\vu\cdot\vn)\vw)\bigg)\\
&=&2\Delta \mathcal{W} +(\partial_t \bar \phi +2\Delta \bar \phi-\bar \phi)div(\vw)-4\sum_{i=1}^3  \partial_i\big((\partial_i \bar \phi) div(\vw)\big)- \bar \phi div((\vu\cdot\vn)\vw),
\end{eqnarray*}
where we used the identity $\displaystyle{\bar \phi \Delta div(\vw)=\Delta(\bar \phi div(\vw))+\Delta \bar \phi div (\vw)-2\sum_{i=1}^3  \partial_i\big((\partial_i \bar \phi) div(\vw)\big)}$. Recall now that we also have the identity (since $div(\vu)=0$):
$$\bar \phi div((\vu\cdot\vn)\vw)=\bar \phi div\big(div(\vw\otimes\vu)\big)=div\big(div(\bar\phi \vw\otimes \vu)\big)-div(\vw\otimes \vu\cdot \vn \bar\phi)-\vn\bar\phi\cdot div(\vw\otimes \vu),$$
and we obtain
\begin{eqnarray*}
\partial_t\mathcal{W}&=&2\Delta \mathcal{W} +(\partial_t \bar \phi +2\Delta \bar \phi-\bar \phi)div(\vw)-4\sum_{i=1}^3  \partial_i\big((\partial_i \bar \phi) div(\vw)\big)- div\big(div(\bar\phi \vw\otimes \vu)\big)\\
&&+div(\vw\otimes \vu\cdot \vn \bar\phi)+\vn\bar\phi\cdot div(\vw\otimes \vu).
\end{eqnarray*}
Thus, since we have $\mathcal{W}(0, \cdot)=0$ (by the properties of the localizing function $\bar\phi$ given in (\ref{Def_FuncAuxiliarPourDivW})), applying the Duhamel formula we can write:
\begin{eqnarray}
\mathcal{W}(t,x)=\underbrace{\int_0^t e^{2(t-s)\Delta}(\partial_t \bar \phi+2\Delta \bar\phi-\bar\phi )div(\vw)ds}_{\mathcal{W}_1}-4\sum_{i=1}^3\underbrace{\int_0^t e^{2(t-s)\Delta}\partial_i \big((\partial_i \bar\phi)div(\vw)\big)ds}_{\mathcal{W}_2}\label{Ws}\\
-\underbrace{\int_0^t e^{2(t-s)\Delta} div\big(div(\bar\phi \vw\otimes \vu)\big)ds}_{\mathcal{W}_3}+\underbrace{\int_0^t e^{2(t-s)\Delta} div(\vw\otimes \vu\cdot \vn \bar\phi)ds}_{\mathcal{W}_4}+\underbrace{\int_0^t e^{2(t-s)\Delta} \vn\bar\phi\cdot div(\vw\otimes \vu)ds}_{\mathcal{W}_5},\notag
\end{eqnarray}
and we will estimate each one of the terms above.

\begin{itemize}
\item For the first term $\mathcal{W}_1$ we write, 
\begin{eqnarray}
|\mathds{1}_{Q_{R_3}}\mathcal{W}_1|&=&\left|\mathds{1}_{Q_{R_3}}\int_0^t e^{2(t-s)\Delta}div\big((\partial_t \bar \phi+2\Delta \bar\phi-\bar\phi )\vw\big)ds\right|\label{Estimation_W1}\\
&&+\left|\mathds{1}_{Q_{R_3}}\int_0^t e^{2(t-s)\Delta}\big(\vn\big(\partial_t \bar \phi+2\Delta \bar\phi-\bar\phi \big)\big)\cdot\vw ds\right|,\notag
\end{eqnarray}
since the convolution kernel of the semi-group $e^{2(t-s)\Delta}$ is the usual 3D heat kernel $\mathfrak{g}_{2t}$, thus by the decay properties of the heat kernel, by the properties of the test function $\bar\phi$ (see (\ref{Def_FuncAuxiliarPourDivW})) and by the definition of the parabolic Riesz potentials $\mathcal{L}_1$ and $\mathcal{L}_2$ given in (\ref{Def_ParabolicRieszPotential}), we can write the estimate
 \begin{eqnarray}
|\mathds{1}_{Q_{R_3}}\mathcal{W}_1|&\leq &C\mathds{1}_{Q_{R_3}}\int_{\mathbb{R}}\int_{\mathbb{R}^3} \frac{|\mathds{1}_{Q_{\rho_a}} \vw(s,y)|}{(|t-s|^{\frac{1}{2}}+|x-y|)^4}dyds+C\mathds{1}_{Q_{R_3}}\int_{\mathbb{R}}\int_{\mathbb{R}^3} \frac{|\mathds{1}_{Q_{\rho_a}} \vw(s,y)|}{(|t-s|^{\frac{1}{2}}+|x-y|)^3}dyds\notag\\
&&\leq C\mathds{1}_{Q_{R_3}}(\mathcal{L}_1(|\mathds{1}_{Q_{\rho_a}} \vw |))(t,x)+C\mathds{1}_{Q_{R_3}}(\mathcal{L}_2(|\mathds{1}_{Q_{\rho_a}} \vw |))(t,x),\label{Estimation_W11}
\end{eqnarray}
and we have
\begin{eqnarray*}
\|\mathds{1}_{Q_{R_3}}\mathcal{W}_1\|_{\M_{t,x}^{\frac{6}{5},\frac{60}{11}}}&\leq& C\left\| \mathds{1}_{Q_{R_3}}(\mathcal{L}_1(|\mathds{1}_{Q_{\rho_a}} \vw |))\right\|_{\M_{t,x}^{\frac{6}{5},\frac{60}{11}}}+C\left\|\mathds{1}_{Q_{R_3}}(\mathcal{L}_2(|\mathds{1}_{Q_{\rho_a}} \vw |))\right\|_{\M_{t,x}^{\frac{6}{5},\frac{60}{11}}}.
\end{eqnarray*}
For the first term above, since $\frac{60}{11}\leq \frac{15}{2}$, we set $p=\frac65$, $q=\frac{9}{2}$ and $\lambda=\frac{1}{10}$ and by Lemma \ref{lemma_locindi} we obtain
$$\left\| \mathds{1}_{Q_{R_3}}(\mathcal{L}_1(|\mathds{1}_{Q_{\rho_a}} \vw |))\right\|_{\M_{t,x}^{\frac{6}{5},\frac{60}{11}}}\leq C\|\mathds{1}_{Q_{R_3}}\mathcal{L}_1(|\mathds{1}_{Q_{\rho_a}}\vw|)\|_{\M_{t,x} ^{\frac{6}{5},\frac{15}{2}}}\leq C\|\mathcal{L}_1(|\mathds{1}_{Q_{\rho_a}}\vw|)\|_{\M_{t,x}^{\frac{6}{\lambda 5},\frac{9}{2\lambda}}},$$
since $\frac{6}{5}<\frac{6}{\lambda 5}$ and $\frac{15}{2}<\frac{9}{2\lambda}$. Thus, applying Lemma \ref{Lemma_Hed} (and Lemma \ref{lemma_locindi}), we have
$$\|\mathcal{L}_1(|\mathds{1}_{Q_{\rho_a}}\vw|)\|_{\M_{t,x}^{\frac{6}{\lambda 5},\frac{9}{2\lambda}}}\le C\|\mathds{1}_{Q_{\rho_a}}\vw \|_{\M_{t,x}^{\frac{6}{5},\frac{9}{2}}}\le C \|\mathds{1}_{Q_{\rho_a}}\vw \|_{L_{t,x}^{6}}<+\infty,$$
since we have the control $\mathds{1}_{Q_{\mathfrak{r}_2}}\vw \in L_{t,x}^{6}(\mathbb{R}\times\R)$ given in (\ref{Information_PremiereEtape}) and we have by (\ref{Info_Radios}) and (\ref{Info_Radios1}) that $\rho_a<\bar R_2<\mathfrak{r}_2$.\\

For the second term that we need to study, we fix $p=\frac65$, $q=\frac{12}{5}$ and $\lambda=\frac{1}{25}$, by applying Lemma \ref{lemma_locindi} and  Lemma \ref{Lemma_Hed} we obtain
\begin{eqnarray*}
\left\|\mathds{1}_{Q_{R_3}}(\mathcal{L}_2(|\mathds{1}_{Q_{\rho_a}} \vw |))\right\|_{\M_{t,x}^{\frac{6}{5},\frac{60}{11}}}&\leq &C\|\mathds{1}_{Q_{R_3}}\mathcal{L}_2(|\mathds{1}_{Q_{\rho_a}}\vw|)\|_{\M_{t,x} ^{\frac65,\frac{15}{2}}}\leq \|\mathcal{L}_2(|\mathds{1}_{Q_{\rho_a}}\vw|)\|_{\M_{t,x}^{\frac{6}{\lambda 5},\frac{12}{\lambda 5}}}\\
&\leq &C\|\mathds{1}_{Q_{\rho_a}}\vw \|_{\M_{t,x}^{\frac{6}{5},\frac{12}{5}}}\le C \|\mathds{1}_{Q_{\rho_a}}\vw \|_{L_{t,x} ^{6}}\le C \|\mathds{1}_{Q_{\mathfrak{r}_2}}\vw \|_{L_{t,x} ^{6}}<+\infty,
\end{eqnarray*}
where we used the information (\ref{Information_PremiereEtape}) and the relationships (\ref{Info_Radios})-(\ref{Info_Radios1}). With these two estimates at hand we conclude that $\|\mathds{1}_{Q_{R_3}}\mathcal{W}_1\|_{\M_{t,x}^{\frac{6}{5}, \frac{15}{2}}}<+\infty$.

\item For the term $\mathcal{W}_2$ of (\ref{Ws}) we need to study, for all $1\le i\le 3$, the quantities
$$\mathbb{I}_i=\left|\mathds{1}_{Q_{R_3}}\int_0^t e^{2(t-s)\Delta}\partial_i \big((\partial_i \bar\phi)div(\vw)\big)ds\right|,$$
and we write 
\begin{equation}\label{Estimation_NoyauBorne}
\mathbb{I}_i\leq \left|\mathds{1}_{Q_{R_3}}\int_0^t e^{2(t-s)\Delta}\partial_i \big(div((\partial_i \bar\phi)\vw)\big)ds\right|+\left|\mathds{1}_{Q_{R_3}}\int_0^t e^{2(t-s)\Delta}\partial_i \big([\vn(\partial_i \bar\phi)]\cdot \vw\big)ds\right|.
\end{equation}
We study the first term above and by the support properties of the function $\bar\phi$ given in (\ref{Def_FuncAuxiliarPourDivW}), we have for $1\leq i,j\leq 3$:
\begin{equation}\label{Estimation_NoyauBorne1}
\left|\mathds{1}_{Q_{R_3}}\int_0^t e^{2(t-s)\Delta}\partial_i \partial_j\big((\partial_i \bar\phi)\vw)\big)ds\right|\leq C\int_{\mathbb{R}}\int_{\R}\frac{\mathds{1}_{Q_{R_3}}(t,x)\mathds{1}_{\mathcal{C}(\rho_b, \rho_a)}(s,y)|\vw(s,y)|}{(|t-s|^{\frac12}+|x-y|)^5} dyds,
\end{equation}
where the set $\mathcal{C}(\rho_b, \rho_a)$ is the corona defined by $Q_{\rho_a}\setminus Q_{\rho_b}$. Noting that $(t,x)\in Q_{R_3}$ and that $(s,y)\in \mathcal{C}(\rho_b, \rho_a)$, since we have $R_3<\rho_b$ by (\ref{Info_Radios1}), the convolution kernel $\frac{\mathds{1}_{Q_{R_3}}(t,x)\mathds{1}_{\mathcal{C}(\rho_b, \rho_a)}(s,y)}{(|t-s|^{\frac12}+|x-y|)^5}$ is bounded and we can write
\begin{equation}\label{Estimation_NoyauBorne2}
\left\|\mathds{1}_{Q_{R_3}}\int_0^t e^{2(t-s)\Delta}\partial_i \partial_j\big((\partial_i \bar\phi)\vw)\big)ds\right\|_{L^\infty_{t,x}}\leq C\|\mathds{1}_{\mathcal{C}(\rho_b, \rho_a)}\vw\|_{L^1_{t,x}}\leq C\|\mathds{1}_{Q_{\mathfrak{r}_2}}\vw\|_{L^6_{t,x}}<+\infty,
\end{equation}
(recall (\ref{Info_Radios})-(\ref{Info_Radios1})), from which we deduce that 
$$\left\|\mathds{1}_{Q_{R_3}}\int_0^t e^{2(t-s)\Delta}\partial_i \big(div((\partial_i \bar\phi)\vw)\big)ds\right\|_{\M_{t,x}^{\frac{6}{5}, \frac{60}{11}}}\leq C\|\mathds{1}_{Q_{\mathfrak{r}_2}}\vw\|_{L^6_{t,x}}<+\infty.$$
The second term of (\ref{Estimation_NoyauBorne}) has the same structure as the first term in (\ref{Estimation_W1}), and thus by the same arguments we can write
$$\left\|\mathds{1}_{Q_{R_3}}\int_0^t e^{2(t-s)\Delta}\partial_i \big([\vn(\partial_i \bar\phi)]\cdot \vw\big)ds\right\|_{\M_{t,x}^{\frac{6}{5}, \frac{60}{11}}}\leq C\|\mathds{1}_{Q_{\mathfrak{r}_2}}\vw\|_{L^6_{t,x}}<+\infty.$$
\item We study now the term $\mathcal{W}_3$ defined in (\ref{Ws}) and we write
$$\|\mathds{1}_{Q_{R_3}}\mathcal{W}_3\|_{\M_{t,x}^{\frac{6}{5}, \frac{60}{11}}}= \left\|\mathds{1}_{Q_{R_3}}\int_0^t e^{2(t-s)\Delta} div\big(div(\bar\phi \vw\otimes \vu)\big)ds\right\|_{\M_{t,x}^{\frac{6}{5}, \frac{60}{11}}},$$
and by the maximal regularity of the heat kernel in Morrey spaces (the Theorem 7.3 of \cite{PGLR0} can be generalized to parabolic Morrey spaces), we have
$$\|\mathds{1}_{Q_{R_3}}\mathcal{W}_3\|_{\M_{t,x}^{\frac{6}{5}, \frac{60}{11}}}\leq C \left\|\bar\phi \vw\otimes \vu\right\|_{\M_{t,x}^{\frac{6}{5}, \frac{60}{11}}},$$
now, using the H\"older inequalities for Morrey spaces stated in Lemma \ref{lemma_Product} (with $\frac{5}{6}=\frac{1}{2}+\frac{1}{3}$ and $\frac{11}{60}=\frac{1}{6}+\frac{1}{60}$) and the properties of the localizing function $\bar\phi$, we obtain
\begin{eqnarray*}
\|\mathds{1}_{Q_{R_3}}\mathcal{W}_3\|_{\M_{t,x}^{\frac{6}{5}, \frac{60}{11}}}&\leq &C \|\mathds{1}_{Q_{\rho_a}} \vw\|_{\M_{t,x}^{2,6}}\|\mathds{1}_{Q_{\rho_a}} \vu\|_{\M_{t,x}^{6,60}}\\
&\leq & C \|\mathds{1}_{Q_{\mathfrak{r}_2}} \vw\|_{\M_{t,x}^{6,6}}\|\mathds{1}_{Q_{\bar R_2}} \vu\|_{\M_{t,x}^{6,60}}\leq  C \|\mathds{1}_{Q_{\mathfrak{r}_2}} \vw\|_{L_{t,x}^{6}}\|\mathds{1}_{Q_{\bar R_2}} \vu\|_{\M_{t,x}^{6,60}}<+\infty,
\end{eqnarray*}
where in the last estimate above we used Lemma \ref{lemma_locindi}, the information available in (\ref{Information_PremiereEtape}) and the relationships (\ref{Info_Radios})-(\ref{Info_Radios1}).

\item For the term $\mathcal{W}_4$ given in (\ref{Ws}) we have, following the same arguments given in (\ref{Estimation_NoyauBorne1}):
$$|\mathds{1}_{Q_{R_3}}\mathcal{W}_4|= \left|\mathds{1}_{Q_{R_3}}\int_0^t e^{2(t-s)\Delta} div(\vw\otimes \vu\cdot \vn \bar\phi)  ds\right|\leq C\int_{\mathbb{R}}\int_{\R}\frac{\mathds{1}_{Q_{R_3}}(t,x)\mathds{1}_{\mathcal{C}(\rho_b, \rho_a)}(s,y)|\vw\otimes \vu(s,y)|}{(|t-s|^{\frac12}+|x-y|)^4} dyds,$$
and thus, by the ideas given in (\ref{Estimation_NoyauBorne1})-(\ref{Estimation_NoyauBorne2}) we can write
\begin{eqnarray*}
\left\|\int_{\mathbb{R}}\int_{\R}\frac{\mathds{1}_{Q_{R_3}}(t,x)\mathds{1}_{\mathcal{C}(\rho_b, \rho_a)}|\vw\otimes \vu(s,y)|}{(|t-s|^{\frac12}+|x-y|)^4} dyds\right\|_{L^\infty_{t,x}}&\leq &C\|\mathds{1}_{\mathcal{C}(\rho_b, \rho_a)}\vw\otimes \vu\|_{L^1_{t,x}}\leq C\|\mathds{1}_{Q_{\mathfrak{r}_2}}\vw\otimes \vu\|_{L^3_{t,x}}\\
&\leq &C\|\mathds{1}_{Q_{\mathfrak{r}_2}}\vw\|_{L^6_{t,x}}\|\mathds{1}_{Q_{\mathfrak{r}_1}}\vu\|_{L^6_{t,x}}<+\infty,
\end{eqnarray*}
where we applied the H\"older inequalities, Lemma \ref{lemma_locindi} (in the Lebesgue space setting) and the relationships (\ref{Info_Radios})-(\ref{Info_Radios1}). With these estimates at hand, we easily deduce that 
$$\|\mathds{1}_{Q_{R_3}}\mathcal{W}_4\|_{\M_{t,x}^{6,60}}<+\infty.$$
\item For the last term $\mathcal{W}_5$ of (\ref{Ws}) we have
$$|\mathds{1}_{Q_{R_3}}\mathcal{W}_5|= \left|\mathds{1}_{Q_{R_3}}\int_0^t e^{2(t-s)\Delta}   \vn\bar\phi\cdot div(\vw\otimes \vu) ds\right|.$$
Thus, for $1\leq i,j,k,l\leq 3$ we need to study the quantities
\begin{eqnarray*}
\mathbb{J}_{i,j,k,l}&=&\left|\mathds{1}_{Q_{R_3}}\int_0^t e^{2(t-s)\Delta}   (\partial_i\bar\phi)\partial_j(w_k u_l) ds\right|\notag\\
&\leq &\left|\mathds{1}_{Q_{R_3}}\int_0^t e^{2(t-s)\Delta}   \partial_j\bigg((\partial_i\bar\phi)(w_k u_l)\bigg) ds\right|+\left|\mathds{1}_{Q_{R_3}}\int_0^t e^{2(t-s)\Delta}   (\partial_j\partial_i\bar\phi)(w_k u_l) ds\right|,
\end{eqnarray*}
where we used the identity $ (\partial_i\bar\phi)\partial_j(w_k u_l)= \partial_j\bigg((\partial_i\bar\phi)(w_k u_l)\bigg) -(\partial_j\partial_i\bar\phi)(w_k u_l)$. Now, due to the properties of the heat kernel and the support properties of the function $\bar\phi$, we obtain the inequality
\begin{eqnarray*}
\mathbb{J}_{i,j,k,l}&\leq &C\int_{\mathbb{R}}\int_{\R}\frac{\mathds{1}_{Q_{R_3}}(t,x)\mathds{1}_{\mathcal{C}(\rho_b, \rho_a)}|w_k u_l(s,y)|}{(|t-s|^{\frac12}+|x-y|)^4} dyds\\
&&+C\int_{\mathbb{R}}\int_{\R}\frac{\mathds{1}_{Q_{R_3}}(t,x)\mathds{1}_{\mathcal{C}(\rho_b, \rho_a)}|w_k u_l(s,y)|}{(|t-s|^{\frac12}+|x-y|)^3} dyds.
\end{eqnarray*}
Now, by the same arguments given in  (\ref{Estimation_NoyauBorne1})-(\ref{Estimation_NoyauBorne2}) we obtain
\begin{eqnarray*}
\left\|\int_{\mathbb{R}}\int_{\R}\frac{\mathds{1}_{Q_{R_3}}(t,x)\mathds{1}_{\mathcal{C}(\rho_b, \rho_a)}|w_k u_l(s,y)|}{(|t-s|^{\frac12}+|x-y|)^4} dyds\right\|_{L^\infty_{t,x}}
&+&\left\|\int_{\mathbb{R}}\int_{\R}\frac{\mathds{1}_{Q_{R_3}}(t,x)\mathds{1}_{\mathcal{C}(\rho_b, \rho_a)}|w_k u_l(s,y)|}{(|t-s|^{\frac12}+|x-y|)^3} dyds\right\|_{L^\infty_{t,x}}\\
&\leq &C\|\mathds{1}_{Q_{\mathfrak{r}_2}}w_k u_l\|_{L^1_{t,x}}+C\|\mathds{1}_{Q_{\mathfrak{r}_2}}w_k u_l\|_{L^1_{t,x}}\\
&\leq&C\|\mathds{1}_{Q_{\mathfrak{r}_2}}\vw\|_{L^6_{t,x}}\|\mathds{1}_{Q_{\mathfrak{r}_1}}\vu\|_{L^6_{t,x}}<+\infty,
\end{eqnarray*}
and with these estimates for $1\leq i,j,k,l\leq 3$, we easily deduce that 
$$\|\mathds{1}_{Q_{R_3}}\mathcal{W}_5\|_{\M_{t,x}^{6,60}}<+\infty.$$
\end{itemize}
With all these controls, Proposition \ref{Propo_divw} is proven. \hfill$\blacksquare$

\section{The end of the proof of Theorem \ref{Teo_HolderRegularity}}\label{Sec_Final}
The key result for obtaining a gain of regularity is the following lemma coming from the theory of parabolic equations (see \cite{Lady68,PGLR1}).
\begin{Lemma}\label{Lemma_parabolicHolder}
Let $\sigma$  be a smooth homogeneous function over $\R\setminus \{0\}$,  of exponent 1 with 
$\sigma(D)$ the Fourier multiplier associated.
Consider a vector field  $\vec \Phi \in \mathcal{M}_{t,x}^{\p_0,\q_0}(\mathbb{R}\times\R)$ and a scalar function $h\in \mathcal{M}_{t,x}^{\p_0,\q_1}(\mathbb{R}\times\R)$ such that $1\le \p_0\le \q_0$, with  $\frac{1}{\q_0}=\frac{2-\alpha}{5}$,$\frac{1}{\q_1}=\frac{1-\alpha}{5}$, for $0<\alpha<1$. Then, the function $\vec {v} $ equal to $0$ for $t\le 0$  and
\begin{equation*}
\vec v(t,x)=\int_0^t e^{(t-s)\Delta }(\vec \Phi (s,\cdot)+\sigma(D) h(s,\cdot))ds,
\end{equation*}
for $t>0$,  is H\"older continuous  of exponent $\alpha$ with respect to the parabolic distance.  
\end{Lemma} 
In order to apply this lemma to the proof of Theorem \ref{Teo_HolderRegularity}, we will first localize the \emph{full} micropolar equations (\ref{MicropolarEquationsIntro}) and then we will show that each term of the corresponding Duhamel formula belongs either to the space $\mathcal{M}_{t,x}^{\p_0,\q_0}(\mathbb{R}\times\R)$ or to the space $ \mathcal{M}_{t,x}^{\p_0,\q_1}(\mathbb{R}\times\R)$.\\

We start by localizing the problem and for this we consider  $\phi: \mathbb{R}\times \R\longrightarrow \mathbb{R}$ a test function such that $supp(\phi)\subset ]-\frac{1}{4},\frac{1}{4}[\times B(0,\frac{1}{2})$ and $\phi \equiv 1\;\; \text{over}\; \; ]-\frac{1}{16},\frac{1}{16}[\times B(0,\frac{1}{4})$. We consider next a radius $\mathbf{R}>0$ such that 
\begin{equation}\label{Info_Radios2}
4\mathbf{R}<R_3<\bar R_2<\mathfrak{r}_2<\mathfrak{r}_1<R_1<R<1,
\end{equation}
where $R_3$ is the radius of Proposition \ref{Propo_divw}, $\bar R_2$ is the radius of Proposition \ref{MorreyPropo2},  $\mathfrak{r}_1, \mathfrak{r}_2$ are the radii from Theorem \ref{Teo_InterdepenciaMicropolar} and $R_1$ is the radius obtained in Proposition \ref{Propo_FirstMorreySpace}. We then write
\begin{equation}\label{Def_Eta}
\eta(t,x)=\phi\left(\frac{t-t_0}{\mathbf{R}^2},\frac{x-x_0}{\mathbf{R}}\right),
\end{equation}
and we consider the variable $\vUU$ defined by the formula
\begin{equation}\label{defU}
\vUU=\eta(\vu+\vw),
\end{equation}
then, by the properties of the auxiliar function $\eta$, we have the identity $\vUU=\vu+\vw$ over a small neighborhood of the point $(t_0,x_0)$, the support of the variable $\vUU$ is contained in the parabolic ball $Q_{\mathbf{R}}$ and moreover we also have $\vUU(0,x)=0$. Thus, if we study the evolution of this variable, following the system (\ref{MicropolarEquationsIntro}), we have
\begin{eqnarray*}
\partial_t\vUU&=&(\partial_t\eta)(\vu+\vw)+\eta \Delta (\vu+\vw)-\eta ((\vu\cdot \vn)\vu)-\eta\vn p+\frac{1}{2}\eta\rot\vw +\eta\vn div(\vw)\\
&&-\eta\vw-\eta((\vu\cdot \vn)\vw)+\frac{1}{2}\eta\rot\vu.
\end{eqnarray*}
We use now the identity $\displaystyle{\eta \Delta (\vu+\vw)=\Delta\vUU- \Delta \eta(\vu+\vw)-2\sum_{i=1}^3 (\partial_i \eta)(\partial_i (\vu+\vw))}$ to obtain the equation
\begin{eqnarray*}
\partial_t\vUU&=&\Delta\vUU+(\partial_t\eta- \Delta \eta)(\vu+\vw)-2\sum_{i=1}^3 (\partial_i \eta)(\partial_i (\vu+\vw))-\eta((\vu\cdot \vn)\vu)-\eta\vn p+\frac{1}{2}\eta\rot\vw \\
&&+\eta\vn div(\vw)-\eta\vw-\eta((\vu\cdot \vn)\vw)+\frac{1}{2}\eta\rot\vu.
\end{eqnarray*}
In the expression above, we need to rewrite six particular terms, indeed, since we have the identities 
\begin{eqnarray*}
\eta \vn div(\vw)&=&\vn(\eta div(\vw))-(\vn \eta) div(\vw),\\
\sum_{i=1}^3 (\partial_i \eta)(\partial_i (\vu+\vw))&=&\sum_{i=1}^3 \partial_i( (\partial_i \eta)(\vu+\vw))-\sum_{i=1}^3 (\partial_i\partial_i \eta)(\vu+\vw),\\
\eta((\vu\cdot\vn)\vu)&=&div(\eta \vu\otimes \vu)-\vu\otimes \vu\cdot \vn \eta,\\
\eta\rot\vw &=& \rot(\eta\vw)-(\vn\eta)\wedge \vw, \quad \mbox{and } \eta\rot\vu = \rot(\eta\vu)-(\vn\eta)\wedge \vu,\\
\mbox{and}\quad \eta((\vu\cdot\vn)\vw)&=&\eta div(\vw\otimes \vu)=div(\eta \vw\otimes \vu)-\vw\otimes \vu\cdot \vn \eta,
\end{eqnarray*}
 we obtain 
\begin{eqnarray*}
\partial_t\vUU&=&\Delta\vUU+(\partial_t\eta+\Delta \eta)(\vu+\vw)-2\sum_{i=1}^3 \partial_i( (\partial_i \eta)(\vu+\vw))-div(\eta \vu\otimes \vu)+\vu\otimes \vu\cdot \vn \eta\\
&&-\eta\vn p+\frac{1}{2}\rot(\eta\vw)-\frac{1}{2}(\vn\eta)\wedge \vw+\vn(\eta div(\vw))-(\vn \eta) div(\vw) \\
&&-\eta\vw-div(\eta \vw\otimes \vu)+\vw\otimes \vu\cdot \vn \eta+\frac{1}{2} \rot(\eta\vu)-\frac{1}{2} (\vn\eta)\wedge \vu.
\end{eqnarray*}
We rewrite this equation in the following form:
\begin{equation}\label{EquationPourUtilisationMorrey}
\begin{cases}
\partial_t \vec {\mathcal{U}}=\Delta \vec {\mathcal{U}}+\vec{\mathcal{A}}+\displaystyle{\sum_{i=1}^3}\partial_i\vec{\mathcal{B}_i}+\vn \mathcal{C}+\rot \vec{\mathcal{D}}+div \mathbb{E},\\[3mm]
\vec {\mathcal{U}}(0,x)=0,
\end{cases}
\end{equation}
where the vector $\vec{\mathcal{A}}$ is given by
\begin{eqnarray}
\vec{\mathcal{A}}&=&(\partial_t\eta+\Delta \eta)(\vu+\vw)+\vu\otimes \vu\cdot \vn \eta -\eta \vn p-\frac{1}{2}(\vn\eta)\wedge \vw \label{A_Morrey}\\
&&-(\vn \eta) div(\vw)-\eta\vw+\vw\otimes \vu\cdot \vn \eta,\notag
\end{eqnarray}
the vector $\vec{\mathcal{B}}_i$ (for $1\leq i\leq 3$) is given by
\begin{equation} \label{B_Morrey}
\vec{\mathcal{B}}_i=2(\partial_i \eta)(\vu+\vw),
\end{equation} 
the scalar function $\mathcal{C}$ is given by
\begin{equation} \label{C_Morrey}
\mathcal{C}=\eta div(\vw),
\end{equation}
the vector $\vec{\mathcal{D}}$ is given by
\begin{equation}\label{D_Morrey}
\vec{\mathcal{D}}=\frac{1}{2}\eta(\vw+\vu),
\end{equation}
and finally, the tensor $\mathbb{E}$ is defined by the expression
\begin{equation}\label{E_Morrey}
\mathbb{E}=-\eta (\vu\otimes \vu +\vw\otimes \vu).
\end{equation}
Thus, by the Duhamel formula, the solution of the equation (\ref{EquationPourUtilisationMorrey}) can be written in the following manner:
\begin{equation}\label{Duhamel_todoMicroPolar}
\vUU(t,x)=\int_0^t e^{(t-s)\Delta }\bigg(\vec{\mathcal{A}}+\displaystyle{\sum_{i=1}^3}\partial_i\vec{\mathcal{B}_i}+\vn \mathcal{C}+\rot \vec{\mathcal{D}}+div \mathbb{E}\bigg)ds,
\end{equation}
thus, in order to apply Lemma \ref{Lemma_parabolicHolder} to this system and obtain a parabolic gain of regularity, we only need to prove that the quantities $\vec{\mathcal{A}}$, $\vec{\mathcal{B}_i}$, $\mathcal{C}$, $\vec{\mathcal{D}}$ and $\mathbb{E}$, defined in (\ref{A_Morrey})-(\ref{E_Morrey}) respectively, satisfy:
\begin{equation}\label{Objetivo_MorreyLadyz}
\vec{\mathcal{A}} \in \M_{t,x}^{\p_0,\q_0}(\mathbb{R}\times\R)\quad \mbox{and}\quad \vec{\mathcal{B}}_i, \mathcal{C}, \vec{\mathcal{D}}, \mathbb{E} \in \M_{t,x}^{\p_0,\q_1}(\mathbb{R}\times\R),
\end{equation}
where $1\le \p_0\leq \frac{6}{5}\le \q_0$, with  $\frac{1}{\q_0}=\frac{2-\alpha}{5}$,$\frac{1}{\q_1}=\frac{1-\alpha}{5}$, for some $0<\alpha<\frac{1}{24}$. \\

Let us start with the quantity $\vec{\mathcal{A}}$ and we have
\begin{Lemma}\label{Lemma_EstimationA}
For the term $\vec{\mathcal{A}}$ defined in (\ref{A_Morrey}) we have 
$$\|\vec{\mathcal{A}}\|_{\M_{t,x}^{\p_0,\q_0}}<+\infty.$$
\end{Lemma}
{\bf Proof.} By definition we have
\begin{eqnarray}
\|\vec{\mathcal{A}}\|_{\M_{t,x}^{\p_0,\q_0}}\leq \underbrace{\|(\partial_t\eta+\Delta \eta)(\vu+\vw)\|_{\M_{t,x}^{\p_0,\q_0}}}_{(1)}+\underbrace{\|\vu\otimes \vu\cdot \vn \eta\|_{\M_{t,x}^{\p_0,\q_0}}}_{(2)}+\underbrace{\|\eta \vn p\|_{\M_{t,x}^{\p_0,\q_0}}}_{(3)} \notag\\
+C\underbrace{\|(\vn\eta)\wedge \vw \|_{\M_{t,x}^{\p_0,\q_0}}}_{(4)}+\underbrace{\|(\vn \eta) div(\vw)\|_{\M_{t,x}^{\p_0,\q_0}}}_{(5)}+\underbrace{\|\eta\vw\|_{\M_{t,x}^{\p_0,\q_0}}}_{(6)}+\underbrace{\|\vw\otimes \vu\cdot \vn \eta\|_{\M_{t,x}^{\p_0,\q_0}}}_{(7)}.\qquad\label{A1_Morrey}
\end{eqnarray}
Each term above is studied separately:
\begin{itemize}
\item For the first term of (\ref{A1_Morrey}), we note that since $\p_0\leq \q_0=\frac{5}{2-\alpha}$ and since $0<\alpha<\frac{1}{24}$, we have $\p_0\leq \q_0<3<6$, and thus by the support properties of the function $\eta$ as well as by the properties of Morrey spaces given in the Lemma \ref{lemma_locindi}, we  obtain
$$ \|(\partial_t\eta+\Delta \eta)(\vu+\vw)\|_{\M_{t,x}^{\p_0,\q_0}}\leq C \|\mathds{1}_{Q_{\mathbf{R}}}(\vu+\vw)\|_{\M_{t,x}^{\p_0,\q_0}}\leq C \|\mathds{1}_{Q_{\mathfrak{r}_2}}(\vu+\vw)\|_{\M_{t,x}^{6,6}}=C \|\mathds{1}_{Q_{\mathfrak{r}_2}}(\vu+\vw)\|_{L_{t,x}^{6}}<+\infty,$$
since we have the controls (\ref{Information_PremiereEtape}) and the relationships (\ref{Info_Radios2}).

\item  The terms (2) and (7) of (\ref{A1_Morrey}) can be treated in a similar manner. Indeed, since $0<\alpha<\frac{1}{24}$ we have $\p_0\leq \q_0<3$ and by the same arguments as above we write for (2):
\begin{eqnarray*}
\|\vu\otimes \vu\cdot \vn \eta\|_{\M_{t,x}^{\p_0,\q_0}}&\leq &C\|\mathds{1}_{Q_{\mathbf{R}}}\vu\otimes \vu \|_{\M_{t,x}^{\p_0,\q_0}}\leq C\|\mathds{1}_{Q_{\mathfrak{r}_1}}\vu\otimes \vu\|_{\M_{t,x}^{3,3}}\\
&\leq &C\|\mathds{1}_{Q_{\mathfrak{r}_1}}\vu\otimes \vu\|_{L_{t,x}^{3}}\leq C\|\mathds{1}_{Q_{\mathfrak{r}_1}}\vu\|_{L_{t,x}^{6}}\|\mathds{1}_{Q_{\mathfrak{r}_1}}\vu\|_{L_{t,x}^{6}}<+\infty,
\end{eqnarray*}
where we used the H\"older inequality in the last estimate as well as the controls (\ref{Information_PremiereEtape}) and the relationships (\ref{Info_Radios2}). The same ideas apply for (7).

\item For the term (3) of (\ref{A1_Morrey}), we recall that by the equation (\ref{Eq_Pression}) over the pressure we have the expression $p = \displaystyle{\sum_{i,j=1}^3}\frac{\partial_i\partial_j}{(-\Delta)}\left(u_i u_j\right)$.   We consider now two auxiliary functions $\widetilde{\phi}$ and $\widetilde{\psi}$ satisfying the same properties stated in (\ref{ProprieteLocalisation_IterationVU}) and such that
\begin{equation*}
\widetilde{\phi} \equiv 1\;\; \text{over}\; \; Q_{r_b}(t_0,x_0),\; \; supp(\widetilde{\phi})\subset Q_{r_b}(t_0,x_0) \quad \mbox{and}\quad\widetilde{\psi} \equiv 1\;\; \text{over}\; \; Q_{r_a}(t_0,x_0),\; \; supp(\widetilde{\psi})\subset Q_{ R_3}(t_0,x_0),
\end{equation*}
where $\mathbf{R}<r_b<r_a<R_3$.\\ 
Thus, by definition of the auxiliary function $\widetilde{\phi}$ we have the identity $\mathds{1}_{Q_{R_3}}=\widetilde{\phi} \mathds{1}_{Q_{R_3}}$ (recall the relationships \ref{Info_Radios2}). Thus the term $\widetilde{\phi} \vn p=\widetilde{\phi}  \displaystyle{\sum^3_{i,j= 1}}\frac{\vn}{(- \Delta )} \partial_i  \partial_j  (u_i u_j )$ can be rewritten in the following manner
\begin{equation}\label{Formule_Pression_Ladyz}
\begin{split}
\widetilde{\phi} \vn p=&\underbrace{\sum^3_{i,j= 1} \widetilde{\phi} \frac{\vn \partial_i  \partial_j }{(- \Delta )} (\widetilde{\psi} u_i u_j )}_{(a)}-\underbrace{\sum^3_{i,j= 1}  \frac{\widetilde{\phi} \vn\partial_i}{(- \Delta )} (\partial_j \widetilde{\psi}) u_i u_j}_{(b)}-\underbrace{\sum^3_{i,j= 1}  \frac{\widetilde{\phi} \vn\partial_j}{(- \Delta )} (\partial_i \widetilde{\psi}) u_i u_j}_{(c)}\\
&+ 2 \underbrace{\sum^3_{i,j= 1} \widetilde{\phi} \frac{\vn}{(-\Delta )}(\partial_i \partial_j\widetilde{\psi}) (u_i u_j)}_{(d)}+ \underbrace{\widetilde{\phi} \frac{\vn \big(  (\Delta \widetilde{\psi}) p \big) }{(-\Delta)}}_{(e)}
-2 \underbrace{\sum^3_{i= 1}\widetilde{\phi}\frac{\vn\big(\partial_i ( (\partial_i \widetilde{\psi}) p )\big)}{(-\Delta)}}_{(f)}
\end{split}
\end{equation}
and since $0<\alpha<\frac{1}{24}$ we have $\p_0\leq \q_0=\frac{5}{2-\alpha}<\frac{120}{47}<3$ and we only need to prove that each one of these terms belong to the Morrey space $\mathcal{M}^{\frac{6}{5}, \frac{120}{47}}_{t,x}(\mathbb{R}\times \R)$.
\begin{itemize}
\item[$*$] The term $(a)$ in (\ref{Formule_Pression_Ladyz}) is treated as follows: since the Riesz transforms are bounded in Morrey spaces we obtain
$$\left\|\widetilde{\phi}\frac{\vn  \partial_i \partial_j }{(- \Delta )}( \widetilde{\psi} u_i u_j )\right\|_{\mathcal{M}^{\frac{6}{5},\frac{120}{47}}_{t,x}}\leq C\left\|\frac{  \partial_i \partial_j \vn}{(- \Delta )}( \widetilde{\psi} u_i u_j )\right\|_{\mathcal{M}^{\frac{6}{5},\frac{120}{47}}_{t,x}}\leq C\left\|\vn( \widetilde{\psi} u_i u_j )\right\|_{\mathcal{M}^{\frac{6}{5}, \frac{120}{47}}_{t,x}},$$
now, for $1\leq k\leq 3$, using all the information available over $\vu$ (see (\ref{Information_PremiereEtape})), by Lemma \ref{lemma_locindi} and by the H\"older inequality in Morrey spaces, we have
$$\left\|(\partial_k\widetilde{\psi}) u_i u_j \right\|_{\mathcal{M}^{\frac{6}{5}, \frac{120}{47}}_{t,x}}\leq C\left\|\mathds{1}_{Q_{\bar R_2}}u_i u_j \right\|_{\mathcal{M}^{\frac{3}{2}, 30}_{t,x}}\leq C\|\mathds{1}_{Q_{\bar R_2}}u_i\|_{\mathcal{M}^{3, 60}_{t,x}}\|\mathds{1}_{Q_{\bar R_2}}u_j\|_{\mathcal{M}^{3, 60}_{t,x}}<+\infty,$$
since $\frac{6}{5}<\frac{3}{2}$, $\frac{120}{47}<30$ and $\frac{2}{3}=\frac{1}{3}+\frac{1}{3}$, $\frac{1}{30}=\frac{1}{60}+\frac{1}{60}$.
By the same arguments (recall the informations over $\vu$ given in (\ref{Information_PremiereEtape})) we have
\begin{eqnarray*}
\|\widetilde{\psi}(\partial_k u_i) u_j \|_{\mathcal{M}^{\frac{6}{5}, \frac{120}{47}}_{t,x}}&\leq &C\|\mathds{1}_{Q_{R_1}}\vn \otimes \vu\|_{\mathcal{M}^{2, \frac{120}{45}}_{t,x}}\|\mathds{1}_{Q_{R_3}} u_j\|_{\mathcal{M}^{3, 60}_{t,x}}\\
&\leq &C\|\mathds{1}_{Q_{R_1}}\vn \otimes \vu\|_{\mathcal{M}^{2, \tau_1}_{t,x}}\|\mathds{1}_{Q_{R_3}} u_j\|_{\mathcal{M}^{3, 60}_{t,x}}<+\infty,\\
 \|\widetilde{\psi} u_i(\partial_k u_j) \|_{\mathcal{M}^{\frac{6}{5}, \frac{120}{47}}_{t,x}}&\leq &C\|\mathds{1}_{Q_{R_3}}u_i\|_{\mathcal{M}^{3, 60}_{t,x}}\|\mathds{1}_{Q_{R_1}} \vn \otimes \vu\|_{\mathcal{M}^{2, \frac{120}{45}}_{t,x}}\\
&\leq& C\|\mathds{1}_{Q_{R_3}}u_i\|_{\mathcal{M}^{3, 60}_{t,x}}\|\mathds{1}_{Q_{R_1}} \vn \otimes \vu\|_{\mathcal{M}^{2, \tau_1}_{t,x}}<+\infty,
\end{eqnarray*}
since $\frac{5}{6}=\frac{1}{2}+\frac{1}{3}$ and $\frac{47}{120}=\frac{45}{120}+\frac{1}{60}$ and $\frac{120}{45}<\tau_1=\frac{30}{11}<\frac{20}{7}$. Thus we can deduce that we have the estimate 
$$\left\|\widetilde{\phi}\frac{\vn  \partial_i \partial_j }{(- \Delta )}( \widetilde{\psi} u_i u_j )\right\|_{\mathcal{M}^{\frac{6}{5}, \frac{120}{47}}_{t,x}}<+\infty.$$

\item[$*$] The terms $(b)$ and $(c)$ of (\ref{Formule_Pression_Ladyz}) can be treated in a similar manner and using the information available in (\ref{Information_PremiereEtape}) we have:
\begin{eqnarray*}
\left\|\frac{\widetilde{\phi} \vn\partial_i}{(- \Delta )} (\partial_j \widetilde{\psi}) u_i u_j\right\|_{\mathcal{M}^{\frac{6}{5}, \frac{120}{47}}_{t,x}}&\leq &C \left\|\frac{\vn\partial_i}{(- \Delta )} (\partial_j \widetilde{\psi}) u_i u_j\right\|_{\mathcal{M}^{\frac{6}{5}, \frac{120}{47}}_{t,x}}\leq C \left\|(\partial_j \widetilde{\psi}) u_i u_j\right\|_{\mathcal{M}^{\frac{6}{5}, \frac{120}{47}}_{t,x}}\\
&\leq &C\|\mathds{1}_{Q_{\bar R_2}}u_i u_j\|_{\mathcal{M}^{\frac{3}{2}, 30}_{t,x}}\leq C\|\mathds{1}_{Q_{\bar R_2}}u_i\|_{\mathcal{M}^{3, 60}_{t,x}}\|\mathds{1}_{Q_{\bar R_2}}u_j\|_{\mathcal{M}^{3, 60}_{t,x}}<+\infty.
\end{eqnarray*}
\item[$*$] The term $(d)$ is treated as follows. By Lemma \ref{lemma_locindi}, since $\frac{6}{5}<\frac{3}{2}$ and $\frac{120}{47}< \frac{15}{4}$, we have
$$\left\|\widetilde{\phi} \frac{\vn}{(-\Delta )}(\partial_i \partial_j \widetilde{\psi}) (u_i u_j)\right\|_{\mathcal{M}^{\frac{6}{5}, \frac{120}{47}}_{t,x}}\leq C\left\|\widetilde{\phi} \frac{\vn}{(-\Delta )}(\partial_i \partial_j \widetilde{\psi}) (u_i u_j)\right\|_{\mathcal{M}^{\frac{3}{2}, \frac{15}{4}}_{t,x}}.$$
Now, by the space inclusion $L^{\frac{3}{2}}_t L^{\infty}_{x}\subset \mathcal{M}^{\frac{3}{2}, \frac{15}{4}}_{t,x}$ we obtain
$$\left\|\widetilde{\phi} \frac{\vn}{(-\Delta )}(\partial_i \partial_j \widetilde{\psi}) (u_i u_j)\right\|_{\mathcal{M}^{\frac{3}{2}, \frac{15}{4}}_{t,x}}\leq C\left\|\widetilde{\phi} \frac{\vn}{(-\Delta )}(\partial_i \partial_j  \widetilde{\psi}) (u_i u_j)\right\|_{L^{\frac{3}{2}}_t L^{\infty}_{x}}$$
Following the same ideas displayed in formulas (\ref{Formula_intermediairevVV80})-(\ref{Formula_intermediairevVV81}), due to the support properties of the auxiliary functions we obtain
$$\left\|\widetilde{\phi} \frac{\vn}{(-\Delta )}(\partial_i \partial_j \widetilde{\psi}) (u_i u_j)\right\|_{L^{\frac{3}{2}}_t L^{\infty}_{x}}\leq C\|\mathds{1}_{Q_{\bar R_2}}u_i u_j\|_{L^{\frac{3}{2}}_{t,x}}\leq C\|\mathds{1}_{Q_{\bar R_2}}\vu\|_{\mathcal{M}^{3,30}_{t,x}}\|\mathds{1}_{Q_{\bar R_2}}\vu\|_{\mathcal{M}^{3,30}_{t,x}}<+\infty.$$

\item[$*$] The term $(e)$ of (\ref{Formule_Pression_Ladyz}) follows the same ideas as previous one, and we have
$$\left\|\widetilde{\phi} \frac{\vn \big(  (\Delta \widetilde{\psi}) p \big) }{(-\Delta)}\right\|_{\mathcal{M}^{\frac{6}{5}, \frac{120}{47}}_{t,x}}\leq C\left\|\widetilde{\phi} \frac{\vn \big(  (\Delta \widetilde{\psi}) p \big) }{(-\Delta)}\right\|_{L^{\frac{3}{2}}_t L^{\infty}_{x}}\leq C \|\mathds{1}_{Q_{R}}p\|_{L^{\frac{3}{2}}_{t,x}}<+\infty,$$
since we have by hypothesis that $\mathds{1}_{Q_{R}}p\in L^{\frac{3}{2}}_{t,x}(\mathbb{R}\times \R)$.

\item[$*$] The last term of (\ref{Formule_Pression_Ladyz}) is estimated in a very similar manner:
$$\left\|\widetilde{\phi}\frac{\vn\big(\partial_i ( (\partial_i \widetilde{\psi}) p )\big)}{(-\Delta)}\right\|_{\mathcal{M}^{\frac{6}{5}, \frac{120}{47}}_{t,x}}\leq C\left\|\widetilde{\phi}\frac{\vn\big(\partial_i ( (\partial_i \widetilde{\psi}) p )\big)}{(-\Delta)}\right\|_{L^{\frac{3}{2}}_t L^{\infty}_{x}}\leq C \|\mathds{1}_{Q_{R}}p\|_{L^{\frac{3}{2}}_{t,x}}<+\infty.$$
\end{itemize}
We have proven that all the terms of (\ref{Formule_Pression_Ladyz}) belong to the Morrey space $\M_{t,x}^{\p_0,\q_0}$ and thus, the term (3) of (\ref{A1_Morrey}) too.\\

\item The terms (4) and (6) of (\ref{A1_Morrey}) are very similar. Indeed, for (4), using the properties of the auxiliar function $\eta$ and with the Lemma \ref{lemma_locindi} we write (recall that $\p_0\leq \q_0<6$ and that we have the controls (\ref{Information_PremiereEtape}))
$$\|(\vn\eta)\wedge \vw \|_{\M_{t,x}^{\p_0,\q_0}}\leq C\|\mathds{1}_{Q_{\mathbf{R}}}\vw \|_{\M_{t,x}^{\p_0,\q_0}}\leq C\|\mathds{1}_{Q_{\mathfrak{r}_2}}\vw \|_{L_{t,x}^{6}}<+\infty.$$
For (6) we have by the same arguments:
$$\|\eta\vw\|_{\M_{t,x}^{\p_0,\q_0}}\leq C\|\mathds{1}_{Q_{\mathbf{R}}}\vw\|_{\M_{t,x}^{\p_0,\q_0}}\leq C\|\mathds{1}_{Q_{\mathfrak{r}_2}}\vw \|_{L_{t,x}^{6}}<+\infty.$$

\item For the term (5) of (\ref{A1_Morrey}), we need to study the quantity $\|(\vn \eta) div(\vw)\|_{\M_{t,x}^{\p_0,\q_0}}$, but by the Proposition \ref{Propo_divw} we know that $\mathds{1}_{Q_{R_3}}div(\vw) \in\M_{t,x}^{\frac{6}{5},\frac{60}{11}}(\mathbb{R}\times\R)$. Since $0<\alpha<\frac{1}{24}$, we have $\p_0\leq \frac{6}{5}\leq \q_0=\frac{5}{2-\alpha}<\frac{60}{11}$, and by the support properties of the function $\eta$ (recall (\ref{Info_Radios2}) and (\ref{Def_Eta})) we have (by Lemma \ref{lemma_locindi})
$$\|(\vn \eta) div(\vw)\|_{\M_{t,x}^{\p_0,\q_0}}\leq C\|\mathds{1}_{Q_{R_3}}div(\vw)\|_{\M_{t,x}^{\frac{6}{5},\frac{60}{11}}}<+\infty.$$
With all these estimates, we can conclude that $\vec{\mathcal{A}} \in \M_{t,x}^{\p_0,\q_0}(\mathbb{R}\times\R)$, and Lemma \ref{Lemma_EstimationA} is proven. \hfill $\blacksquare$\\
\end{itemize}
We study now the quantity $\vec{\mathcal{B}}_i$ defined in (\ref{B_Morrey}). Following (\ref{Objetivo_MorreyLadyz}), we shall obtain that  $\vec{\mathcal{B}}_i\in \M_{t,x}^{\p_0,\q_1}(\mathbb{R}\times\R)$ where $1\le \p_0\le \q_0$, with  $\frac{1}{\q_0}=\frac{2-\alpha}{5}$,$\frac{1}{\q_1}=\frac{1-\alpha}{5}$, for some $0<\alpha<1$. Since $0<\alpha<\frac{1}{24}$, we have $\q_1=\frac{5}{1-\alpha}<6$ and we thus write
$$\|\vec{\mathcal{B}}_i\|_{\M_{t,x}^{\p_0,\q_1}}= \|2(\partial_i \eta)(\vu+\vw)\|_{\M_{t,x}^{\p_0,\q_1}}\leq C\|\mathds{1}_{Q_{\mathfrak{r}_2}}(\vu+\vw)\|_{\M_{t,x}^{\p_0,\q_1}}\leq C \|\mathds{1}_{Q_{\mathfrak{r}_2}}(\vu+\vw)\|_{L_{t,x}^{6}}<+\infty,$$
where we used the support properties of the test function $\eta$, the Lemma \ref{lemma_locindi} and the controls (\ref{Information_PremiereEtape}).\\ 
We thus obtain that $\vec{\mathcal{B}_i} \in \M_{t,x}^{\p_0,\q_1}(\mathbb{R}\times\R)$.\\

For the term $\mathcal{C}$ given in (\ref{C_Morrey}) we have $\|\mathcal{C}\|_{\M_{t,x}^{\p_0,\q_1}}=\|\eta div(\vw)\|_{ \M_{t,x}^{\p_0,\q_1}}$. Since $1\leq \p_0\leq \frac{6}{5}$ and $\q_1=\frac{5}{1-\alpha}<\frac{60}{11}$ (since $0<\alpha<\frac{1}{24}$), by the support properties of the function $\eta$ and by Lemma \ref{lemma_locindi} we obtain
$$\|\eta div(\vw)\|_{ \M_{t,x}^{\p_0,\q_1}}\leq C\|\mathds{1}_{Q_{R_3}} div(\vw)\|_{ \M_{t,x}^{\frac{6}{5},\frac{60}{11}}}<+\infty.$$
With this estimate we obtain $\mathcal{C} \in \M_{t,x}^{\p_0,\q_1}(\mathbb{R}\times\R)$.\\

The term $\vec{\mathcal{D}}$ given in (\ref{D_Morrey}) can be treated just as the terms $\vec{\mathcal{B}}_i$ above. Indeed, using the controls (\ref{Information_PremiereEtape}) we write:
$$\|\vec{\mathcal{D}}\|_{\M_{t,x}^{\p_0,\q_1}}=\|\frac{1}{2}\eta(\vw+\vu)\|_{\M_{t,x}^{\p_0,\q_1}}\leq C\|\mathds{1}_{Q_{\mathfrak{r}_2}}(\vu+\vw)\|_{\M_{t,x}^{\p_0,\q_1}}\leq C \|\mathds{1}_{Q_{\mathfrak{r}_2}}(\vu+\vw)\|_{L_{t,x}^{6}}<+\infty.$$
We have $\vec{\mathcal{D}} \in \M_{t,x}^{\p_0,\q_1}(\mathbb{R}\times\R)$.\\

For the tensor $\mathbb{E}$ defined in (\ref{E_Morrey}), since $1\leq \p_0\leq \frac{6}{5}$ and $\q_1=\frac{5}{1-\alpha}<\frac{60}{11}$ we obtain, by Lemma \ref{lemma_locindi}:
\begin{eqnarray*} 
\|\mathbb{E}\|_{\M_{t,x}^{\p_0,\q_1}}&=&\|\eta (\vu\otimes \vu +\vw\otimes \vu)\|_{\M_{t,x}^{\p_0,\q_1}}\leq \|\eta \vu\otimes \vu\|_{\M_{t,x}^{\p_0,\q_1}} +\|\eta\vw\otimes \vu\|_{\M_{t,x}^{\p_0,\q_1}}\\
&\leq &C\|\eta \vu\otimes \vu\|_{\M_{t,x}^{\frac{6}{5},\frac{60}{11}}} +C\|\eta\vw\otimes \vu\|_{\M_{t,x}^{\frac{6}{5},\frac{60}{11}}}
\end{eqnarray*}
and by the H\"older inequalities in Morrey spaces (see Lemma \ref{lemma_Product}) with $\frac{5}{6}=\frac{1}{2}+\frac{1}{3}$ and $\frac{11}{60}=\frac{1}{6}+\frac{1}{60}$, we can write :
\begin{eqnarray*}
\|\mathbb{E}\|_{\M_{t,x}^{\p_0,\q_1}}&\leq &C\|\eta \vu\|_{\M_{t,x}^{2,6}}\|\eta \vu\|_{\M_{t,x}^{3,60}}+C\|\eta \vw\|_{\M_{t,x}^{2,6}}\|\eta \vu\|_{\M_{t,x}^{3,60}}\\
&\leq & C\|\mathds{1}_{Q_{\bar R_2}} \vu\|_{L_{t,x}^{6}}\|\mathds{1}_{Q_{\bar R_2}} \vu\|_{\M_{t,x}^{3,60}}+C\|\mathds{1}_{Q_{\bar R_2}} \vw\|_{L_{t,x}^{6}}\|\mathds{1}_{Q_{\bar R_2}} \vu\|_{\M_{t,x}^{3,60}}<+\infty.
\end{eqnarray*}
We thus have $\mathbb{E} \in \M_{t,x}^{\p_0,\q_1}(\mathbb{R}\times\R)$.\\

With all the previous computations we have proven all the information stated in (\ref{Objetivo_MorreyLadyz}), which applied in the integral representation formula (\ref{Duhamel_todoMicroPolar}) allows us, with Lemma \ref{Lemma_parabolicHolder}, to conclude that  $\vUU\in \dot{\mathcal{C}}^\alpha(\mathbb{R}\times \R)$ with $0<\alpha<\frac{1}{24}$, and since by (\ref{defU}) we have $\vUU=(\vu+\vw) $ over a small neighborhood of the point $(t_0, x_0)$, we deduce that $\vu$ and $\vw$ are also H\"older regular and this finishes the proof of Theorem \ref{Teo_HolderRegularity}. \hfill $\blacksquare$\\



\begin{thebibliography}{2}
\bibitem{Adams} 
D. R. \textsc{Adams},  J. \textsc{Xiao}. \emph{Morrey spaces in harmonic analysis}.  Ark. Mat. Vol. 50, N. 2, 201-230. (2012).
\bibitem{Bat00}
G. \textsc{Batchelor}. \emph{An Introduction to Fluid Dynamics}.  Cambridge University Press. (2000).
\bibitem{Brezis}
H. \textsc{Brezis}, P. \textsc{Mironescu}. \emph{Gagliardo-Nirenberg inequalities and non-inequalities: the full story}. Annales de l’Institut Henri Poincaré (C) Non Linear Analysis, Elsevier, 35(5):1355-1376. (2018).
\bibitem{LibroBrezis}
H. \textsc{Brezis}. \emph{Functional Analysis, Sobolev Spaces and Partial Differential Equations}. Universitext, Springer. (2011).
\bibitem{CKN} 
L. \textsc{Caffarelli}, R. \textsc{Kohn} \& L. \textsc{Nirenberg}. \emph{Partial regularity of suitable weak solutions of the Navier--Stokes equations}. Communications on Pure and Applied Mathematics, 35(6):771-831. (1982).
\bibitem{ChHe21}
D. \textsc{Chamorro},  J. \textsc{He}.   \emph{On the partial regularity theory for the MHD equations}. J. Math. Anal. Appl., 494(1):124449-124487. (2020).
\bibitem{ChLl21}
D. \textsc{Chamorro}, D. \textsc{Llerena}.   \emph{Interior e-regularity theory for the solutions of the magneto-micropolar equations with a perturbation term}.
J. Elliptic Parabol. Equ., 8:555–616. (2022). 
\bibitem{ChLl22}
D. \textsc{Chamorro}, D. \textsc{Llerena}. \emph{ A crypto-regularity result for the micropolar fluids equations}. Journal of Mathematical Analysis and Applications, 520(2). (2023).
\bibitem{Cru19}
F. \textsc{Cruz}. \emph{Global strong solutions for the incompressible micropolar fluids equations}. Arch. Math., 113: 201–212 (2019). 
\bibitem{Cru22}
F.W. \textsc{Cruz}, C.F. \textsc{Perusato}, M.A. \textsc{Rojas-Medar}, \& P.R. \textsc{Zingano}.  \emph {Large time behavior for MHD micropolar fluids in $\mathbb{R}^n$}. Journal of Differential Equations, 312:1-41. (2022).
\bibitem{Eri66}
A.C. \textsc{Eringen}, \emph{Theory of micropolar fluids}. J. Math. Mech. 16:1–18. (1966).
\bibitem{Gal97}
G.P. \textsc{Galdi}, S. \textsc{Rionero}. \emph{A note on the existence and uniqueness of solutions of the micropolar fluid equations}. Internat. J. Eng. Sci. 15:105–108. (1977).
\bibitem{Gal11}
G.P. \textsc{Galdi}. \emph{An Introduction to the Mathematical Theory of the Navier-Stokes Equations Steady-State Problems}. Springer Monographs in Mathematics. Springer. (2011).
\bibitem{Kukavica} 
I. \textsc{Kukavica}. \emph{On partial regularity for the Navier--Stokes equations}. Discrete and continuous dynamical systems, 21:717-728. (2008).
\bibitem{Lady68}
O. \textsc{Ladyzhenskaya}, V. \textsc{Solonnikov}, \& N. \textsc{Uraltseva}. \emph{Linear and quasilinear equations of parabolic type}. English translation : American Math. Society. (1968).
\bibitem{PGLR0} P.G. \textsc{Lemarié-Rieusset}. \emph{Recent Developments in the Navier-Stokes problem}. Chapman \& Hall/CRC. (2002).
\bibitem{PGLR1} P.G. \textsc{Lemarié-Rieusset}. \emph{The Navier-Stokes problem in the 21st century}. Chapman \& Hall/CRC. (2016).
\bibitem{Loa16} 
M. \textsc{Loayza},  M.A. \textsc{Rojas-Medar}, \emph{A weak-Lp Prodi-Serrin type regularity criterion for the micropolar fluid equations}. Journal of Mathematical Physics, 57(2). (2016). 
\bibitem{LoMelo}
J. \textsc{Lorenz}, W. G. \textsc{Melo} \& S. C. P. \textsc{de Souza}. \emph
{Regularity criteria for weak solutions of the Magneto-micropolar equations}.
Electronic Research Archive, 29(1):1625–1639. (2021).
\bibitem{Luka99}
G. \textsc{Łukaszewicz}. \emph{Micropolar Fluids, Theory and Applications}. Birkhäuser Boston.(1999).
\bibitem{Nic22}
C. J. \textsc{Niche}, C. F. \textsc{Perusato}. \emph{Sharp decay estimates and asymptotic behaviour for 3D magneto-micropolar fluids}. Zeitschrift Angew. Math. Phys., 73(2). (2022)
\bibitem{OLeary}
M. \textsc{O'Leary}. \emph{Conditions for the local boundedness of solutions of the Navier--Stokes system in three dimensions}. Comm. Partial Differential Equations, 28:617-636. (2003).
\bibitem{Pro59}
G. \textsc{Prodi}. \emph{Un teorema di unicità per le equazioni di Navier-Stokes. Annali di Matematica}. 48:173–182. (1959).
\bibitem{RagWU20}
M.A. \textsc{Ragusa}, F. \textsc{Wu}. \emph{A regularity criterion for three-dimensional micropolar fluid equations in Besov spaces of negative regular indices}. Anal.Math.Phys. 10(3). (2020).
\bibitem{RemTic21}
A. \textsc{Remond-Tiedrez}, I. \textsc{Tice}. \emph{Anisotropic Micropolar Fluids Subject to a Uniform Microtorque: The Unstable Case.} Commun. Math. Phys. 381:947–999. (2021). 
\bibitem{Robinson}
J. \textsc{Robinson}. \emph{An introduction to the classical theory of the Navier– Stokes equations.} Lecture notes, IMECC-Unicamp. (2010).
\bibitem{Scheffer}
V. \textsc{Scheffer}. \emph{Hausdorff measure and the Navier-Stokes equation}. Comm. Math. Phys., 55:97-112. (1977).
\bibitem{Serr62}
J. \textsc{Serrin}. \emph{On the interior regularity of weak solutions of the Navier--Stokes equations}. Arch. Rat. Mech. Anal., 9:187-195. (1962).
\bibitem{Serr63}
J. \textsc{Serrin}. \emph{The initial value problem for the Navier–Stokes equations.}
In Nonlinear Problems  (Rudolph E. Langer, ed.), 69–98. Madison: The University of Wisconsin press. (1963).
\bibitem{Yam05}
N. \textsc{Yamaguchi}. \emph{Existence of global solution to the micropolar fluid system in a bounded domain}. Math. Method Appl. Sci., 28:1507-1526. (2005).
\bibitem{Yua08}
B. \textsc{Yuan}. \emph{Regularity of weak solutions to magneto-micropolar fluid equations}. Math. Sci. Ser., 30:1469–1480. (2010).
\end{thebibliography}
\end{document}